\documentclass[11pt]{article}
\usepackage{color}
\usepackage[colorlinks=true,citecolor=blue]{hyperref}
\usepackage{amssymb,amsmath,amsfonts,amsthm,times,graphicx,subfig,amsthm}
\usepackage{mathtools}
\usepackage[english]{babel}
\usepackage[utf8]{inputenc}
\usepackage[margin=0.90in]{geometry}
\usepackage[margin=0.2in]{caption}
\usepackage{amsmath}
\numberwithin{equation}{section}
\usepackage{mathpazo,times}
\graphicspath{{figures/}}

\makeatletter
\renewcommand\section{\@startsection {section}{1}{\z@}%
  {-2ex \@plus -1ex \@minus -.2ex}%
  {1ex \@plus.1ex}%
  {\normalfont\bf\sffamily}}
\renewcommand\subsection{\@startsection{subsection}{2}{\z@}%
  {-1.75ex\@plus -0.4ex \@minus -.2ex}%
  {0.6ex \@plus .1ex}%
  {\normalfont\small\bf\sffamily}}
\renewcommand\subsubsection{\@startsection{subsubsection}{3}{\z@}%
  {-0.6ex\@plus -0.2ex \@minus -.2ex}%
  {0.4ex \@plus .1ex}%
  {\normalfont\normalsize\it}}
\renewcommand\paragraph{\@startsection{paragraph}{4}{\z@}%
  {0.2ex \@plus0.2ex \@minus0.1ex}{-0.5em}%
  {\normalfont\normalsize\bfseries}}
\def\ps@headings{%
  \let\@oddfoot\@empty
  \let\@evenfoot\@empty
  \def\@evenhead{\small\sffamily\thepage\hfil\slshape\leftmark}%
  \def\@oddhead{\small\sffamily{\slshape\rightmark}\hfil\thepage}%
      \let\@mkboth\markboth
    \def\chaptermark##1{\markboth{{\ifnum \c@secnumdepth >\m@ne
          \if@mainmatter \@chapapp\ \thechapter. \ \fi \fi ##1}}{}}%
    \def\sectionmark##1{\markright {{\ifnum \c@secnumdepth >\z@
          \thesection. \ \fi ##1}}}}
\makeatother

\advance\textwidth -8mm
\advance\textheight -10mm
\advance\hoffset 4mm
\advance\voffset 4mm
\advance\textheight 0mm
\advance\parskip 4pt
\newdimen\figwidth
\newtheorem{theorem}{Theorem}[section]
\newtheorem{definition}[theorem]{Definition}
\newtheorem{corollary}[theorem]{Corollary}
\newtheorem{lemma}[theorem]{Lemma}

\def\maketitle{\par\noindent{\LARGE\bf\sffamily\thetitle}\\[1.4ex]
{\large\theauthor}\\[0.6ex]
\textit{\thetextinfo}\\[0.2ex]
{\small\today}\par\vglue1.4\bigskipamount}
\def\title#1{\def\thetitle{#1}}
\def\author#1{\def\theauthor{#1}}
\def\textinfo#1{\def\thetextinfo{#1}}

\def\be{\begin{equation}}
\def\ee{\end{equation}}
\def\bse{\begin{subequations}}
\def\ese{\end{subequations}}

\definecolor{deeppurple}{rgb}{0.5, 0, 0.7}

\def\doibase{http://dx.doi.org/}
\def\reftitle#1{``#1''}
\def\gl{\mathrel{\mathpalette\overl@ss>}}
\def\txtfrac#1#2{{\textstyle\frac{#1}{#2}}}
\def\half{{\textstyle\frac12}}
\def\halfep{{\textstyle\frac\ep 2}}
\def\sech{\mathop{\rm sech}\nolimits}
\def\diag{\mathop{\rm diag}\nolimits}

\def\dn{\mathop{\rm dn}\nolimits}
\def\Natural{\mathbb{N}}
\def\Real{\mathbb{R}}
\def\R{\mathbb{R}}
\def\Complex{\mathbb{C}}
\def\I{\mathbb{I}}

\def\Z{\mathbb{Z}}
\def\i{\text{i}}
\def\halfi{{\textstyle\frac\i2}}
\def\ep{\epsilon}
\def\halfe{{\textstyle\frac\ep2}}
\def\halfN{{\textstyle\frac N2}}
\def\n{|\hspace{-0.5mm}|}

\def\Re{\mathop{\rm Re}\nolimits}
\def\Im{\mathop{\rm Im}\nolimits}

\def\tr{\mathop{\rm tr}\nolimits}

\def\d{\mathrm{d}}

\def\e{\mathop{\rm e}\nolimits}
\def\@#1{{\mathbf{#1}}}
\def\_#1{{\mathsf{#1}}}
\let\~=\tilde
\let\==\bar
\let\^=\hat
\let\<=\langle
\let\>=\rangle

\def\Leps{\mathop{\mathcal{L}^\epsilon}\nolimits}

\def\min{\mathop{\rm min}\nolimits}
\def\max{\mathop{\rm max}\nolimits}

\def\Lax{{\mathrm{Lax}}}
\def\Dir{{\mathrm{Dir}}}

\def\C{{\mathbb C}}

\def\R{{\mathbb R}}
\def\N{{\mathbb N}}

\def\Z{{\mathbb Z}}
\def\delt{\delta}
\def\D{\Delta_\epsilon}
\def\1{{\bf 1}}
\def\s{\sigma}
\def\g{\gamma}
\makeatother

\let\trueparagraph=\paragraph
\def\paragraph#1{\par\smallskip\trueparagraph{\rm\textbf{#1}}}
\allowdisplaybreaks

\usepackage{epsfig,cancel,amsmath,amsthm,amssymb}
\usepackage{tikz}
\usetikzlibrary{decorations.fractals}
\usepackage{wrapfig}
\usepackage{graphicx,framed}
\usetikzlibrary{decorations.markings}

\begin{document}
\pagestyle{plain}

{\title{On the spectrum of the periodic focusing Zakharov-Shabat operator}
\author{Gino Biondini$^1$, Jeffrey Oregero$^1$ and Alexander Tovbis$^2$}
\textinfo{{\small1:} Department of Mathematics, State University of New York at Buffalo, Buffalo, New York 14260\\
{\small2:} Department of Mathematics, University of Central Florida, Orlando, Florida 32816}
\maketitle}

\begin{abstract}
\noindent
The spectrum of the focusing Zakharov-Shabat operator on the circle is studied,
and its explicit dependence on the presence of a semiclassical parameter is also considered.
Several new results are obtained. 
In particular: 
(i) it is proved that the resolvent set is comprised of two connected components, 
(ii) new bounds on the location of the Floquet and Dirichlet spectra are obtained, 
some of which depend explicitly on the value of the semiclassical parameter,
(iii) it is proved that the spectrum localizes to a ``cross" in the spectral plane in the semiclassical limit.
The results are illustrated by discussing several examples in which the spectrum is computed analytically or numerically.
\end{abstract}


\section{Introduction} 
\label{s:intro}

In this work we investigate the spectrum of the non-self-adjoint Zakharov-Shabat (ZS) scattering problem with a periodic potential. 
The ZS scattering problem is given by the first-order coupled system of ordinary differential equations (ODEs) \cite{APT2004, FT1987, NMPZ1984, ZS1972}:
\be
\epsilon \@v' = (-\i z\sigma_3 + Q(x))\@v\,,
\label{e:zs}
\ee
where $\@v(x;z,\epsilon) = (v_1, v_2)^{T} $ (the superscript $T$ denoting matrix transpose),
$\sigma_{3} = \diag (1,-1)$ is the third Pauli matrix, $Q(x)$ is the matrix-valued function
\be
Q(x) = \begin{pmatrix} 0 & q(x) \\ -\overline{q(x)} & 0 \end{pmatrix}\,,
\label{e:Q}
\ee
$z \in \Complex$ is the spectral parameter, 
prime denotes differentiation with respect to the independent variable (here $x$),
overbar denotes complex conjugation,
and $0 < \epsilon \leq 1$ is the semiclassical parameter.
Unless stated otherwise, throughout this work, the ``potential'' 
$q: \R \to \C$ is a complex-valued function with minimal period~$L$,
i.e.,
\be
q(x+L) = q(x)\,,\quad \forall x\in\Real\,.
\label{e:Qperiodic}
\ee
Unless stated otherwise we assume $q \in L^{\infty}(\Real)$, 
the space of essentially bounded Lebesgue measurable functions with essential supremum norm.
We say a function $\@v(x;z,\epsilon)$ is a solution of equation~\eqref{e:zs}, 
if it is locally absolutely continuous,
i.e.,
if $\@v\in AC_{\rm loc}(\R)$,
and the equality~\eqref{e:zs} holds almost everywhere
(cf.~\ref{s:notations}).

The main motivation for studying~\eqref{e:zs} derives from its role in  the analysis of the focusing nonlinear Schr\"odinger (NLS) equation on the circle,
which is given by
\be
\i\ep \partial_t\psi + \ep^{2}\partial_x^2\psi + 2|\psi|^{2}\psi = 0\,,
\label{e:nls}
\ee
where $\psi : \Real \times \Real_+ \to \Complex$ is the slowly-varying complex envelope of a quasi-monochromatic, 
weakly dispersive nonlinear wave packet, 
and the physical meaning of the variables depends on the context. 
(E.g., in nonlinear fiber optics, $t$ represents propagation distance while $x$ is a retarded time.)
In general, the parameter $\epsilon$ quantifies the relative stength of dispersion compared to nonlinearity. 
(In the quantum-mechanical setting, $\epsilon$ is also proportional to Planck's constant $\hbar$.)

Specifically, 
it was shown by Zakharov and Shabat in 1972 that~\eqref{e:zs} makes up the first half of the Lax pair for the focusing NLS equation~\cite{ZS1972}. 
This observation is the key to solving the initial value problem for~\eqref{e:nls} by means of the inverse scattering method.
The solution $\psi(x,t;\ep)$ of~\eqref{e:nls} with initial data $\psi(x,0;\ep):=q(x)$ is constructed by computing suitable scattering data generated by the potential $q(x)$ in~\eqref{e:zs}. 
Time evolution according to~\eqref{e:nls} corresponds to an isospectral deformation of the potential in~\eqref{e:zs},
and the time evolution of the scattering data can be computed trivially in certain cases.
This allows one to obtain the solution $\psi(x,t;\ep)$ of~\eqref{e:nls} by solving an inverse scattering problem, 
i.e., 
by reconstructing the potential of the ZS scattering problem from the knowledge of the time-evolved
scattering data.
As a result, 
analysis of ZS systems such as~\eqref{e:zs} has become an active area of research 
(e.g., see~\cite{BiondiniOregero, BL2018, bronski, djakovmityagin, fujiiewittsten2018, kapmity, kappeler, KS2002, KS2003, mclaughlinli, mclaughlinoverman, miller2001, SatsumaYajima,tkachenko,tovbisvenakides}) 
and a natural starting point in the study of solutions to the focusing NLS equation.

Importantly,~\eqref{e:zs} and~\eqref{e:nls} depend on the semiclassical parameter $\ep$. 
Letting $\ep \downarrow 0$ in~\eqref{e:nls} is referred to as the ``semiclassical limit'' since it allows one to establish a connection
between quantum and classical mechanics~\cite{Messiah}.
Similar phenomena occur in physical applications when dispersive effects are weak compared to nonlinear ones. 
These situations, referred to as small dispersion limits, can produce a wide variety of physical effects such as supercontinuum generation, dispersive shocks and wave turbulence, to name a few 
(e.g., see \cite{az2014, DudleyTaylor, elhoefer, onoratoosborne, randoux2014, solli, Whitham, ZabuskyKruskal, zakharov2009} and references therein).
Semiclassical, or small dispersion, limits are also of mathematical interest. Solutions to equations such as focusing NLS have rich structure which becomes more evident as the semiclassical parameter tends to zero.
As a result, 
much effort has been devoted to the analysis of integrable nonlinear evolution equations in the semiclassical limit
(e.g., see \cite{bertolatovbis, DVZ, deiftmclaughlin, pre2017deng, elkhamistovbis, jenkins2013,  KMM2003, LaxLevermore, lyngmiller, millerkamvissis1998, trillo2018naturephoton, TVZ2004, tvz2006} and references therein).

Since~\eqref{e:zs} is a non-self-adjoint eigenvalue problem, 
its spectrum is in general quite complicated.
Much of the research in this area has been devoted to studying~\eqref{e:zs} with zero-background potentials, 
or more precisely, $q\in L^1(\Real)$, i.e., Lebesgue integrable.
In 1974, 
Satsuma and Yajima proved if $q(x) = A\sech x$ with $A \in \Real$, then~\eqref{e:zs}  has purely imaginary point spectrum \cite{SatsumaYajima}. 
Moreover, 
when $A\in\N$ there are precisely $2A$ eigenvalues located at
$z = \pm\i(n-1/2)$, $n=1, \dots , A$. 
More recently, 
one of the authors and S. Venakides extended these results to a one-parameter family of potentials of the form $q(x;\ep)=\sech x\,\exp(\i S(x)/\ep)$, 
where $S'(x) = \kappa \tanh x$ with $\kappa\in\R$~\cite{tovbisvenakides}.
More generally, 
Klaus and Shaw proved that all real, 
piecewise-smooth, 
single-lobe potentials with zero background have purely imaginary point spectrum~\cite{KS2002, KS2003}. 
Here, 
the term single-lobe means $q(x)$ is nondecreasing when $x<0$, 
and is nonincreasing when $x>0$.
Recently, 
one of the authors and X.-D.\ Luo generalized the Klaus-Shaw result to real, single-lobe potentials on non-zero background, 
that is, 
$q(x)\to q_{0}$ as $|x| \to \infty$ and $q(x) > q_{0}$ for all $x\in \Real$ \cite{BL2018}. 
Considerable work has also been devoted to potentials with rapid phase variations, 
namely, 
$q(x;\ep)=A(x) \exp(\i S(x)/\epsilon)$, 
where $A(x)$ and $S(x)$ are both real.
Careful numerical experiments by Bronski showed that, in the semiclassical limit, the point spectrum can accumulate on a ``Y-shaped'' set of curves in the spectral plane, 
and the number of eigenvalues scales like $O(1/\epsilon)$ as $\ep \downarrow 0$~\cite{bronski}.
Also, 
bounds on the point spectrum were derived by Deift, 
Venakides, 
and Zhou under the assumptions $A(x) \to 0 $, 
and $S'(x) \to 0$ as $|x| \to \infty$~\cite{bronski}.
(For details of the proof see the work by DiFranco and Miller~\cite{difranco}.)
Later, 
Miller put forth a formal WKB based asymptotic analysis~\cite{miller2001}.
Moreover,
one of the authors and S. Venakides used the idea of the semiclassical limit of the scattering transform to derive spectral information of various potentials~\cite{TV2012}.
In particular, 
the Y-shaped spectral curve from~\cite{bronski} was confirmed there. 
The semiclassical limit of solutions to the focusing NLS equation generated by zero-background potentials was then studied in~\cite{bertolatovbis, elkhamistovbis, jenkins2013, KMM2003, lyngmiller, TVZ2004}. 
Even so, 
there are still many important open questions.

Many studies have been devoted to analyzing the spectrum of~\eqref{e:zs} with potentials in $L^{1}(\Real)$.
Many works have also been devoted to studying the spectrum of Hill operators with complex potentials
(e.g. see~\cite{djakovmityagin,Gesztesy_Hill1996,rofebek,rofebek2,serov,shin_2004, tkachenkannals} and references therein).
However, 
less is known for the non-self-adjoint ZS system~\eqref{e:zs} with a periodic potential.
Some localization results corresponding to (anti)periodic eigenvalues of \eqref{e:zs}
were obtained in~\cite{djakovmityagin,mclaughlinli}, and ``gap estimates'' in weighted Sobolev spaces were obtained in~\cite{djakovmityagin,kapmity,kappeler}.
The spectrum was also rigorously studied in~\cite{gesztesyweikard_acta1998,gesztesyweikard_bams1998,mclaughlinoverman,tkachenko}, where, 
among other results, 
asymptotic statements as $z\to\infty$ were obtained.
Moreover, 
for real-analytic periodic potentials, 
it was shown in~\cite{fujiiewittsten2018} that the set of periodic eigenvalues was discrete and clustered on the real and imaginary axes of the spectral variable in the 
limit $\ep\downarrow 0$.

The semiclassical limit of \eqref{e:zs} with real periodic potentials
was recently studied by two of the authors using a formal WKB approach, 
and numerical simulations~\cite{BiondiniOregero}.
This work is partly motivated by that study.
The main results of the present work are several statements about the spectrum of~\eqref{e:zs}, 
both in general and in the semiclassical limit. 
Specifically, 
the work is organized as follows:
We begin in Section~\ref{s:prelim} by briefly recalling various preliminary notions concerning the spectral problem~\eqref{e:zs}.
In Section~\ref{s:general} we present various general properties of the spectrum, 
several of which are new to the best of our knowledge.
For example, 
we show that the resolvent set of~\eqref{e:zs} with a periodic potential is comprised of two connected components,
in contrast with the case of Hill's equation, 
whose resolvent set is path connected.
In Section~\ref{s:semiclassical}
we show that, 
under fairly general assumptions on the potential, 
the spectrum localizes to a subset of the real and imaginary axes of the spectral variable in the limit $\epsilon \downarrow 0$,
thus putting the asymptotic and numerical results of~\cite{BiondiniOregero} into a rigorous mathematical setting.
In Section~\ref{s:real}
we discuss some further properties acquired by the spectrum when the potential is real and/or symmetric.
The proofs of all results discussed in sections~\ref{s:general}, 
\ref{s:semiclassical} and \ref{s:real} are given in Section~\ref{s:proofs}.
In Section~\ref{s:examples}, 
specific examples are analyzed exactly and via careful numerical simulations, 
and we comment on the similarities between the periodic problem and the infinite line problem when $0<\epsilon \ll 1$. 
Various details about the numerical simulations and some of the calculations are relegated to the appendices.

\section{Preliminaries}
\label{s:prelim}
Equation~\eqref{e:zs} is equivalent to the eigenvalue problem
\be
\Leps\@v = z\@v\,,
\label{e:op}
\ee
where
\be
\Leps := \i\sigma_{3}(\ep \partial_{x} - Q)\,,
\label{e:operator}
\ee
is a one-dimensional Dirac operator acting in $L^2(\Real,\C^2)$ with dense domain $H^{1}(\Real,\C^2)$,
and $Q$ is given by~\eqref{e:Q}.
Here $H^1$ is the Sobolev space of square integrable functions with square integrable first derivative.
(For a brief discussion on notation and function spaces see~\ref{s:notations}.)

Regarding the inverse spectral method for the focusing NLS equation on the circle,
one is concerned with the Lax spectrum of the operator $\Leps$, 
namely,
the set of $z\in\C$ for which a non-trivial solution of~\eqref{e:zs} exists which is bounded for all $x \in \Real$.  
Specifically,
the \textit{Lax spectrum}, or simply spectrum, is defined as:
\be
\Sigma_{\mathrm{Lax}} := 
\{z \in \Complex : 
\exists~\@v \not\equiv 0\in \textrm{AC}_\mathrm{loc}(\Real) ~
~{\rm s.t.}~ ~ \Leps \@v = z\@v
~ ~\&~ ~ \text{sup}_{x\in\Real} \|\@v(x;z,\ep)\| < \infty\}\,.
\label{e:laxspec}
\ee
It is well known that $\Sigma_{\Lax}$
is purely continuous, 
that is,
essential without any eigenvalues and empty residual spectrum~\cite{danilov,gesztesyweikard_acta1998,hislop,rofebek,rofebek2}.
Likewise, 
consider the operator~\eqref{e:operator},
now acting in $L^{2}([0,L],\C^2)$ with dense domain $H^1([0,L],\C^2)$. 
Then, 
for each $\nu\in\Real$, 
the associated \textit{Floquet spectrum} is defined as:
\be
\Sigma_{\nu} := \{z \in \Complex : \exists~\@v \not\equiv 0 \in H^1([0,L],\C^2)~
~{\rm s.t.}~ ~ \Leps \@v = z\@v ~ ~\&~ ~ \@v(L;z,\ep) = \e^{\i\nu L}\@v(0;z,\ep)\}\,.
\label{e:floquetspec}
\ee
Importantly,
$\nu=2n\pi/L$ corresponds to the periodic spectrum,
and $\nu=(2n-1)\pi/L$ corresponds to the antiperiodic spectrum,
where $n\in\Z$. 
Any $z\in\Sigma_{\nu}$ will be referred to as an eigenvalue of~\eqref{e:op}, 
and the corresponding $\@v(x;z,\ep)$ as a Floquet eigenfunction. 
Clearly the eigenvalues depend on $\ep$,
i.e.,
$z=z(\ep)$.

Basic properties of the Lax spectrum that follow from the theory of linear homogeneous ODEs with periodic coefficients 
are reviewed next to introduce some relevant concepts and to set the notation.
\begin{theorem}
\label{t:floquet}
(Floquet,~\cite{Floquet_int,Floquet})
Consider the system of linear homogeneous ODEs given by
\be
\@y' = A(x)\@y\,,
\label{e:yprime} 
\ee
where $A\in L^1_{{\rm loc}}(\R)$ is a  
$n \times n$ matrix-valued function such that $A(x+L) = A(x)$. 
Then any fundamental matrix solution $Y(x)$ of \eqref{e:yprime} can be written in the Floquet normal form
\be
Y(x) = \Psi(x)\e^{Rx}\,,
\label{e:canonical}
\ee
where $\Psi(x+L) = \Psi(x)$, 
$\Psi$ is nonsingular, 
and $R$ is a constant matrix. 
\end{theorem} 
Without loss of generality one can take $R$ to be in Jordan normal form. 
Since $\Psi(x)$ is 
locally absolutely continuous and periodic, 
the behavior of solutions as $x \to \pm \infty$ is determined by the eigenvalues, 
called Floquet exponents, 
of the matrix $R$.
In particular:
(i) If the Floquet exponent has non-zero real part then the solution grows exponentially as $x \to  \infty$, 
or as $x \to -\infty$; 
(ii) If the Floquet exponent has zero real part, but $R$ has non-trivial Jordan blocks then the solution is algebraically growing;
(iii) Otherwise the Floquet exponent is purely imaginary and the solution remains bounded for all $ x \in \Real $. 
 
By Theorem~\ref{t:floquet} all Floquet eigenfunctions of the ZS system~\eqref{e:zs} have the form
\be
\@v(x;z,\ep) = \e^{\i\nu x} \@w(x;z,\ep)\,,
\label{e:bounded}
\ee
where 
$\@w(x+L;z,\ep) = \@w(x;z,\ep)$,
and $\nu=\nu(\ep)\in\Real$ is the quasi-momentum.
More generally one has the so called normal solutions, 
that is, 
solutions of \eqref{e:zs} such that 
\be
\@v(x+L;z,\ep) = \sigma\,\@v(x;z,\ep)\,,
\label{e:bdsol}
\ee 
and $\sigma=\sigma(\ep)$ is the Floquet multiplier. 
Thus, 
a solution of~\eqref{e:zs} is bounded for all $x\in\R$ if and only if $|\sigma|=1$ in which case one has 
$\sigma = \e^{\i\nu L}$.
Moreover, 
the Floquet multipliers are the eigenvalues of a monodromy matrix, 
which is defined as
\be
Y(x+L;z,\ep)=Y(x;z,\ep)M(z;\ep)\,, 
\label{e:monodromy}
\ee
where $Y(x;z,\ep)$ is any fundamental matrix solution of~\eqref{e:zs}.
Let $\Phi:=\Phi(x;z,\ep)$ be the principal matrix solution of~\eqref{e:zs},
that is,
the solution of~\eqref{e:zs} normalized so that $\Phi(0;z,\ep)\equiv\mathbb{I}$,
where $\I$ is the $2\times 2$ identity matrix.
Since all monodromy matrices are similar, 
for the remainder of this work we fix
\be
M(z;\ep) = \Phi(L;z,\ep)\,.
\label{e:monodromy2}
\ee
Note also that the well known symmetries of solutions of~\eqref{e:zs} imply
\be
\overline{M(\overline{z};\ep)} = \sigma_{2} M(z;\ep) \sigma_2\,,
\label{e:msymmetry1}
\ee
where $\sigma_2$ is the second Pauli matrix (cf.~\ref{s:notations}).
In turn,~\eqref{e:msymmetry1} implies 
$M(z;\ep)$ can be written as
\be
M(z;\ep) = \begin{pmatrix} c(z;\ep) & -\overline{s(\overline{z};\ep)} \\[0.2ex] s(z;\ep) & \overline{c(\overline{z};\ep)} \end{pmatrix}\,,
\label{e:Mrepresentation}
\ee
where, by~\eqref{e:monodromy2}, $c(z;\epsilon)$ and $s(z;\epsilon)$ are the components at $x=L$ of the solution $\@v(x;z,\epsilon)$ of~\eqref{e:zs} equaling $(1,0)^T$ at $x=0$.
Also, 
since~\eqref{e:zs} is traceless, 
it follows from Abel's formula that $\det M(z;\ep) \equiv 1$.  
Hence the eigenvalues of $M(z;\ep)$ are given by roots of the quadratic equation
%
$\sigma^2 - (2 \D)\,\sigma + 1 = 0$,
%
where we introduced the Floquet discriminant
\be
\D:=\D(z)=\half\tr M(z;\ep)\,,
\label{e:DeltaM}
\ee
and ``$\tr$'' is the matrix trace.
Thus~\eqref{e:zs} has bounded solutions if and only if the following two conditions are simultaneously satisfied:
\bse
\label{e:deltaconditions}
\begin{gather}
\Im \D(z) = 0\,,
\label{e:ImDelta=0}
\\
-1\leq \Re \D(z) \leq 1\,.
\label{e:ReDelta_in[-1,1]}
\end{gather}
\ese
(Here ``$\Re$'' and ``$\Im$'' denote the real and imaginary components of a complex function, respectively.)
Thus, 
one gets the following equivalent representation of the Lax spectrum, namely,
\be
\Sigma_{\Lax} = \{ z \in \Complex : \D(z) \in [-1,1]\}.
\label{e:laxspec2}
\ee
Note~\eqref{e:Mrepresentation} and~\eqref{e:DeltaM}
immediately imply $\D$ is real-valued along the real $z$-axis.
Further,
recall the following:
\begin{lemma}(Entire,~\cite{MaAblowitz, mclaughlinoverman})
\label{l:entire}
The Floquet discriminant $\D(z)$ is entire.
\end{lemma}
Hence,
$\D$ satisfies the Schwarz reflection principle,
i.e.,
$\D(\overline{z})=\overline{\D(z)}$.

Additionally, 
the Floquet spectrum~\eqref{e:floquetspec} has the equivalent representation
(see~\cite{MW1966} for further details)
\be
\Sigma_{\nu}=\{z \in \Complex : \D(z) = \cos(\nu L)\,,\,\,\nu\in\R\}\,.
\label{e:bloch}
\ee
Clearly,
for $\nu\in\R$ one has 
$\Sigma_{\nu + 2\pi/L} = \Sigma_{\nu}$.
Importantly, 
the above arguments show that 
\be
\Sigma_{\mathrm{Lax}} = \cup_{\nu \in [0,2\pi/L)} \Sigma_{\nu}\,,
\label{e:specunion}
\ee
where the eigenfunctions of~\eqref{e:floquetspec}  
extend naturally to $\Real$.
That is,       
the Lax spectrum is the union of all Floquet spectra.
(Note however that the terminology used in the literature varies somewhat.)

Recall that $\Leps$ is isospectral with respect to time deformations of the potential that obey the NLS equation~\eqref{e:nls}. 
Thus, 
$\Sigma_{\nu}$ and $\Sigma_{\mathrm{Lax}}$,
and therefore also $\D$,
are all conserved with respect to the flow of the NLS equation~\eqref{e:nls}. 

Also recall that for the defocusing NLS equation
[i.e., \eqref{e:nls} with a negative sign in front of the nonlinear term],
the corresponding Zakharov-Shabat spectral problem 
[i.e., \eqref{e:zs} with the entry $-\overline{q(x)}$ in $Q(x)$ replaced with $\overline{q(x)}$]  
is self-adjoint, 
and therefore the spectrum is confined to the real $z$-axis.
In that case, 
for periodic potentials the real $z$-axis decomposes into a (possibly infinite) number of spectral bands~\cite{Floquet_int, Eastham, MW1966}. 
Thus, 
in the self-adjoint case, 
knowing the periodic/antiperiodic spectrum and their geometric multiplicities is enough to completely characterize the Lax spectrum.

In contrast, 
the focusing NLS equation corresponds to the non-self-adjoint case,
for which there is no restriction on the location of the spectral bands in the complex plane (apart from constraints such as the Schwarz symmetry).  
This greatly complicates the analysis;
and implies the entire Floquet spectrum is relevant.
Nonetheless, 
one can still introduce the concept of spectral bands and gaps like with self-adjoint problems, 
as discussed in section~\ref{s:general}.

Importantly,
the inverse spectral theory for~\eqref{e:zs} also involves the Dirichlet spectrum.
The \textit{Dirichlet spectrum} is defined as:
\be
\Sigma_{\Dir}(x_o)
:=\{\zeta \in \C : \exists~\@v \not\equiv 0 \in H^1([x_o,x_o+L],\C^2)~ ~\text{s.t.}~ ~ \Leps\@v = \zeta\@v ~ ~\&~ ~ 
\@v\in {\rm BC}_{{\rm Dir}}(x_o)
\}\,,
\label{e:dirichlet}
\ee
where ``BC$_{{\rm Dir}}(x_o)$'' are Dirichlet boundary conditions (BCs) with base point $x_o$,
i.e,
\bse
\label{e:Dirbcs}
\begin{gather}
v_1(x_o;\zeta,\ep) + v_2(x_o;\zeta,\ep) = 0\,,
\\
v_1(x_o+L;\zeta,\ep) + v_2(x_o+L;\zeta,\ep) = 0\,.
\end{gather}
\ese
Any value $\zeta\in\Sigma_{\Dir}(x_o)$ will be referred to as a Dirichlet eigenvalue of~\eqref{e:op}.
Similarly to the Floquet spectrum,
one can define $\Sigma_{\Dir}(x_o)$ as the zero set of an analytic function.
Consider the following modified fundamental matrix solution of~\eqref{e:zs}:
\be
\widetilde{\Phi}(x;z,\ep) = \Phi(x;z,\ep)\,C\,,
\qquad
C = \frac{1}{\sqrt{2}} \begin{pmatrix*}[r] -1 & \i \\ 1 & \i \end{pmatrix*}\,.
\label{e:unitary}
\ee
The monodromy matrix $\widetilde{M}^\ep(z)$ associated with $\widetilde{\Phi}$ 
is $\widetilde{M}^\ep(z)=C^{-1} M(z;\ep)C$.
Then it follows easily that 
(see~\ref{e:defndirichlet} for details)
\be
\Sigma_{\Dir}(0) = \{\zeta \in \C : \widetilde{M}^{\ep}_{21}(\zeta) = 0\}\,,
\label{e:dirzeros}
\ee
where the subscript ``21'' denotes the second row first column entry of the corresponding matrix.
Unlike the Floquet spectrum, however,
the Dirichlet spectrum is dependent on the base point $x_o$. 
Similarly, 
$\Sigma_{\Dir}(x_o)$ is not conserved by the flow of the NLS equation~\eqref{e:nls}.
Indeed,
the Dirichlet spectrum provides angle information in the ``action--angle'' formalism of an integrable system.
Moreover,
Floquet and Dirichlet spectra together comprise the set of scattering data from which one can reconstruct the potential~\eqref{e:Qperiodic}.
Namely, 
one has the following:
\begin{theorem}
(Trace formulae,~\cite{mclaughlinoverman})\,\ 
Let $\{z_n\}_{n\in\Z}$ be the sequence of periodic and antiperiodic eigenvalues,
$\{\zeta_n(x_o)\}_{n\in\Z}$ the sequences of Dirichlet eigenvalues, 
and $\{\widehat{\zeta}_n(x_o)\}_{n\in\Z}$ the sequences of auxiliary Dirichlet eigenvalues, 
defined respectively as follows:
\be
\D^2(z_n)-1=0\,, \qquad
\hat{M}^{\ep}_{21}(\zeta_n;x_o)=0\,,\qquad
\check{M}^{\ep}_{21}(\widehat{\zeta}_n;x_o)=0\,,
\ee
where $\hat{M}^{\ep}(z;x_o)$ and $\check{M}^{\ep}(z;x_o)$ 
are respectively the modified monodromy matrices $\widetilde{M}^{\ep}(z)$ associated with the translated potentials 
$q(x_o+x)$ and ${\rm i} q(x_o+x)$.
Then
\bse
\label{e:trace}
\begin{gather}
q(x) - \overline{q(x)} = 2 \sum_{j\in\mathbb{Z}} (z_{2j} + z_{2j+1} - 2\zeta_{j}(x))\,,
\label{e:trace1}
\\
q(x) + \overline{q(x)} = -2{\rm i} \sum_{j\in\mathbb{Z}} (z_{2j} + z_{2j+1} - 2\widehat{\zeta}_{j}(x))\,.
\label{e:trace2}
\end{gather}
\ese
\end{theorem}
The numbering of the eigenvalues in~\eqref{e:trace} is such that, 
for $|j|$ sufficiently large, 
take $z_{2j+1}=\overline{z}_{2j}$ and the $\zeta_j$ and $\widehat{\zeta}_j$ are the Dirichlet eigenvalues associated with $z_{2j}$ and $z_{2j+1}$
(cf.\ Theorem~\ref{theo-spines} and~\cite{djakovmityagin,kappeler}).

\section{General properties of the spectrum}
\label{s:general}

In this section we discuss some general properties of the Lax spectrum.  
Some of the results below were known, 
but here are proven for a broader class of potentials.
Other results are new to the best of our knowledge.
In this section we assume $\epsilon>0$ is fixed.
Proofs of all the results in this section are given in Section~\ref{s:mainproofs}.

Owing to~\eqref{e:ImDelta=0},
the Lax spectrum~\eqref{e:laxspec} is located along the contour lines 
\be
\Gamma:=\{z\in\Complex : \Im \D(z) = 0 \}\,.
\label{e:contours}
\ee
Moreover, 
$\Gamma$ is the union of an at most countable set of regular analytic curves $\Gamma_n$ in the complex $z$-plane~\cite{gesztesyweikard_acta1998}, 
each starting from infinity and ending at infinity:
\be
\Gamma=\cup_{n \in \Z} \Gamma_n\,.
\ee
(The precise details of the map $n\mapsto\Gamma_n$ are not important for the present purposes.)
\begin{definition}
\label{s:bandsandgaps} 
A spectral band is a maximally connected regular analytic arc along~$\Gamma_n$ determined by~\eqref{e:ReDelta_in[-1,1]}.
Each finite portion of $\Gamma_n$ where $|\Re \D| > 1$ and which is delimited by a band on either side is called a spectral gap. 
\end{definition}

\begin{figure}[t!]	
\vglue-4ex
\centerline{\includegraphics[width=8.5cm]{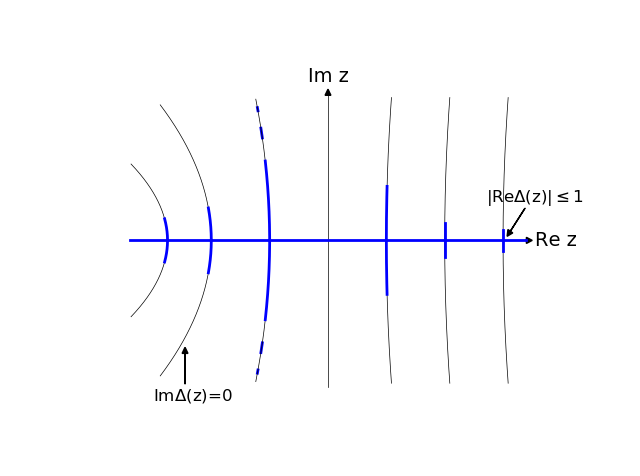}
\hspace{-8.5mm}
\includegraphics[width=8.5cm]{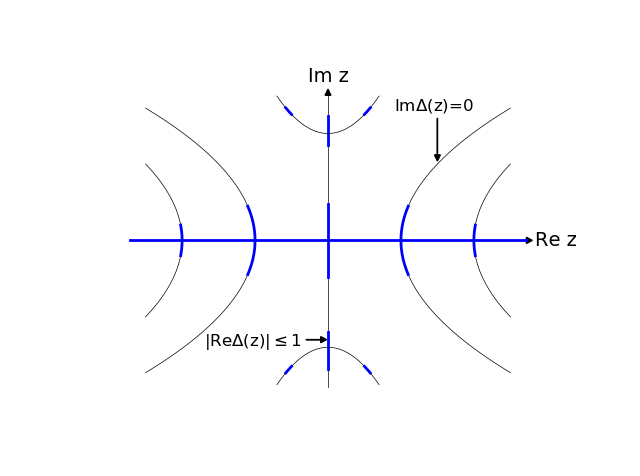}}
\kern-\bigskipamount
\caption{
Left: Schematic illustration of the Lax spectrum of the focusing Zakharov-Shabat scattering problem with a generic periodic potential.
Right: Schematic illustation of the Lax spectrum of the focusing Zakharov-Shabat scattering problem with a real, even, or odd periodic potential. In these cases elements of the Floquet spectrum come in quartets, i.e., $\{z, \overline{z}, -z, -\overline{z} \}$ (see Lemma~\ref{l:symmetries}).} 
\label{f:0}
\end{figure}

With the above definition, 
one can talk about bands and gaps along each $\Gamma_n$ as in a self-adjoint problem.  
The difference is of course that the bands and gaps are not restricted to lie along the real $z$-axis
as they would be in a self-adjoint problem,
but lie instead along arcs of
$\Gamma_n$.
Moreover, 
different curves $\Gamma_i \ne \Gamma_j$ (and therefore different spectral bands) can 
intersect at saddle points of $\D$.
Figure~\ref{f:0} provides a schematic illustration of the Lax spectrum.
Note that 
two curves in $\Gamma$ can intersect at most once, 
as a result of the following:
\begin{lemma}
\label{l:closedband}
The set $\Gamma$, 
and thus the Lax spectrum $\Sigma_{\mathrm{Lax}}$, cannot contain any closed curves in the finite $z$-plane.
\end{lemma}
\kern-\smallskipamount
The equivalent statement to Lemma~\ref{l:closedband} for non-self-adjoint Sturm-Liouville operators is well known \cite{rofebek,rofebek2}.
As in the self-adjoint case, 
the band edges still correspond to the periodic and antiperiodic eigenvalues. 
More precisely, 
for $\D(z) = 1$, 
the corresponding eigenfunctions are periodic with period $L$, 
while for $\D(z) = -1$, 
the corresponding eigenfunctions are antiperiodic (i.e., $2L$-periodic). 
As with Hill's equation, 
some of the spectral gaps might be closed,
in which case, 
for the purposes of this work, 
we will count two adjacent bands 
along a single $\Gamma_n$ as a single one.
On the other hand, 
two intersecting spectral bands lying on different curves $\Gamma_{i} \neq \Gamma_{j}$ 
will be counted as separate.
Given this convention, 
we say $q$ is a \textit{finite-band potential} if its Lax spectrum is composed of finitely many spectral bands. 
   
We note that the terminology ``finite-gap potential'' is much more common in the literature,
especially in the context of special solutions of infinite-dimensional integrable systems (e.g., see~\cite{BBEIM}).
In the self-adjoint case, every finite-gap potential is also a finite-band potential and vice versa, 
so the two concepts are equivalent.
However, 
this is not true in the non-self-adjoint case.

The following property is well known, 
but a proof of it 
is provided in Section~\ref{s:mainproofs} for convenience:
\begin{lemma}
\label{l:infband}
(Infinite band\,,~\cite{ForestLee,MaAblowitz,mclaughlinoverman})
The real $z$-axis is an infinitely long spectral band; that is, $\Real \subset \Sigma_{\mathrm{Lax}}$.   
\end{lemma}
A key to characterize many properties of the spectrum is the asymptotic behavior of the Floquet discriminant as $z\to\infty$.
\begin{lemma}
\label{l:asympDelta}
If $q \in L^{\infty}(\Real)$, 
then for each fixed $\ep>0$ the Floquet discriminant $\D$ has the following asymptotic behavior:
\be
\label{e:Delta_asymp_1}
\D(z) = \begin{cases}
\half \e^{-{\rm i} zL/\ep}(1 + \e^{2{\rm i} zL/\ep} +\, o(1))\,, &z\to \infty\,, \hspace{1mm} \Im z \ge 0\,, \\
\half \e^{-{\rm i} zL/\ep}(1 + \e^{2{\rm i} zL/\ep} +\, O(1/\Im z))\,, &\Im z\to \infty\,, \hspace{1mm} \Im z > 0\,.
\end{cases}
\ee
If $q'\in L^1([0,L])$,
then for each fixed $\ep>0$
\be
\D(z) = 
\half \e^{-{\rm i} zL/\ep} \Big(1 + \e^{2{\rm i} zL/\ep} +\, \frac{1}{2{\rm i} z\ep}(1 - \e^{2{\rm i} zL/\ep})\|q \|_{2}^{2}\, +\, o\Big(\frac{1}{z}\Big)\Big)\,, \quad z\to \infty\,, \hspace{1mm} \Im z \ge 0\,.
\label{e:Delta_asymp_2}
\ee
Moreover,
\be
\D(z) = \cos(zL/\ep) - \frac{1}{2z\ep}\sin(zL/\ep)\| q\|_{2}^2 + o\bigg(\frac{\e^{L\Im z }}{z}\bigg)\,, \quad z \to \infty\,, \hspace{1mm} \Im z \ge 0\,.
\label{e:Delta_asymp_3}
\ee
If, 
additionally,
$q'\in L^{\infty}(\R)$,
then,
\be
\D(z)=\half \e^{-{\rm i} zL/\ep} \Big(1 + \e^{2{\rm i} zL/\ep} +\, \frac{1}{2{\rm i} z\ep}(1 - \e^{2{\rm i} zL/\ep})\|q\|_{2}^{2} + O(1/(\Im z)^2)\Big)\,, \quad \Im z\to \infty\,, \hspace{1mm} \Im z > 0\,.
\label{e:Delta_asymp_qprime}
\ee
\end{lemma}
Further, 
by differentiating~\eqref{e:zs} with respect to $z$ one also has the following:
\begin{lemma}
If $q\in L^{\infty}(\Real)$, 
then for each fixed $\ep>0$
\be
\D'(z) = -\frac{L}\ep \sin(zL/\ep) +  \e^{L\Im z/\ep}o(1)\,, 
\quad z \to \infty\,, 
\hspace{2mm} \Im z \ge 0\,. 
\label{D'-ass}
\ee
If, 
additionally, 
$q' \in L^1([0,L])$,
then for each fixed $\ep>0$
\be
\D'(z) = -\frac{L}{\ep}\sin(zL/\ep) - \frac{L}{2z\ep^2}\cos(zL/\ep)\| q\|_{2}^{2} + 
o\bigg(\frac{\e^{L\Im z}}{z}\bigg)  \,, 
\quad z \to \infty\,, 
\hspace{2mm} \Im z \ge 0\,.
\label{e:D'-ass2}
\ee
\label{l:Dprime}
\end{lemma} 
Note the asymptotic behavior for $\Im z<0$ follows from the Schwarz reflection principle.
Expansions similar to those in Lemma~\ref{l:asympDelta}
were given in~\cite{gesztesyweikard_acta1998,mclaughlinoverman},
but were obtained under the assumptions that the potential is twice differentiable.
Lemmas~\ref{l:asympDelta} and~\ref{l:Dprime} extend the validity of the corresponding results to potentials in $L^\infty(\Real)$.
Lemmas~\ref{l:entire} and~\ref{l:asympDelta} also immediately imply:
\begin{corollary}
\label{l:essentialsing}
If $q \in L^{\infty}(\Real)$, the Floquet discriminant $\D(z)$ 
has an essential singularity at infinity.
\end{corollary}
The fact that $\D$ is entire has a further consequence.
Recall~\eqref{e:bloch}, 
which says that the Floquet spectrum for a given value of $\nu\in\Real$ is the set of $z\in\C$ for which $\D(z)=\cos(\nu L)$. Then one easily gets:
\begin{lemma}
\label{l:discretespec}
For each $\nu \in \Real$, the corresponding Floquet spectrum $\Sigma_\nu$ is discrete.
\end{lemma}
Similarly, 
by~\eqref{e:dirzeros} one has that $\Sigma_{\Dir}(x_o)$ is also discrete.  
Next we derive an upper bound on the imaginary component of points in the spectrum.
\begin{lemma}
\label{l:bound1}
Let $q\in L^\infty(\Real)$. 
Then for all $z \in \Sigma_\Lax$, 
\be
|\Im z| \leq \| q\|_{\infty}\,.
\label{e:inequality1}
\ee
Moreover, 
for all $\zeta\in\Sigma_\Dir(x_o)$, 
\be
|\Im \zeta| \leq \| q\|_{\infty}\,.
\label{e:dirbound1}
\ee
\end{lemma}
The proof of Lemma~\ref{l:bound1} actually yields a stronger estimate: for all $z\in\Sigma_\Lax$, 
\be
|\Im z| \leq | \Re \langle qv_2, v_1 \rangle |/ (\| v_1\|_{2} \| v_2\|_{2}),
\label{e:stronger}
\ee
where 
$\@v(x;z,\ep)=(v_1,v_2)^T$ is an associated Floquet eigenfunction of the ZS problem~\eqref{e:zs},
and $\langle \cdot,\cdot \rangle$ is the $L^2([0,L])$ inner product
(cf.\ \ref{s:notations}).
Inequality~\eqref{e:inequality1} refines the bound obtained in Theorem~4.2 in \cite{gesztesyweikard_acta1998}.
Importantly, 
there exists potentials such that estimate~\eqref{e:inequality1} is sharp (see Section~\ref{s:const} for an example).

Next,
Lemma~\ref{l:asympDelta} leads to the following result:
\begin{lemma}
\label{l:zeros}
Let $q \in L^{\infty}(\Real)$.
Then for each fixed $\ep>0$, the Floquet discriminant $\D(z)$ 
has infinitely many simple real zeros,
and at most finitely many nonreal zeros. 
The same is true for the zeros of $\D'(z)$. 
\end{lemma}
In turn, 
Lemma~\ref{l:closedband}, 
combined with Lemma~\ref{l:asympDelta} provides the key to proving the next result. 
Recall that the resolvent set of an operator is the complement of its spectrum.  
Thus we give the following definition:
\be
\rho_{\Lax} := \Complex\setminus\Sigma_{\Lax}\,.
\label{e:resolvent}
\ee
\begin{theorem}
\label{t:resolvent}
Let $q \in L^{\infty}(\Real)$.
Then the resolvent set $\rho_{\Lax}$ is comprised of two connected components.
\end{theorem}
Importantly, 
Theorem~\ref{t:resolvent} should be compared to that for Hill's equation, 
whose resolvent set is connected~\cite{rofebek,rofebek2}.

By Lemma~\ref{l:entire} and Corollary~\ref{l:essentialsing} one gets that $\D$ is entire as a function of $z$ with an essential singularity at infinity. 
It then follows from Picard's theorem that $\D$ takes on every value infinitely many times with at most one exception.
The following result,
which is a consequence of Lemma~\ref{l:zeros},
clarifies that the lone exception (if one exists) cannot be a value corresponding to the Lax spectrum:
\begin{theorem}
\label{t:infeigenvalues} 
Let $q \in L^{\infty}(\Real)$.
Then for every $\nu\in\R$ the Floquet spectrum $\Sigma_\nu$ is countably infinite. 
\end{theorem}
Next we turn our attention to the number of spectral bands.
As discussed above, 
adjacent bands belonging to the same contour $\Gamma_n$ with a degenerate gap are counted as a single one (and therefore the real axis counts as a single infinitely long band),
but bands belonging to different contours $\Gamma_{i} \neq \Gamma_{j}$ (which can intersect at most once) are counted as separate.
Note that a degenerate gap occurs at a multiple point,
i.e.,
a point $z_{n}^{\rm m}\in\C$ such that $\D^2(z_n^{\rm m})=1$, 
and $\D'(z_n^{\rm m})=0$.
Thus a multiple point is a periodic,
or antiperiodic,
eigenvalue that is also a critical point of the Floquet discriminant.
An important distinction is that in the self-adjoint case necessarily $\D''(z_n^{\rm m})\ne 0$,
while in the non-self-adjoint case one may have higher order zeros.

Recall that,
if the number of spectral bands is finite, 
we call $q$ 
a finite-band potential.
Equivalently, 
$q$ is a finite-band potential if the number of points  
$z_n^{\rm c}\in\R$ such that $\D'(z_n^{\rm c}) = 0$, 
and $\D(z_n^{\rm c}) \in (-1,1)$ is finite.

Next,
let $R_{N}$ denote the rectangle with vertices $ \pm N \pm \i \n q\n_{\infty}$ where $N \in \N$. 
An important consequence of Lemmas~\ref{l:bound1} and~\ref{l:zeros} is the following:
\begin{theorem}
\label{t:finiteband1}
Let $q\in L^\infty(\Real)$ and fix $\ep>0$.
Then $q$ is a finite-band potential if and only if 
$\,\exists N=N(q;\ep)\in\N$ such that $(\Sigma_{\Lax}\setminus\Real)\subset R_{N}$.
\end{theorem}
As a special case, 
Theorem~\ref{t:finiteband1} implies that, 
if the Lax spectrum is contained in the ``cross'' $\Real\cup \i\Real$, 
the potential is finite-band. 
Importantly, 
one can show that, 
if all periodic and antiperiodic eigenvalues are in $\Real\cup \i\Real$, that is enough to conclude that the potential is finite-band.
Specifically: 
\begin{theorem}
\label{t:finiteband2}
Let $q \in L^{\infty}(\Real)$ and suppose
\be
\Sigma_{\pi/L} \cup \Sigma_{2\pi/L} \subset \Real \cup {\rm i} \Real\,,
\label{e:setnotation}
\ee
that is, the periodic/antiperiodic spectrum is only real and purely imaginary. 
Then $q$ is a finite-band potential. 
\end{theorem}
Note that the converse of Theorem~\ref{t:finiteband2} does not hold.  
That is, 
there are finite-band potentials whose periodic and antiperiodic eigenvalues are not only real or imaginary.
(For example, 
any Galilean transformation of a potential obviously preserves its finite-band nature, 
but shifts the Lax spectrum horizontally,
and therefore moves any imaginary eigenvalues off the imaginary $z$-axis.)
Nonetheless, 
Theorem~\ref{t:finiteband2} will be relevant in Section~\ref{s:real}, 
where we study potentials whose spectrum possesses
additional symmetries.

\begin{figure}[t!]
\vglue-4ex
\centerline{\includegraphics[width=8.5cm]{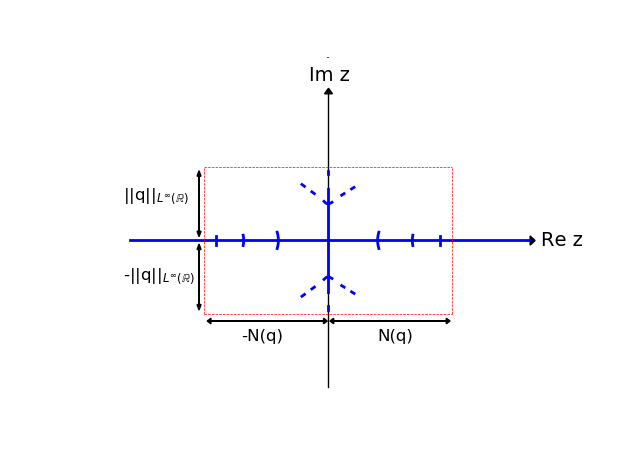}%
\kern-7.4mm
\includegraphics[width=8.5cm]{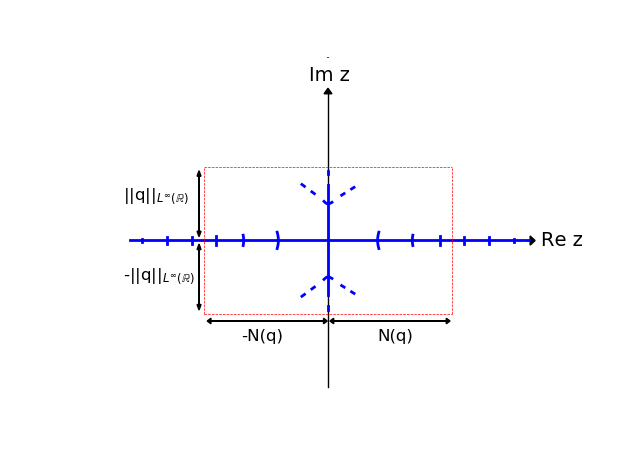}}
\vglue-4ex
\caption{Left: Schematic diagram of the spectrum for a finite-band potential. 
Right: Schematic diagram of the spectrum for an infinite-band potential}
\label{f:1}
\end{figure}

Next we discuss the asymptotic distribution of bands as $z\to\infty$. 
For this purpose, 
following \cite{mclaughlinoverman}, 
we introduce the concept of ``spine'', 
defined as a spectral band that intersects the real axis transversally and does not intersect any other band. 
Recall that any intersection point between two or more bands is a saddle point (or critical point) of the discriminant, 
that is, 
a point $z_n^{\rm c}\in\C$ such that $\D'(z_n^{\rm c})=0$.
\begin{theorem}\label{theo-spines}
For any $q \in L^{\infty}(\Real)$ and each fixed $\ep>0$ there exists some $N=N(q;\ep)\in\N$, such that 
all but finitely many bands of the Lax spectrum are spines located outside $R_{N}$.
Moreover: 
i) for each of these spines there exists $n\in\Z$ such that  
the intersection point between the spine and the real axis is $o(1)$-close to the point~${n \pi\ep}/L$
as $n\to \pm\infty$; 
ii) only one spine can be  $o(1)$-close to ${n \pi\ep}/L$ as $n\to \pm\infty$.
Further, 
if $q' \in L^1([0,L])$, the intersection point is $O(1/z)$-close. 
\end{theorem}
Related asymptotic results were given in \cite{djakovmityagin,gesztesyweikard_acta1998,mclaughlinoverman}.
As a consequence of Theorem~\ref{theo-spines}, 
we have that the potential $q \in L^{\infty}(\Real)$ is a finite-band potential if and only if $\exists N\in\N$ such that $\Sigma_{\rm Lax}$ contains no spines outside $R_N$.  
An illustration of Theorem~\ref{theo-spines} is given in Fig.~\ref{f:1}.

Theorem~\ref{theo-spines} has an important consequence.
Since a spectral gap is delimited by the presence of a spectral band on either side of it,
the number of gaps associated to any given potential is always finite,
since for any potential there is a finite region of the complex plane outside of which all bands are spines
(for which there are no associated gaps), 
and the number of spectral bands (and therefore gaps) inside any finite region of the complex plane is always finite
(as a consequence of the fact that $\D(z)$ is entire).
In other words, 
\textit{for the focusing Zakharov-Shabat operator on the circle, 
every potential is a finite-gap potential}
(cf.\ \cite{gesztesyweikard_bams1998}).
This situation is in marked contrast with the case of the spectral problem for the Hill operator and that for
the self-adjoint Zakharov-Shabat operator, 
for which finite-gap potentials are a special subset.
In contrast, 
for the focusing Zakharov-Shabat operator the only meaningful distinction is that
between finite-band and infinite-band potentials
(cf.\ \cite{gesztesyweikard_bams1998}).
The key distinction is the existence of finitely many (as opposed to infinitely many) simple periodic/antiperiodic eigenvalues,
which defines a Riemann surface of finite-genus~\cite{BBEIM,gesztesyweikard_acta1998, gesztesyweikard_bams1998}. 

\section{Properties of the spectrum in the semiclassical limit}
\label{s:semiclassical}

In this section we describe the Lax spectrum in the semiclassical limit for complex-valued potentials that do not depend on $\ep$ (see Section~\ref{s:examples}).
Before we do so, 
however, 
we discuss further rigorous bounds on the location of the Lax spectrum under fairly weak assumptions on the potential.

\begin{lemma}
Suppose $q' \in L^{\infty}(\Real)$. 
If $z \in \Sigma_{\Lax}$,
then
\be
|\Re z||\Im z| \leq \frac{\epsilon}{2} \| q'\|_{\infty}\,.
\label{e:inequality2}
\ee
\label{l:bound2}
\end{lemma}
The proof of Lemma~\ref{l:bound2} actually yields the stronger estimate
\be
|\Re z||\Im z| \leq \halfep |\Im \langle q'v_2, v_1 \rangle |/ (\| v_1\|_2 \| v_2\|_2),
\ee
where 
$\@v(x;z,\ep) = (v_1,v_2)^T$ 
is an associated Floquet eigenfunction of the ZS system~\eqref{e:zs},
and $\langle \cdot,\cdot \rangle$ denotes the $L^{2}([0,L])$ inner product (cf.~\ref{s:notations}).
Thus, 
if $z\in\Sigma_{\Lax}$, 
then Lemmas~\ref{l:bound1} and~\ref{l:bound2} together give the bound
\be
| \Im z | \leq \min \Big \{ \| q\|_{\infty} \;, \; \frac{\ep\| q'\|_{\infty}}{2|\Re z|} \Big \}\,.
\label{e:maininequality}
\ee
%

\begin{figure}[t!]
\vglue-2ex
\centerline{\includegraphics[width=7.5cm]{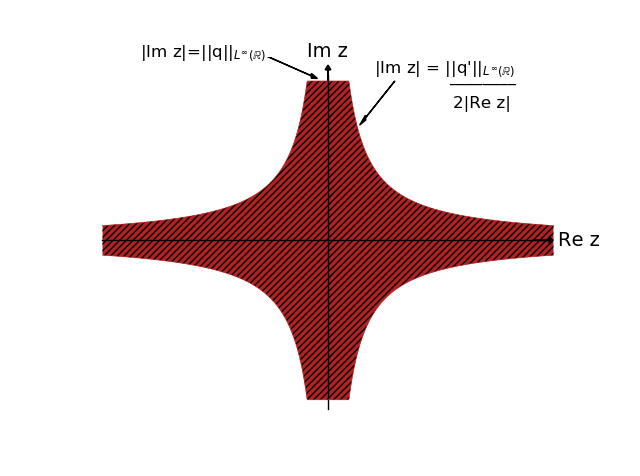}
\includegraphics[width=7.5cm]{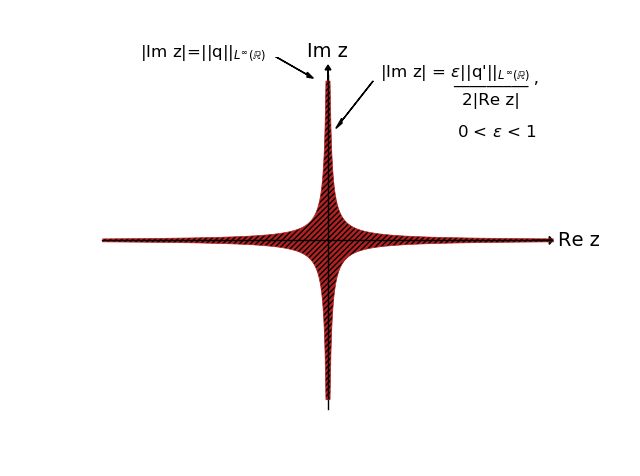}}
\kern-\bigskipamount
\caption{Left: Bounds on the Lax spectrum for $\ep=1$ (dark red). Right: Bounds on the Lax spectrum when $0<\epsilon<1$ (dark red).}
\label{f:2}
\end{figure}

\begin{theorem}
\label{t:set1}
Suppose $q'\in L^{\infty}(\Real)$. 
Then
\be
\Sigma_{\Lax} \subset \Lambda^{\ep}(q)\,,
\label{setinclusion}
\ee
where
\be
\Lambda^{\ep}(q):=\Big\{z \in \Complex : |\Im z| \leq \| q\|_{\infty} \Big\} \cap \Big\{z \in \Complex : |\Re z| |\Im z| \leq \frac{\ep}{2}\| q'\|_{\infty} \Big\}\,.
\label{e:eigenvalue_set}
\ee
\end{theorem}
That is, 
the Lax spectrum is contained in the set $\Lambda^{\ep}(q)$ depending explicitly on $\ep$,
the semiclassical parameter.
The region $\Lambda^{\ep}(q)$ is shown in dark red in Fig.~\ref{f:2}. 

Theorem~\ref{t:set1} complements the localization result of~\cite{djakovmityagin}
for the periodic and antiperiodic eigenvalues.
Moreover, 
the results in Lemmas~\ref{l:bound1} and~\ref{l:bound2} apply to the entire Lax spectrum, not just the periodic and antiperiodic eigenvalues.

Also, 
recall that cardinality results for the 
periodic, 
antiperiodic, 
and Dirichlet eigenvalues were obtained in~\cite{kapmity,kappeler,mclaughlinoverman}.
Namely, 
$\exists N=N(q)\in\N$ such that in the disc $D(0;(N-1/2)\pi/L)$ there are exactly $2N-2$ periodic eigenvalues, 
$2N$ antiperiodic eigenvalues,
and $2N-1$ Dirichlet eigenvalues.  
Importantly,
for $q\in H^1([0,L])$ one can get an explicit estimate for $N=N(q)$ (see~\cite{mclaughlinoverman} for details).
Theorem~\ref{t:set1} allows one to strengthen those results.
For brevity, 
we let $\ep = 1$ in Theorem~\ref{t:Pi&Xi} and Corollary~\ref{t:eigvalsorder}.
\begin{theorem}
\label{t:Pi&Xi}
Suppose $q' \in L^{\infty}(\Real)$ and $\ep=1$. 
Consider the sets 
\bse
\begin{gather}
\Pi_{\rm Lax}(q):=\Lambda^{\ep}(q) \cap D(0;(N-1/2)\pi/L)\,,
\label{e:card1}
\\
\Pi_{\rm Dir}(q):=\big(\R \times {\rm i}[-\| q\|_{\infty}, \| q\|_{\infty}]\big) \cap D(0;(N-1/2)\pi/L)\,.
\label{e:card2}
\end{gather}
\ese
The set $\Pi_{\rm Lax}(q)$ contains exactly $2N-2$ periodic eigenvalues and $2N$ antiperiodic eigenvalues. 
The set $\Pi_{\rm Dir}(q)$ contains exactly $2N-1$ Dirichlet eigenvalues.
\end{theorem}
A consequence of Lemma~\ref{l:asympDelta} and Lemma~\ref{l:bound2} is the following: 
\begin{corollary}
\label{t:eigvalsorder}
Fix $\nu\in\Real$.
Suppose $q' \in L^{\infty}(\Real)$ and $\ep=1$.
Then the Floquet eigenvalues $\{z_n(\nu)\}_{n\in \Z}$ can be partially ordered such that
$\Re z_n(\nu) \to \pm\infty$ as $n \to \pm\infty$.
Moreover, 
\begin{gather}
\Im z_n(\nu) = O(1/\Re z_n(\nu))\,, \quad {\rm as} \hspace{1mm} n\to\pm\infty\,, \hspace{1mm} \Re z \ne 0\,,
\label{e:eigenlims2}
\end{gather}
\end{corollary}

We now turn to the semiclassical limit.
Theorem~\ref{t:set1} allows us to give a rigorous characterization of the Lax spectrum as $\ep \downarrow 0$.
This is important since for certain initial data the focusing NLS equation on the circle appears to admit a coherent structure in the semiclassical limit~\cite{BiondiniOregero}.
Note that spectral confinement to the real and imaginary axes together with the fact that the real $z$-axis is an infinitely long band implies that all nonlinear excitations emerging from the input have zero velocity.
\begin{corollary}
\label{c:mainresult}
Suppose $q'\in L^{\infty}(\Real)$. 
Define
\be
\label{S-infty}
\Sigma_\infty := \Real \cup {\rm i}[-\| q\|_{\infty},\,\| q\|_{\infty}]\,.
\ee
Moreover, 
let $N_\delta(\Sigma_\infty)$ be a $\delta$-neighborhood of $\Sigma_\infty$.
Then for any $\delta>0$, 
\be
\Sigma_{\Lax} \subset N_{\delta}(\Sigma_{\infty})\,, 
\label{e:spectrumlimit}
\ee
for all sufficiently small values of $\ep$.
That is, 
for any $\delta>0$ there exists an $\ep_*>0$ such that \eqref{e:spectrumlimit} holds for all $0<\ep \leq \ep_*$.
\end{corollary}

\begin{figure}[t!]
\vglue-2ex
\centerline{\includegraphics[width=7.5cm]{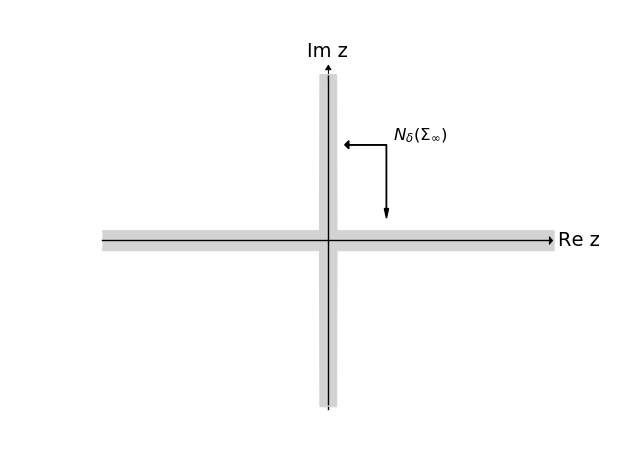}
\includegraphics[width=7.5cm]{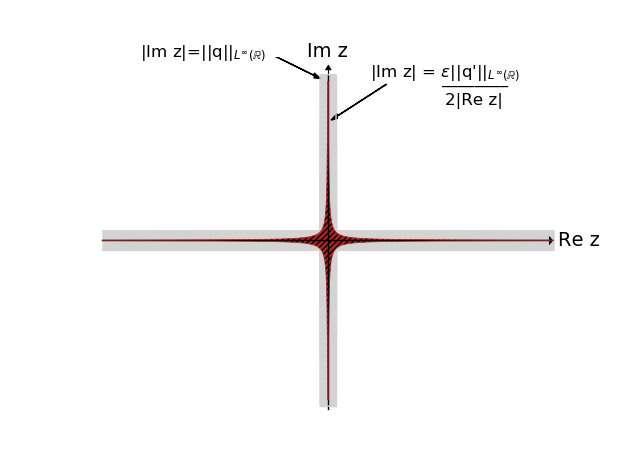}}
\kern-\bigskipamount
\caption{
Left: The $\delta$-neighborhood of $\Sigma_{\infty}$ (light gray). Right: Bounds on the Lax spectrum for $0 < \epsilon \leq \epsilon_{*}$ (dark red) with 
the $\delta$-neighborhood of $\Sigma_{\infty}$ (light gray). 
}
\label{f:3}   
\end{figure}

Corollary~\ref{c:mainresult} is a direct consequence of Theorem~\ref{t:set1}.
Moreover, 
Corollary~\ref{c:mainresult} corroborates and puts on a rigorous footing the numerical results of~\cite{BiondiniOregero}.
Roughly speaking, 
one can interpret Corollary~\ref{c:mainresult} as saying that the Lax spectrum of $\Leps$ ``converges'' to $\Sigma_\infty$ in the limit $\epsilon\downarrow 0$.
Note, 
however, 
that in general the semiclassical limit of \eqref{e:zs} is a singular limit.
This is why it is necessary to state the result in a more indirect way. 
(It was also shown in~\cite{BiondiniOregero} through numerical and asymptotic calculations that the number of bands is $O(1/\ep)$ as $\ep\downarrow 0$.
However, 
no rigorous proof of this result exists to the best of our knowledge.)
An illustration of the $\delta$-neighborhood $N_\delta(\Sigma_\infty)$ and its relation to the Lax spectrum is given in Fig.~\ref{f:3}.
The proof of Theorem~\ref{c:mainresult} and of all results in this section is given in Section~\ref{s:semiclassicalproofs}.

\section{Properties of the spectrum for real or symmetric potentials}
\label{s:real}
In addition to Schwarz symmetry, 
the Lax spectrum acquires additional symmetries when the potential in~\eqref{e:zs} is itself symmetric. 
In particular, 
if the potential $q$ is real then it follows that the monodromy matrix satisfies
\be
\overline{M(-\overline{z};\ep)} = M(z;\ep)\,.
\label{e:msym2}
\ee
Next we list some additional symmetries of the monodromy matrix when the potential is symmetric. 
These symmetries are easily deduced from the symmetries of the ZS system~\eqref{e:zs} (cf.~\ref{s:symmetries}).
Suppose that the potential satisfies a generalized reflection symmetry,
that is,
$\exists\theta\in\Real$ such that
$q(-x) = \e^{2\i\theta}q(x)$ for all $x\geq 0$. 
(Obviously for $\theta = 0\; \text{mod} \;\pi$, and $\theta = \pi/2\; \text{mod} \;\pi$ one has the cases of even and odd potentials, respectively.)
Then one gets
\be
\overline{M(-\overline{z};\ep)^{-1}} = (\cos \theta \sigma_1 + \sin \theta \sigma_2) M(z;\ep) (\cos \theta \sigma_{1} + \sin \theta \sigma_{2})\,,
\label{e:msym3}
\ee
where $\sigma_{1}$ and $\sigma_2$ are the first and second Pauli matrices, respectively (cf. \ref{s:notations}). 

If the potential in~\eqref{e:zs} is PT-symmetric, 
that is,
if $q(-x) = \overline{q(x)}$, 
then it follows that the monodromy matrix satisfies
\be
\overline{M(\overline{z};\ep)^{-1}} = \sigma_3M(z;\ep)\sigma_3\,.
\label{e:ptmsym}
\ee
Using the above symmetries, 
for odd real, 
or even real potentials 
we obtain
\be
M(z;\ep) = \begin{pmatrix} c(z;\ep) & \pm s(z;\ep) \\ s(z;\ep) & (1 \pm s^2(z;\ep))/{c(z;\ep)} \end{pmatrix}\,,
\label{e:ptsymmetry5}
\ee
where $c = M_{11}(z;\ep)\not\equiv 0$, 
and $s = M_{21}(z;\ep)$,
and where the ``+'' sign holds for odd real potentials and the ``$-$'' sign for even real potentials.
Using the above symmetries we get the following:
\begin{lemma}
\label{l:symmetries}
Let the potential $q$ satisfy at least one of the following conditions: 
a) it is real;
b) it is even;
c) it is odd. 
Then $\Sigma_{\Lax}$ is symmetric with respect to the imaginary $z$-axis,
and $\Gamma$ [cf.~\eqref{e:contours}] includes the $\Im z$-axis.
\end{lemma}
Next, using Lemma~\ref{l:symmetries} we get the following:
\begin{theorem}
\label{t:reimag}
Let the potential $q$ satisfy at least one of the following conditions: 
a) it is real;
b) it is even;
c) it is odd. 
If the periodic and antiperiodic spectra are real and purely imaginary only,  
then the entire Lax spectrum is contained within the real and imaginary axes, 
that is, 
$\Sigma_{\Lax} \subset \Sigma_\infty\subset\Real \cup {\rm i}\Real$, 
[cf.~\eqref{S-infty}].
\end{theorem}

Further, 
by assuming the potential is real, 
or symmetric one is able to obtain stronger bounds on the Dirichlet and Floquet spectra:
\begin{lemma}
\label{l:dirbounds2}
Suppose $q' \in L^{\infty}(\R)$. 
Moreover, 
let the potential $q$ satisfy at least one of the following conditions: 
a) it is real;
b) it is odd;
c) it is {\rm PT}-symmetric.
If $\zeta \in \Sigma_{\Dir}(0)$, 
then
\be
|\Re \zeta|\Im \zeta| \leq \frac{\ep}{2}\n q'\n_{\infty}\,.
\label{e:estimate4}
\ee
\end{lemma}
Importantly,
for $q$ real the result holds more generally for $\zeta\in\Sigma_{\Dir}(x_o)$.
That is,
the bound~\eqref{e:estimate4} is independent of the base point $x_o$ when $q$ is real.
\begin{theorem}
\label{t:dirset2}
Assume $q'\in L^{\infty}(\R)$. 
Let the potential $q$ satisfy at least one of the following conditions: 
a) it is real;
b) it is odd;
c) it is {\rm PT}-symmetric.
Then
\be
\Sigma_{\Dir}(0) \subset \Lambda^{\ep}(q)\,,
\label{e:dirsubset}
\ee
where
\be
\Lambda^{\ep}(q) := \{\zeta \in \C : |\Im \zeta| \leq \| q\|_{\infty} \} \cap \{\zeta \in \C : |\Re \zeta||\Im \zeta| \leq \halfe\| q'\|_{\infty}\}.
\label{e:dirregion}
\ee
\end{theorem}
Thus, 
for these classes of potentials the Dirichlet spectrum localizes to the real and imaginary axes in the semiclassical limit.
An interesting open question is whether the bound obtained in Lemma~\ref{l:dirbounds2} holds more generally for complex-valued potentials.
These results should be compared to those pertaining to the Lax spectrum.
Finally,
we have the following interesting bound:
\begin{lemma}
\label{l:wkbbound}
Let $q$ be real,
strictly positive,
and suppose that $q'\in L^{\infty}(\R)$.
Then,
\be
|\Im z| \leq \frac{\ep}{2} \|(\ln q)'\|_{\infty}\,,
\label{e:wkbbound}
\ee
for all $z\in\Sigma_{\Lax}\setminus \rm{i}\R$.
Moreover,
\be
|\Im \zeta| \leq \frac{\ep}{2}\|(\ln q)'\|_{\infty}\,,
\label{e:wkbboundDir}
\ee
for all $\zeta\in\Sigma_{\Dir}(x_o)\setminus {\rm i}\R$.
\end{lemma} 
Thus,
for strictly positive real potentials
$\Im z(\nu;\ep)$ (resp. $\Im \zeta(\ep)$) is uniformly $O(\ep)$ as $\ep\downarrow 0$ for non-purely imaginary $z\in\Sigma_{\Lax}$ (resp. $\zeta\in\Sigma_{\Dir}(x_o)$).
(Note by phase invariance this result can easily be extended to strictly negative potentials as well.)
This result provides strong justification of the WKB analysis employed by two of the authors in~\cite{BiondiniOregero}. 
Moreover, 
Lemma~\ref{l:wkbbound} is the analogue for periodic potentials of the classic results
about purely imaginary discrete eigenvalues for potentials on the line~\cite{BL2018,KS2002}.
Finally,
localization along the imaginary axis of the spectral variable as in Lemma~\ref{l:wkbbound} together with the WKB results obtained in~\cite{BiondiniOregero} find application in the study of soliton and breather gases in focusing nonlinear media (see~\cite{eltovbisSG}).

When the potential in the ZS system~\eqref{e:zs} is real, 
one can map~\eqref{e:zs} to a non-self-adjoint Hill's equation.
Specifically, 
the invertible change of variables 
\be
y_{\pm}(x;\lambda,\ep) = v_{1}(x;z,\ep) \pm \i v_{2}(x;z,\ep)\,
\label{e:changevariab}
\ee
maps~\eqref{e:zs} into linear second-order ODEs
\be
\big(-\ep^{2} \partial_{x}^{2} + W_\pm(x;\ep) \big)y_\pm = \lambda y_\pm\,,
\label{e:schrodinger}
\ee
which are of course the celebrated
Hill's equation, 
with spectral parameter $\lambda:=z^{2}$ and complex potentials 
\be
W_\pm(x;\ep) := -q^{2}(x) \mp \i\ep q'(x)\,.
\label{e:potentials}
\ee
Because $W_\pm(x;\ep)$ are not real, 
in general the eigenvalue problems~\eqref{e:schrodinger} are non-self-adjoint.
Note however that, 
when $q$ is real and even, 
the potentials in~\eqref{e:potentials} are PT-symmetric, 
i.e., 
invariant under the combined action of space reflections and complex conjugation~\cite{RamyEl}. 
So the results above for real and even potentials also serve as a further study of spectral problems for the Hill operator with a PT-symmetric potential.
Note, 
estimate~\eqref{e:inequality1} gives $\Re \lambda \ge -\n q\n_{\infty}^2$,
and the estimate~\eqref{e:inequality2} gives $\Im \lambda = O(\ep)$ as $\ep \downarrow 0$ in~\eqref{e:schrodinger}. 
Theorem~\ref{t:reimag} implies that, 
even though real potentials are associated with a Hill operator with a complex potential,
if the periodic and antiperiodic Floquet eigenvalues are all real or purely imaginary,
then the corresponding non-self-adjoint Hill equation~\eqref{e:schrodinger} has a purely real spectrum.
Moreover,
this implies that $W_{\pm}$ is a finite-band potential (recall that for the ZS system the real $z$-axis is an infinitely long band). 
Hence, 
if $q\in L^{\infty}(\Real)$, 
then along the real $\lambda$-axis the spectral bands and gaps are confined to the interval $[-\|q\|_{\infty}^2, 0 \big)$ followed by one infinitely long spectral band along the interval $[0, \infty)$.
Note however the thesis of Theorem~\ref{t:reimag} is \textit{not} true if the potential is not real or even or odd, 
because in that case the spectrum possesses no left-right symmetry.
The proofs of all the results in this section are given in Section~\ref{s:symmetricproofs}.

\section{Proofs}
\label{s:proofs}

\subsection{Proofs: General properties of the spectrum} 
\label{s:mainproofs}

We refer the reader to~\ref{s:notations} for definitions.

\textit{Proof of Lemma~\ref{l:closedband}.}
By Lemma~\ref{l:entire}, 
$\D(z)$ is an entire function differing from a constant. 
Thus, 
$\Im \D(z)$ is harmonic.
By the maximum principle the set $\Gamma$ cannot contain any closed curve. 
Further, 
along any spectral band we have $\Im\D(z)=0$.
Hence, $\Sigma_{\Lax}$ contains no closed bands in the finite complex $z$-plane.
\qed

\textit{Proof of Lemma~\ref{l:infband}.}
By $\eqref{e:Mrepresentation}$,
$\D(z) \in \Real$ along the real $z$-axis,
and
$|c(z;\ep)|^2 + |s(z;\ep)|^2 \equiv 1$.
Further, along the real $z$-axis
$\D(z) = \Re c(z;\ep)$.
Thus, 
$|\D(z)| \leq (1 - |s(z;\ep)|^2)^{1/2} \leq 1$.
Hence, 
it follows that $\Real \subset \Sigma_{\Lax}$. 
\qed

In order to prove Lemma~\ref{l:asympDelta}, 
we need two intermediate results.
Recall that one can write
\be
M(z;\ep) = \Phi(L;z,\ep)\,,
\ee
where $\Phi(x;z,\ep)$ is the principal matrix solution of~\eqref{e:zs}.
Thus \eqref{e:DeltaM} implies $\D(z)=\half\tr \Phi(L;z,\ep)$.
Accordingly,~Lemma~\ref{l:asympDelta} is a consequence of the asymptotic behavior of $\Phi(x;z,\ep)$ as $z\to\infty$.
Introducing the change of dependent variable
\be
\mu(x;z,\ep) := \Phi(x;z,\ep)\e^{\i zx \sigma_{3}/\ep}\,,
\label{e:modified}
\ee
as well as the notation $\mu_{ij}$ for $i,j=1,2$ to denote matrix entries,
we next prove:
\begin{lemma}
\label{l:asympmu}
Let $q \in L^{\infty}([0,L])$,
and let $\Im z \geq 0$.
Then for each fixed $\ep>0$
\bse
\be
|\mu_{11}(x;z,\ep)-1| \leq \frac{\e^{M}L\|q\|_{\infty}^2}{2\ep\Im z}\,, \quad \Im z > 0\,,
\label{e:mu11asymp_imz>0}
\ee
where $M = (L/\ep)^2\|q\|_{\infty}$. Moreover,
\be
\mu_{11}(x;z,\ep) = 1 + o(1)\,, \quad z\to\infty\,, \quad \Im z \ge 0\,,
\label{e:mu11asymp}
\ee
\ese
uniformly for $x\in[0,L]$.
%
If $q'\in L^1([0,L])$,
then for each fixed $\ep>0$
\be
\mu_{11}(x;z,\ep) = 1 + \frac{1}{2{\rm i} z\ep} \int_{0}^{x} |q(y)|^2\,\d y + o(1/z)\,, \quad z \to \infty\,, \quad \Im z \ge 0\,.
\label{e:mu11asymp2}
\ee
If,
additionally,
$q'\in L^{\infty}([0,L])$,
then,
\be
\mu_{11}(x;z,\ep) = 
1 + \frac{1}{2{\rm i} z\ep} \int_{0}^{x} |q(y)|^2\,\d y + O(1/(\Im z)^2)\,, \quad \Im z \to \infty\,, \quad \Im z > 0\,.
\label{e:mu11asymp2star}
\ee
\end{lemma}
\textit{Proof of Lemma~\ref{l:asympmu}.}
Using~\eqref{e:modified}, 
it follows that $\mu = (\mu_{ij})$ satisfies the following system of ODEs: 
\be
\ep\mu' = -\i z [\sigma_3, \mu] + Q(x)\mu\,,
\label{e:ode_mu}
\ee
where $[\sigma_3, \mu] := \sigma_3\mu - \mu\sigma_3$. 
Thus, 
$\mu$ satisfies the matrix Volterra linear integral equation
\be
\mu(x;z,\ep) = \mathbb{I} + \int_{0}^{x} \e^{-\i z(x-y)\sigma_{3}/\ep}Q(y)\mu(y;z,\ep)\e^{\i z(x-y)\sigma_{3}/\ep} \d y/\ep\,.
\label{e:inteq}
\ee
Explicitly~\eqref{e:inteq} yields coupled integral equations for the individual entries of $\mu(x;z,\ep)$.
In particular,
\bse
\begin{align}
\mu_{11}(x;z,\ep) &= 1 + \int_{0}^{x} q(y)\mu_{21}(y;z,\ep)\, \d y/\ep \,,
\label{e:mu11inteqn}
\\
\mu_{21}(x;z,\ep) &= -\int_{0}^{x} \e^{2\i z(x-y)/\ep}\overline{q(y)}\mu_{11}(y;z,\ep)\, \d y/\ep \,,
\label{e:mu21inteqn}
\end{align}
\ese
with similar equations for $\mu_{12}(x;z,\ep)$ and $\mu_{22}(x;z,\ep)$.
Eliminating $\mu_{21}$,
one then obtains
\bse
\begin{gather}
\mu_{11}(x;z,\ep) = 1 + K(\mu_{11})(x;z,\ep)\,,
\label{e:mu11}
\\
\noalign{where}
K(f)(x;z,\ep) := \int_{0}^{x} k(x,t;z,\ep) f(t)\, \d t/\ep\,,
\\
k(x,t;z,\ep) := -\overline{q(t)}\int_{t}^{x} \e^{2\i z(y-t)/\ep} q(y)\, \d y/\ep\,, \hspace{3mm} x>t\,.
\label{e:kern}
\end{gather}
\ese
Note that, $\forall x\in[0,L]$ and $\forall t\in[0,x]$, $|k(x,t;z,\ep)|\leq L\n q \n_{\infty}^{2}/\ep$.
Then,
introducing the Neumann series
\bse
\begin{gather}
\mu_{11}(x;z,\ep) = \sum_{n=0}^\infty\mu_{11}^{(n)}(x;z,\ep)\,,
\\
\mu_{11}^{(0)} = 1,\qquad \mu_{11}^{(n+1)}(x;z,\ep) = K(\mu_{11}^{(n)}) = K^{n+1}(\mu_{11}^{(0)})\,,
\end{gather}
\ese
we have the bound
\begin{align*}
|K^{n}(f)(x;z,\ep)| \leq 
\int_{0}^{x}|k(x,t_1;z,\ep)| 
\cdots \int_{0}^{t_{n-1}}|k(t_{n-1},t_n;z,\ep)||f(t_n)| \d t_n \cdots \d t_1/\ep^n
\leq \frac{M^n}{n!} \|f\|_{\infty} \,,
\end{align*}
where $M:=(\txtfrac{L}{\ep})^{2}\|q\|_{\infty}^{2}$. 
Hence, the series  is absolutely convergent,
and 
\be 
\sup_{x \in [0,L], \Im z \geq 0} |\mu_{11}(x;z,\ep)| \leq \sum_{n=0}^{\infty} |K^{n}(1)(x;z,\ep)| \leq \e^{M}\,,
\label{e:mu11bound}
\ee
since $\mu_{11}^{(0)}\equiv 1$. 
Next, 
from~\eqref{e:mu21inteqn} and~\eqref{e:mu11bound} we immediately have 
$|\mu_{21}(x;z,\ep)| \leq \e^{M}\n q\n_{\infty}/(2\Im z)$ $\forall \Im z > 0$. 
Thus, using~\eqref{e:mu11inteqn} one gets~\eqref{e:mu11asymp_imz>0}. 

Next we consider the case $z\to \infty$ while $\Im z$ remains bounded. 
The Riemann-Lebesgue lemma applied to~\eqref{e:kern}
implies that, 
$\forall x\in[0,L]$ and $\forall t\in[0,x]$, $k(x,t,z,\ep) = o(1)$ as $z\to\infty$ from $\Im z \ge 0$.
Moreover, it is straightforward to show that, as in \cite{Eastham},
$q\in L^\infty([0,L])$ is sufficient to ensure that the $o$ symbol is uniform with respect to $x\in[0,L]$.
(E.g., this can be done by approximating $q$ in~\eqref{e:kern} with smooth functions, which are
dense in $L^1([0,L])$.)
Note that \cite{Eastham} considers a Sturm-Liouville problem, which involves not only $q$ but also~$q'$,
and absolute continuity of $q$ implies the existence of $q'\in L^1([0,L])$.
Here, $q\in L^\infty([0,L])$ implies $q\in L^1([0,L])$.
Together with the uniform boundedness of $\mu_{11}$ for $\Im z \ge 0$, this
gives
\vspace*{-1ex}
\be
|K(\mu_{11})(x;z,\ep)| \leq \e^{M} \int_{0}^{x} |k(x,t;z,\ep)|\, \d t/\ep = o(1)\,, \quad z \to \infty\,, \quad \Im z \ge 0\,.
\label{e:kernest}
\ee
Thus, 
$\mu_{11}(x;z,\ep) = 1 + o(1)$ as $z \to \infty$ from $\Im z \ge 0$,
i.e.,~\eqref{e:mu11asymp}.
In turn, using~\eqref{e:mu21inteqn} and recalling~\eqref{e:mu11bound} one gets immediately
\be
\mu_{21}(x;z,\ep) = 
o(1)\,, \quad z \to \infty\,, \quad \Im z \ge 0\,.
\label{e:mu21estimate}
\ee
Now assume additionally that $q' \in L^1([0,L])$. From~\eqref{e:ode_mu} it follows $\mu_{11}' = q(x)\mu_{21}/\ep$.
Then integrate~\eqref{e:mu21inteqn} by parts:
\be
\mu_{21}(x;z,\ep) = \frac{1}{2\i z}\big(\overline{q(x)}\mu_{11}(x;z,\ep) - \e^{2\i zx/\ep}\overline{q(0)}\big) - \frac{1}{2\i z} \int_{0}^{x} \e^{2\i z(x-y)/\ep}(\overline{q(y)}\mu_{11}(y;z,\ep))' \d y\,,
\label{e:parts1}
\ee
where $\mu_{11}(0;z,\ep) \equiv 1$. 
This implies,
for $q'\in L^{\infty}([0,L])$,
\be
\mu_{21}(x;z,\ep) = \frac{1}{2\i z} \big(\overline{q(x)}\mu_{11}(x;z,\ep) - \e^{2\i zx/\ep}\overline{q(0)}\big) + O(1/(\Im z)^2)\,, \quad \Im z \to \infty\,, \hspace{1mm} \Im z>0\,.
\label{e:parts2}
\ee
Plugging~\eqref{e:parts2} into~\eqref{e:mu11inteqn} gives~\eqref{e:mu11asymp2star}.
Finally, note,
by~\eqref{e:mu21estimate},
and since $\mu_{11}'$ is proportional to $\mu_{21}$, 
it follows that $\mu_{11}'(x;z,\ep) = o(1)$ as $z \to \infty$ from $\Im z \ge 0$.
Hence, 
by~\eqref{e:parts1} one gets
\be
\mu_{21}(x;z,\ep) = \frac{1}{2\i z} \big( \overline{q(x)}\mu_{11}(x;z,\ep) - \e^{2\i zx/\ep}\overline{q(0)} \big) + o\Big(\frac{1}{z}\Big)\,, \hspace{2mm} z \to \infty\,, \hspace{2mm} \Im z \ge 0\,.
\label{e:mu21parts}
\ee
Plugging~\eqref{e:mu21parts} into~\eqref{e:mu11inteqn}
gives~\eqref{e:mu11asymp2}.
\qed

Next, 
since $\mu_{22}(x;z,\ep)$ is unbounded as $\Im z \to \infty$ from $\Im z > 0$ we make a further change of variables, 
namely,
\be
(\tilde{\mu}_{12}, \tilde{\mu}_{22})^{T} := \e^{2\i zx/\ep}(\mu_{12}, \mu_{22})^{T}\,,
\label{e:changevar}
\ee
and examine the asymptotic behavior of $\tilde{\mu}_{22}(x;z,\ep)$ as 
$z \to \infty$ from $\Im z \ge 0$. 
\begin{lemma}
\label{l:mutildasym}
Let $q \in L^{\infty}([0,L])$,  
and let $\Im z \geq 0$.
Then for each fixed $\ep>0$ 
\bse
\be
|\tilde{\mu}_{22}(x;z,\ep) - \e^{2{\rm i} zx/\ep}| \leq \frac{\e^{M}L\|q\|_{\infty}^{2}}{2\ep\Im z}\,, \quad \Im z > 0\,,
\label{e:mu22tildasymp_imz>0}
\ee
where $M= ({L}/{\ep})^2\|q\|_{\infty}^2$. Moreover,
\be
\tilde{\mu}_{22}(x;z,\ep) = \e^{2{\rm i} zx/\ep} +\, o(1)\,, \quad z\to\infty\,, \quad \Im z \geq 0\,,
\label{e:tildmu22asymp}
\ee
\ese
uniformly for $x\in[0,L]$. 
If $q' \in L^1([0,L])$,
then for each fixed $\ep>0$
\be
\tilde{\mu}_{22}(x;z,\ep) = \e^{2 {\rm i} zx/\ep} - \frac{\e^{2{\rm i} zx/\ep}}{2 {\rm i} z\ep} \int_{0}^{x} |q(y)|^2\, \d y + o\big(1/z)\,, \quad z \to \infty\,, \quad \Im z \ge 0\,.
\label{e:tildmu22asymp2}
\ee
If,
additionally,
$q'\in L^{\infty}([0,L])$,
then,
\be
\tilde{\mu}_{22}(x;z,\ep) = 
\e^{2{\rm i} zx/\ep} - \frac{\e^{2 {\rm i} zx/\ep}}{2 {\rm i} z\ep} \int_{0}^{x} |q(y)|^2\, \d y + O(1/(\Im z)^2)\,, \quad \Im z \to \infty\,, \quad \Im z > 0\,.
\label{e:tildemu22asymp2star}
\ee
\end{lemma}
\textit{Proof of Lemma~\ref{l:mutildasym}.}
Using~\eqref{e:inteq} and~\eqref{e:changevar} we have
\bse
\begin{align}
\tilde{\mu}_{12}(x;z,\ep) &= \int_{0}^{x} q(y)\tilde{\mu}_{22}(y;z,\ep)\, \d y/\ep \,,
\label{e:mu12tildinteqn}
\\
\tilde{\mu}_{22}(x;z,\ep) &= \e^{2\i zx/\ep} -\int_{0}^{x} \e^{2\i z(x-y)/\ep}\overline{q(y)}\tilde{\mu}_{12}(y;z,\ep)\, \d y/\ep \,.
\label{e:mu22tildinteqn}
\end{align}
\ese
Eliminating $\tilde{\mu}_{12}$, one then obtains
\bse
\begin{gather}
\tilde{\mu}_{22}(x;z,\ep) = \e^{2\i zx/\ep} + G(\tilde{\mu}_{22})(x;z,\ep)\,,
\label{e:mu22tilde}
\\
\noalign{where}
G(f)(x;z,\ep) := \int_{0}^{x} g(x,t;z,\ep)f(t)\, \d t/\ep\,,
\\
g(x,t;z,\ep) := -q(t)\int_{t}^{x}\e^{2\i z(x-y)/\ep}\overline{q(y)}\, \d y/\ep\,, \hspace{3mm} x>t\,.
\label{e:kern2}
\end{gather}
\ese
Note that, $\forall x\in[0,L]$ and $\forall t\in[0,x]$, $|g(x,t;z,\ep)|\leq L\n q \n_{\infty}^{2}/\ep$.
Then, 
introducing the Neumann series
\bse
\begin{gather}
\tilde{\mu}_{22}(x;z,\ep) = \sum_{n=0}^\infty\tilde{\mu}_{22}^{(n)}(x;z,\ep)\,,
\\
\tilde{\mu}_{22}^{(0)} = \e^{2\i zx/\ep},\qquad \tilde{\mu}_{22}^{(n+1)}(x;z,\ep) = G(\tilde{\mu}_{22}^{(n)}) = G^{n+1}(\tilde{\mu}_{22}^{(0)})\,,
\end{gather}
\ese 
we have the bound
\begin{align*}
|G^{n}(f)(x;z,\ep)| \leq 
\int_{0}^{x}|g(x,t_1;z,\ep)| 
\cdots \int_{0}^{t_{n-1}}|g(t_{n-1},t_n;z,\ep)||f(t_n)| \d t_n \cdots \d t_1/\ep^n \leq
\frac{M^n}{n!} \|f\|_{\infty}\,,
\end{align*}
where $M:=(L/\ep)^2\|q\|_{\infty}^{2}$.  
For $\Im z \ge 0$ it follows $|\e^{2\i zx/\ep}|\leq 1$ so $\tilde{\mu}_{22}(x;z,\ep) = \sum_{n=0}^{\infty} G^{n}(\e^{2\i zx/\ep})$ is absolutely convergent, and 
\be
\sup_{x \in [0,L], \Im z \geq 0} |\tilde{\mu}_{22}(x;z,\ep)| \leq \sum_{n=0}^{\infty} |G^{n}(\e^{2\i zx/\ep})(x;z,\ep)| \leq \e^{M}\,.
\label{e:mu22tildebound}
\ee
From~\eqref{e:mu12tildinteqn} and~\eqref{e:mu22tildebound} we have $|\tilde{\mu}_{12}(x;z,\ep)| \leq \e^{M}L\n q \n_{\infty}/\ep$.
Using~\eqref{e:mu22tildinteqn} one gets~\eqref{e:mu22tildasymp_imz>0}.

Next we consider the case $z\to \infty$ while $\Im z$ remains bounded.
The Riemann-Lebesgue lemma applied to~\eqref{e:kern2} implies that,
$\forall x \in [0,L]$ and $\forall t \in [0,x]$,
$g(x,t;z,\ep) = o(1)$ as $z \to \infty$ from $\Im z \ge 0$.
Together with the uniform boundedness of $\tilde{\mu}_{22}$ for $\Im z \ge 0$, 
this gives
\be
|G(\tilde{\mu}_{22})(x;z,\ep)| \leq \e^{M} \int_{0}^{x}|g(x,t;z,\ep)|\,\d t = o(1)\,, \hspace{2mm} z \to \infty\,, \quad \Im z \ge 0\,.
\label{e:kernest2}
\ee
Thus, 
$\tilde{\mu}_{22}(x;z,\ep) = \e^{2\i zx/\ep} +\, o(1)$ as $z \to \infty$ from $\Im z \ge 0$,
i.e.,~\eqref{e:tildmu22asymp}. 
Inserting this expression into~\eqref{e:mu12tildinteqn} then yields, by the Riemann-Lebesgue lemma,
\be
\tilde{\mu}_{12}(x;z,\ep) = \int_{0}^{x} q(y)(\e^{2\i zy/\ep} + o(1))\, \d y/\ep = o(1)\,, \hspace{2mm} z \to \infty\,, \quad \Im z \ge 0\,.
\label{e:asympest}
\ee
Now assume additionally that $q' \in L^1([0,L])$.
From~\eqref{e:ode_mu} it follows $\tilde{\mu}_{12}' = q(x)\tilde{\mu}_{22}/\ep$.
Then integrate~\eqref{e:mu22tildinteqn} by parts:
\be
\tilde{\mu}_{22}(x;z,\ep) = \e^{2\i zx/\ep} + \frac{1}{2\i z} \overline{q(x)}\tilde{\mu}_{12}(x;z,\ep) - \frac{1}{2\i z} \int_{0}^{x} \e^{2\i z(x-y)/\ep}(\overline{q(y)}\tilde{\mu}_{12}(x;z,\ep))' \d y\,,
\label{e:parts3}
\ee
where $\tilde{\mu}_{12}(0;z,\ep) \equiv 0$. 
Expanding the derivative, 
using~\eqref{e:mu12tildinteqn},
and $q'\in L^{\infty}([0,L])$ gives
\be
\tilde{\mu}_{22}(x;z,\ep) = \e^{2\i zx/\ep} - \frac{1}{2\i z} \int_{0}^{x} \e^{2\i z(x-y)/\ep}\overline{q(y)}\tilde{\mu}_{12}'(y;z,\ep)\, \d y + O(1/(\Im z)^2)\,,
\label{e:parts4}
\ee
as $\Im z \to \infty$ from $\Im z > 0$. 
Then using $\tilde{\mu}_{12}' = q(x)\big(\e^{2\i zx/\ep} + O(1/(\Im z)\big)$ as $\Im z \to \infty$ gives~\eqref{e:tildemu22asymp2star}. 
Moreover, 
$\tilde{\mu}_{12}' = q(x)\big(\e^{2\i zx/\ep} + o(1)\big)$ as $z \to \infty$ from $\Im z \ge 0$.
Hence, 
using~\eqref{e:parts3} one gets 
\be
\tilde{\mu}_{22}(x;z,\ep) = \e^{2\i zx/\ep} +
\frac{1}{2\i z} \overline{q(x)}\tilde{\mu}_{12}(x;z,\ep) - \frac{\e^{2\i zx/\ep}}{2\i z} \int_{0}^{x}|q(y)|^2\, \d y/\ep + o\Big(\frac{1}{z}\Big)\,,
\label{e:tildmu22parts}
\ee
as $z \to \infty$ from $\Im z \ge 0$.
Using~\eqref{e:asympest} gives~\eqref{e:tildmu22asymp2}.  
\qed

We are now ready to prove Lemma~\ref{l:asympDelta}.

\textit{Proof of Lemma~\ref{l:asympDelta}.}
Now that we have the asymptotic behavior of $\mu_{11}$, 
and $\tilde{\mu}_{22}$, 
Lemma~\ref{l:asympDelta} follows by expressing $\D$ as 
\begin{align*}
\D(z) &= \half \tr \Phi(L;z,\ep) = \half \e^{-\i zL/\ep}(\mu_{11}(L;z,\ep) + \tilde{\mu}_{22}(L;z,\ep))\,,
\end{align*}
and combining the above results using the relevant function spaces.
\qed

\textit{Proof of Lemma~\ref{l:Dprime}.}
One of the consequences of Lemmas~\ref{l:asympmu} and~\ref{l:mutildasym} is that the normalized matrix fundamental solution $\Phi(x;z,\ep)$, 
introduced right above equation~\eqref{e:monodromy}, 
has the asymptotics
\be
\label{ass-Phi}
\Phi(x;z,\ep)
= (\I+o(1))\e^{-\i zx\s_3/\ep}\,, \quad {\rm as} 
\hspace{2mm} z\to\infty\,, 
\ee
from $\Im z \ge 0$ and bounded,
and for each fixed $\ep>0$.

To prove  the remaining formula \eqref{D'-ass} we start with the equation
\be\label{e:Phi_z}
\partial_z\Phi (x;z,\ep)= -\i \Phi(x;z,\ep) \int _0^x   \Phi^{-1}(y;z,\ep)  \sigma_3 \Phi(y;z,\ep) \d y,
\ee
obtained by differentiating~\eqref{e:zs} (with $\bf v$ replaced by $\Phi$) with respect to $z$  and then solving the resulting nonhomogeneous ODE for $\partial_z\Phi$. 
Equation \eqref{D'-ass} follows after substituting~\eqref{ass-Phi} into~\eqref{e:Phi_z}.
A further consequence of Lemmas~\ref{l:asympmu} and~\ref{l:mutildasym} is that if 
$q' \in L^1([0,L])$, 
then the normalized matrix fundamental solution $\Phi(x;z,\ep)$ has the asymptotics
\begin{multline}
\Phi(x;z,\ep) = \bigg[ \,\mathbb{I} + \frac{1}{2\i z}\sigma_3\int_{0}^{x}|q(y)|^2\, \d y/\ep
+ \frac{1}{2\i z}\sigma_3 ( Q(x) - \e^{-2\i zx \sigma_3/\ep}Q(0) )
+ o\Big(\frac{1}{z}\Big)\, \bigg]\,\e^{-\i zx\sigma_3/\ep}\,,
\label{e:dprime2}
\end{multline}
as $z \to \infty$ from $\Im z \ge 0$ and bounded,
and for each fixed $\ep>0$.
Then substituting~\eqref{e:dprime2} into~\eqref{e:Phi_z} gives
\begin{multline}
\partial_z\Phi(x;z,\ep) = \Big(-\i\frac{x}{\ep}\sigma_3 - \frac{x}{2\ep^2 z}\int_{0}^{x}|q(y)|^2 \d y \Big)\e^{-\i zx\sigma_3/\ep} +
\\ 
\frac{x}{2\ep z} \Big( \e^{-\i zx\sigma_3/\ep} Q(0) + \e^{\i zx\sigma_3/\ep}Q(x) \Big) +  o\bigg(\frac{\e^{L\Im z}}{z}\bigg) \,,
\label{e:phiprime}
\end{multline}
as $z \to \infty$ from $\Im z \ge 0$. 
Finally,~\eqref{e:D'-ass2} is given by $\D'(z) = \half\tr\partial_z\Phi(L;z,\ep)$. 
\qed 

\textit{Proof of Corollary~\ref{l:essentialsing}.} 
From Lemmas~\ref{l:entire} and~\ref{l:asympDelta} we know $\D$ is an entire function differing from a constant. 
Thus $\D$ must have either a pole, 
or an essential singularity at infinity. 
Suppose $\D$ has a pole at infinity. 
Then it must be a polynomial. 
This implies there exists an integer $n\ge1$ such that $\D(z) = O(z^{n})$ as $z \to \infty$. 
By Lemma~\ref{l:asympDelta} one has $\D(z) = O(\e^{L\Im z})$ as $\Im z \to \infty$ which is a contradiction. 
Thus $\D$ must have an essential singularity at infinity. 
\qed

\textit{Proof of Lemma~\ref{l:discretespec}:}
Fix $\nu_o \in \Real$.
Define $f_{\ep}(z) := \D(z) - \cos(\nu_o L)$.
By Lemma~\ref{l:entire} $f_\ep$ is an entire function of $z$. 
Further,
using \eqref{e:bloch} 
one gets $\Sigma_{\nu_o} = \{z \in \Complex : f_\ep(z) = 0\}$.
Hence,
$\Sigma_{\nu_o}$ is the set of zeros of an entire function. 
Thus, 
the Floquet spectrum is discrete. 
\qed

\textit{Proof of Lemma~\ref{l:bound1}:} 
Let $z \in \Sigma_{\Lax}$. 
This implies $\@v \in L^{\infty}(\Real)$ is a non-trivial solution of~\eqref{e:zs} bounded $\forall x\in\R$. 
By Floquet's theorem
$\@v = \e^{\i\nu x}\@w$,
where $\@w(x+L;z,\epsilon) = \@w(x;z,\epsilon)$, 
and $\nu \in \Real$.
Plugging this expression for $\@v$ into \eqref{e:zs} gives a modified scattering problem, 
namely,
\be
\i\ep \@w' = \begin{pmatrix} z + \ep \nu & \ \i q(x) \\ -\i \overline{q(x)} & -z + \ep \nu \end{pmatrix} \@w\,.
\label{e:zs2}
\ee
Write~\eqref{e:zs2} in component form:  
\bse
\label{e:zss}
\begin{align}
\i \ep w'_{1} - \i q(x)w_{2} &= (z + \ep \nu)w_{1}\,,
\label{e:zs21}
\\
\i \ep w'_{2} + \i \overline{q(x)}w_{1} &= (-z + \ep \nu)w_{2}\,.
\label{e:zs22}
\end{align}
\ese
Multiply~\eqref{e:zs21} by $\overline{w}_{1}$ and take the complex conjugate. 
This gives two equations which we integrate over a full period.
Thus we arrive at the integrated expressions 
\bse
\begin{align}
\i\int_{0}^{L} q(x)\overline{w}_{1}w_2\, \d x &= -\i\ep\langle w_1, w_1' \rangle - (z+\ep\nu) \langle w_1,w_1\rangle \,, 
\label{e:starstar}
\\
\i \int_{0}^{L} \overline{q(x)}w_{1}\overline{w}_{2}\, \d x &= \i\ep \langle w_1, w_1' \rangle + (\overline{z}+\ep\nu) \langle w_1,w_1 \rangle \,,
\label{e:star}
\end{align}
\ese
where $\langle \cdot,\,\cdot\rangle$ is the $L^2([0,L])$ inner product of a scalar function (cf.~\ref{s:notations}),
and boundary terms vanish since $w_{1}(x+L;z,\ep)=w_{1}(x;z,\ep)$.
Adding~\eqref{e:starstar} to~\eqref{e:star} one gets 
\be
-\Im z \| w_{1}\|_{2}^{2} =  
\Re \langle q w_{2}, w_{1} \rangle\,.
\label{e:eq1}
\ee
Similarly, 
multiply \eqref{e:zs22} by $\overline{w}_{2}$, 
take the complex conjugate,
and integrate.
Thus we have the integrated expressions
\bse
\begin{align}
\i\int_{0}^{L} \overline{q(x)}w_{1}\overline{w}_{2}\, \d x &= \i\ep\langle w_2, w_2' \rangle + (-z+\ep\nu) \langle w_2,w_2\rangle\,,
\label{e:bstarstar}
\\
\i \int_{0}^{L} q(x)\overline{w}_{1}{w}_{2}\, \d x &= -\i\ep\langle w_2, w_2' \rangle + (\overline{z}-\ep\nu) \langle w_2,w_2\rangle \,,
\label{e:bstar}
\end{align}
\ese
where again boundary terms vanish since $w_{2}(x+L;z,\ep) = w_{2}(x;z,\ep)$. 
Adding~\eqref{e:bstarstar} and~\eqref{e:bstar} one gets 
\be
-\Im z \| w_{2}\|_{2}^{2} = 
\Re \langle qw_2, w_1 \rangle\,.
\label{e:eq2}
\ee
Equating~\eqref{e:eq1} and~\eqref{e:eq2} we conclude:
\be
z\in\Sigma_{\Lax}\setminus\R \implies \|w_1\|_{2} = \|w_2\|_{2}\,.
\label{e:eqnorm}
\ee
Thus, 
the Cauchy-Schwarz inequality implies
\begin{align*}
0 < |\Im z|\| w_2\|_{2}^{2} \leq |\langle qw_2, w_1 \rangle|\,
\leq \|q\|_{\infty} \| w_2\|_{2}^{2}\,.
\end{align*}
Hence, 
one gets~\eqref{e:inequality1}. 

Next, 
let $\zeta\in\Sigma_{\Dir}(x_o)$.
Without loss of generality take $x_o=0$.
Note that the Dirichlet boundary conditions~\eqref{e:Dirbcs} imply:
\be
|v_1(L;\zeta,\ep)|^2 - |v_1(0;\zeta,\ep)|^2 = |v_2(L;\zeta,\ep)|^2 - |v_2(0;\zeta,\ep)|^2\,.
\label{e:bcs}
\ee
Write~\eqref{e:zs} in component form:
\bse
\begin{align}
\i \ep v'_{1} - \i q(x)v_{2} &= \zeta v_{1}\,,
\label{e:zs11}
\\
\i \ep v'_{2} + \i \overline{q(x)}v_{1} &= -\zeta v_{2}\,.
\label{e:zs12}
\end{align}
\ese
Multiply~\eqref{e:zs11} by $\overline{v}_1$ and take the complex conjugate. This gives two equations which we integrate over a full period. 
Thus we arrive at the integrated expressions 
\bse
\begin{align}
\i\int_{0}^{L} q(x)\overline{v}_{1}v_2\, \d x &= \i\ep(|v_1(L;\zeta,\ep)|^2 - |v_1(0;\zeta,\ep)|^2)
-\i\ep \langle v_1, v_1' \rangle - \zeta \| v_1\|_2^2\,,
\label{e:dir11}
\\
\i \int_{0}^{L} \overline{q(x)}v_{1}\overline{v}_{2}\, \d x &= \i\ep \langle v_1, v_1' \rangle + \overline{\zeta} \| v_1\|_2^2 \,.
\label{e:dir12}
\end{align}
\ese
Adding~\eqref{e:dir11} to~\eqref{e:dir12} one gets
\be
\int_0^L q(x)\overline{v}_1v_2 + \overline{q(x)}v_1\overline{v}_2\, \d x = \ep(|v_1(L;\zeta,\ep)|^2 - |v_1(0;\zeta,\ep)|^2) - 2\Im \zeta \|v_1\|_2^2\,.
\label{e:dirstar}
\ee
Similarly, 
multiply~\eqref{e:zs12} by $\overline{v}_2$ and take the complex conjugate. 
This gives two equations which we integrate over a full period. 
Then add to get
\be
\int_0^L q(x)\overline{v}_1v_2 + \overline{q(x)}v_1\overline{v}_2\, \d x = -\ep(|v_2(L;\zeta,\ep)|^2 - |v_2(0;\zeta,\ep)|^2) - 2\Im \zeta \| v_2\|_2^2\,.
\label{e:dirstar2}
\ee
Adding~\eqref{e:dirstar} to~\eqref{e:dirstar2} gives
\be
-\Im \zeta \langle \@v,\@v\rangle = \int_0^L q(x)\overline{v}_1v_2 + \overline{q(x)}v_1\overline{v}_2\, \d x\,,
\label{e:dirstar3}
\ee
where $\langle \@v,\@v\rangle$ is the $L^2([0,L],\C^2)$ inner product of a two-component vector function (cf.~\ref{s:notations}).
Thus, 
since $2|v_1||v_2| \leq |v_1|^2 + |v_2|^2$ it follows
\begin{align*}
|\Im \zeta| \|\@v\|_2^2 \leq \|q\|_{\infty} \int_{0}^{L} |\overline{v}_1v_2| + |v_1\overline{v}_2|\, \d x\,
\leq \|q\|_{\infty}\|\@v\|_{2}^2\,.
\end{align*}
Hence, 
one gets~\eqref{e:dirbound1}.
\qed

In the following proofs we use the fact that one can rewrite~\eqref{e:Delta_asymp_1} as:
\be
\D(z) = \cos(zL/\ep) + o(1)\,, \quad z \to \infty\,, \hspace{2mm} 0 \le \Im z \le \n q\n_{\infty}\,,
\label{e:proofasymp}
\ee
for each fixed $\ep>0$,
and that $\D(\overline{z})=\overline{\D(z)}$.
\begin{lemma}
\label{lem-rouche}
Let $R_{n}$ denote the rectangular region with vertices $\pm \frac{n\pi\ep} L \pm {\rm i} \|q\|_{\infty}$,  
where $n\in \N$, 
and $\ep>0$ is fixed. 
Define $R_{k,n} := R_{k}\setminus R_{n}$ where $k\in\N$ with $k>n$.
Then for $k,n$ sufficiently large the functions $\D(z)$ and $\cos(zL/\ep)$ 
have the same number of zeros inside  $R_{k,n}$.
\end{lemma}
\textit{Proof of Lemma~\ref{lem-rouche}.}
Let us fix some arbitrary $\delt\in(0,1)$ such that $\delt< \sinh(\|q\|_{\infty})$ and assume that 
$k,n\in \N$ are so large that,
according to~\eqref{e:proofasymp}
\be
\label{ineq-1}
|\D(z) - \cos(zL/\ep)|\leq \delt\,,  \quad z\in\partial R_{k,n}\,,
\ee
where $\partial R_{k,n}$ is the boundary of $R_{k,n}$.
Then, 
using $|\cos(x+\i y)|^2=\cos^2x+\sinh^2y$, 
we see that $|\cos(zL/\ep)|\geq 1$, 
and $|\cos(zL/\ep)|\geq \sinh( \|q\|_{\infty})$ on the vertical and the horizontal sides of $\partial R_{k,n}$,
respectively. 
Thus, 
by inequality~\eqref{ineq-1},
\be\label{ineq-2}
0<|\D(z) - \cos(zL/\ep)| < |\cos(zL/\ep)|\,,  \quad z\in\partial R_{k,n}\,.
\ee
We now apply Rouch\'{e}'s theorem (cf.~\cite{titchmarsh}) for $g(z):=\D(z) - \cos(zL/\ep)$
and $f(z):=\cos(zL/\ep)$ to complete the argument.
\qed

\textit{Proof of Lemma~\ref{l:zeros}.}
By Lemma~\ref{l:bound1},
all zeros of $\D(z)$ are confined to $\mathbf{S}:=\{z\in\C:\,|\Im z| \leq \n q\n_{\infty}\}$ (cf.~\eqref{e:laxspec2}).
Take $\delt\in(0,1)$ as in Lemma \ref{lem-rouche} and choose $M>0$ so large
that the inequality~\eqref{ineq-1} holds for all $z\in \mathbf{S}$ such that $|\Re z|>M$. Denote by $\mathbf{S}_M$ the set of such $z$.
To prove the lemma, 
we first prove that $\D(z)$ has infinitely 
many real zeros and then show that it has no complex (non-real) zeros in $\mathbf{S}_M$.

Consider the points $z_n={n\pi \ep}/L$ and $z_{n+1}={(n+1)\pi \ep}/L$, $n\in \N$,
where, 
say, 
$z_n>M$.
The remaining case $z_{-n}<-M$ can be worked out similarly. 
Since $\D(z)$ is real on $\R$,
\eqref{ineq-1} shows that $\D(z_n)\D(z_{n+1})<0$ and, 
thus, 
there exists a zero of $\D(z)$ on $(z_n,z_{n+1})$. 
Thus, 
there are infinitely many zeros of $\D(z)$ on $\R$.

We can now use Lemma~\ref{lem-rouche} to show that there is exactly one zero of $\D(z)$ on
the rectangle $R_{n,n+1}$ and, 
as it was just shown above this zero is real. 
Thus, 
there are no non-real zeros in $\mathbf{S}_M$ and we have completed the proof about the zeros of $\D(z)$. 
In view of~\eqref{D'-ass},
the proof of zeros of $\D'(z)$ is similar.
\qed

\textit{Proof of Theorem~\ref{t:resolvent}.}
Consider the set $\rho:=\mathbf{H}\setminus\Sigma_{\Lax}$, 
where $\mathbf{H}:=\{z \in \C : \Im z > 0 \}$ is the upper half-plane.
Since $\Sigma_{\Lax}$ is Schwarz symmetric, 
and $\Real \subset \Sigma_{\Lax}$ (cf. Lemma~\ref{l:infband}) one gets
$\rho_{\Lax} = \rho \cup \overline{\rho}$.
To prove the result it is only necessary to show that $\rho$ is connected.
First, 
by Lemma~\ref{l:closedband} no spectral band can be closed in the finite $z$-plane. 
Next, 
we show that the only spectral band extending to infinity is the real $z$-axis.
Suppose to the contrary that there exists a band extending to infinity and intersecting the real $z$-axis at most once.
By Lemma~\ref{l:bound1}, 
and~\eqref{e:proofasymp}, 
it is sufficient to assume the band is confined to $\{z\in\C:\Re z \geq 0\,, 0 \le \Im z \le \n q\n_{\infty}\}$.
Recall that along a spectral band one necessarily has $\Im \D \equiv 0$.
Fix $\delta \in (0,1)$ such that $\delta < \sinh(\n q\n_{\infty})$. 
Then there exists $N=N(q;\ep)>0$ such that 
\be
|\Re \D(z) - \cos(\Re z L/\ep)\cosh(\Im z L/\ep)| < \delta\,, 
\label{e:realimag}
\ee
for all $\Re z > N$. 
Let $n\in\N$ be such that $N < n\pi\ep/L$.
As in Lemma~\ref{l:zeros} it then follows that there exists a zero of $\D(z)$ 
for $\Re z \in (n\pi\ep/L, (n+1)\pi\ep/L)$.
Further, 
by assumption $\Im z = 0$ for at most one point on the band.
This implies that $\D$ has infinitely many complex zeros which contradicts Lemma~\ref{l:zeros}.
Thus,
we have shown that the real $z$-axis is the only spectral band extending to infinity. 
Hence, 
the set $\rho$ is connected which completes the proof.
\qed

\textit{Proof of Theorem~\ref{t:infeigenvalues}.}
Recall that $\D(z)$ is real-valued along the real $z$-axis.
By Lemma~\ref{l:zeros} there exists a sequence $\{z_n\}$ of simple real  zeros of $\D$ for $|n|$ sufficiently large.
Suppose that $z_{k-1}$ and $z_{k}$ are two consecutive zeros such that $\D'(z_{k-1}) > 0$, 
and $\D'(z_{k}) < 0$. 
Then by Rolle's theorem there exists $z^{\rm c}_{k} \in (z_{k-1}, z_{k})$ such that 
$\D'(z^{\rm c}_k) = 0$.
Moreover,
$\{z_n^{\rm c}\}$ are simple zeros of $\D'$ by Lemma~\ref{l:zeros}.
Also, 
by Lemma~\ref{l:infband} necessarily $0 < \D(z^{\rm c}_k) \le 1$.
Suppose  $\D(z^{\rm c}_k) < 1$.  
Since $\Re \D$ is harmonic it follows that
$z^{\rm c}_k$ is a saddle point; 
and a band emanates from $z^{\rm c}_k$ into the complex plane along a steepest ascent curve.
Further,
$\Im \D \equiv 0$ along this steepest ascent curve. 
Thus, 
by Theorem~\ref{t:resolvent} and continuity there must exist a band edge along this steepest ascent curve at which one gets $\D(z_k^{+})=1$.
Otherwise $\D(z^{\rm c}_k) = 1$ (i.e., $z^{\rm c}_k$ is a double point) and no further analysis is required.
Next,
$\D'(z_{k}) < 0$ implies $\D'(z_{k+1}) > 0$.
In this case $-1\leq \D(z_{k+1}^{\rm c}) < 0$,
and the argument is completely analogous to that above.
Finally,
by Theorem~\ref{theo-spines} the bands emanating from the real $z$-axis do not intersect; 
and by Lemma~\ref{l:discretespec} the Floquet spectrum is discrete.
This completes the proof.
\qed

\textit{Proof of Theorem~\ref{t:finiteband1}.}
Suppose $q$ is a finite-band potential.
Recall $z\in\Sigma_{\Lax}$ implies $|\Im z| \le \n q\n_{\infty}$ (cf. Lemma~\ref{l:bound1}).
By Theorem~\ref{t:resolvent} the real $z$-axis is the only band extending to infinity.
So, 
there are finitely many bands and each band (except $\Real$) is bounded.
Thus, 
$\Sigma_{\Lax}\setminus \R$ is bounded and the result follows trivially.
Next let $N=N(q;\ep) > 0$ be such that $(\Sigma_{\Lax}\setminus\Real) \subset R_{N}$.
Suppose $\Sigma_{\Lax}$ is comprised of infinitely many bands.
Then infinitely many periodic eigenvalues (which correspond to band edges) exist in the closure of $R_{N}$.
This impies the set of periodic eigenvalues must have a finite limit point which is a contradiction (cf. Lemma~\ref{l:discretespec}).
\qed  

\textit{Proof of Theorem~\ref{t:finiteband2}.}
Recall that the periodic and antiperiodic Floquet eigenvalues correspond to the band edges of $\Sigma_{\Lax}$.
Suppose that there exist infinitely many spectral bands.
Then by Theorem~\ref{t:finiteband1} for any $N>0$ the set $(\Sigma_{\mathrm{Lax}} \setminus \Real) \cap R_{N} \neq \emptyset$.
This implies that there exist infinitely many periodic,
or antiperiodic, 
eigenvalues along the imaginary $z$-axis which is a contradiction. 
The contradiction follows from the fact that the Floquet eigenvalues are discrete with no finite accumulation points, 
and $|\Im z| \leq \n q\n_{\infty}$.
\qed  

\textit{Proof of Theorem~\ref{theo-spines}.}
Let $\g_j$ denote a generic spectral band.  
Two possibilities arise:  Either $\g_j$ has a periodic and an antiperiodic endpoint, 
or both of its endpoints are periodic (or antiperiodic).
If one endpoint is periodic and the other is antiperiodic, $\D(z)=0$ at some point on $\g_j$.
Conversely, 
if both endpoints are periodic or antiperiodic, 
there is a point on $\g_j$, 
where  $\D'(z)=0$. 
In either case, 
according to Lemma~\ref{l:zeros}, 
there can be no more than finitely many such bands that do not intersect $\R$. 
In order to prove that these bands are actually spines, 
it remains to show that, 
as $z\to\infty$, 
these bands cannot intersect any other bands.
So, 
let $\g_j$ be such a band, 
which intersects $\R$ at some point $z_o$.  
Then $\D'(z_o)=0$. 
Again, by Lemma \ref{l:zeros} there can only be finitely many such bands that have more than one 
point where $\D'(z)=0$ (since the number of complex zeros of $\D'$ is finite). 
Hence there can be only finitely many bands that, in addition to  $\R$, also intersects 
some other band.
Thus, we proved that all but finitely many bands are the spines.
Moreover,
by Lemma~\ref{l:Dprime} $\D'(z)=-L\sin(zL/\ep)/\ep+o(1)$ as $z\to\infty$ from $|\Im z|\leq \|q\|_{\infty}$. 
Thus,
spines are $o(1)$-close to $n\pi\ep/L$ for $|n|\in\N$ sufficiently large.
That this is the only $o(1)$-close spine follows from Rouch\'e's theorem.
This completes the proof.
\qed  

\subsection{Proofs: Semiclassical limit} 
\label{s:semiclassicalproofs}

\textit{Proof of Lemma~\ref{l:bound2}.} 
Let $z \in \Sigma_{\Lax}$ be such that $|\Im z| > 0$, and $|\Re z|>0$. 
Again, 
we can write $\@v = \e^{\i \nu x}\@w$, 
where $\@w(x+L;z,\ep) = \@w(x;z,\ep)$, 
and $\nu \in \Real$. 
Plugging this expression for $\@v$ into \eqref{e:zs} gives the modified
ZS system~\eqref{e:zs2}. 
Multiply~\eqref{e:zs21} by $\overline{w}'_{1}$ and integrate by parts
\be
\i\ep \|w_1'\|_2^2 + \i \int_{0}^{L} (q(x)w_2)'\overline{w}_{1}\, \d x = (z+\ep \nu)\langle w_1, w_1' \rangle \,.
\label{e:one}
\ee
The complex conjugate to \eqref{e:one} is
\be
-\i\ep \|w_1'\|_2^2 - \i \int_{0}^{L}(\overline{q(x)} \overline{w}_{2})'{w}_{1}\, \d x = (\overline{z}+\ep\nu) \langle w_1', w_1 \rangle\,.
\label{e:two}
\ee
Integrate the right hand side of~\eqref{e:two} by parts, 
then add to~\eqref{e:one} and multiply by $\ep$ getting
\be
2\i\ep \Im \langle q'w_2, w_1 \rangle =
2\Im z \langle w_1, \ep w_1'\rangle + 
\int_{0}^{L} \big( \overline{q(x)}w_{1}(\ep\overline{w}_{2}') - q(x)\overline{w}_{1}(\ep w_{2}') \big) \d x\,.
\label{e:three}
\ee
Using~\eqref{e:zs21} one gets
\be
2\Im z \langle w_1, \ep w_1' \rangle =
2\Im z \langle w_1, qw_2 \rangle + 2\i \Im z(\overline{z} + \ep \nu)\|w_2\|_{2}^{2}\,.
\label{e:rhs1}
\ee
Further, from~\eqref{e:zs22} and recalling equations~\eqref{e:eq1}, and~\eqref{e:eq2}, one gets 
\be
\int_{0}^{L} \big( \overline{q(x)}w_{1}(\ep\overline{w}_{2}') - q(x)\overline{w}_{1}(\ep w_{2}') \big) \d x\ =
2\i \Im z \Im \langle qw_2, w_1 \rangle - 2\i \Im z(\ep \nu - \Re z)\| w_2\|_{2}^{2}\,.
\label{e:rhs2}
\ee
Substitute~\eqref{e:rhs1} and~\eqref{e:rhs2} into~\eqref{e:three} and simplify to get
\be
4\i \Im z \Re z \| w_{2}\|_{2}^{2} =
2\i\ep \Im \langle q'w_2, w_1 \rangle 
\label{e:four}
\ee
Applying the Cauchy-Schwarz inequality gives 
\begin{align*}
0 < 2|\Re z| |\Im z| \| w_{2} \|_{2}^{2} \leq \ep|\langle q'w_2, w_1 \rangle| \leq \ep \| q' \|_{\infty} \| w_{2} \|_{2}^{2} \,.
\end{align*}
Hence, one gets~\eqref{e:inequality2}. 
\qed

\textit{Proof of Theorem~\ref{t:set1}.}
This result follows immediately from Lemmas~\ref{l:bound1} and~\ref{l:bound2}.
\qed

\textit{Proof of Theorem~\ref{t:Pi&Xi}.}
This result follows from the cardinality results found in~\cite{kapmity, kappeler} together with Lemmas~\ref{l:bound1} and~\ref{l:bound2}.
\qed

\textit{Proof of Corollary~\ref{t:eigvalsorder}.}
Fix $\nu \in \Real$. 
Denote the corresponding countably infinite set of eigenvalues by $\{z_n(\nu)\}_{n\in\Natural}$.
Without loss of generality assume that $\Re z > 0$ and $\Im \ge 0$.
By Lemma~\ref{l:entire} and~\eqref{e:bloch} it follows that 
in any neighborhood of infinity $\Sigma_{\nu}$ is infinite.
Further $\Sigma_{\nu} \cap D(0;r)$ is finite for any $r>0$, 
where $D(0;r) := \{z \in \Complex : |z| \le r\}$.
This implies that infinity is the only accumulation point of $\Sigma_{\nu}$,
and that,  due to Lemma~\ref{l:bound1}, there exists a partially ordered increasing sequence $\{\Re z_{n}(\nu)\}_{n\in\Natural}$
such that $\Re z_{n}(\nu) \to \infty$ as $n \to \infty$. 
Finally, that $\Im z_n(\nu) = O(1/\Re z_n(\nu))$ now follows easily from Lemma~\ref{l:bound2}.
\qed

It now remains to show that Lemmas~\ref{l:bound1} and~\ref{l:bound2} together imply Corollary~\ref{c:mainresult}. 

\textit{Proof of Corollary~\ref{c:mainresult}.} 
Fix $\delta > 0$. 
Consider the $\delta$-neighborhood defined by $N_{\delta}(\Sigma_{\infty})$.
Without loss of generality assume that $\Im z > 0$, 
and that $|\Re z| > 0$. 
Next consider the curve in the spectral plane defined by
$ |\Im z| = \min \big\{\n q\n_{\infty} \hspace{0.5mm}, \hspace{0.5mm} \ep \| q'\|_{\infty} / 2|\Re z| \big\}$ which bounds the set $\Lambda^{\ep}(q)$.
Then it is easily seen that there exists $\ep_{*}>0$  such that
$\Lambda^{\ep_*}(q) \subset N_{\delta}(\Sigma_{\infty})$. 
Thus, 
if $ z \in \Sigma_{\Lax}$,
then $z \in N_{\delta}(\Sigma_{\infty})$ 
whenever $0 < \ep \leq \ep_{*}$. 
Hence, 
$\Sigma_{\Lax} \cap (\Complex \setminus  N_{\delta}(\Sigma_{\infty})) = \emptyset$ which completes the proof.
\qed

\subsection{Proofs: Real or symmetric potentials}
\label{s:symmetricproofs}

\textit{Proof of Lemma~\ref{l:symmetries}.}
Let $\@v(x;z,\ep)$ be a bounded solution of the ZS system~\eqref{e:zs}.
Then, 
\eqref{e:sym1} establishes the $z \mapsto \overline{z}$ symmetry of the spectrum.  
Similarly, 
\eqref{e:sym2} and~\eqref{e:sym3} establish the $z \mapsto -\overline{z}$ symmetry in the spectrum when the potential is real, 
even, 
or odd.
Suppose that the potential is real. 
By~\eqref{e:sym2} it follows that the monodromy matrix satisfies the symmetry~\eqref{e:msym2}.
Then noting $z=-\overline{z}$ if and only if $z = \i\eta$ we have $\tr M(z;\epsilon)$ is real-valued for $z \in \i\Real$. 
Similarly, 
using~\eqref{e:sym3} one gets the symmetry~\eqref{e:msym3}.
Thus, 
if the potential is even ($\theta = 0 \; \text{mod}\;\pi$), 
or odd ($\theta = \pi/2 \; \text{mod}\;\pi$),
then we have $\tr M(z;\ep)$ is real-valued for $z \in \i\R$.
Finally, 
recalling that $\D(z) = \tr M(z;\ep)/2$ it follows that
when $q$ is real,
even, 
or odd one gets $\i\R \subset \Gamma$.   
\qed

\textit{Proof of Theorem~\ref{t:reimag}.}
Let $z_n^{\pm}$ be a periodic, 
or antiperiodic eigenvalue,
respectively. 
This implies $\Im \D(z_n^{\pm}) = 0$, and $|\Re \D(z_n^{\pm})| = 1$.
By assumption $z_n^{\pm}\in\Real \cup \i\Real$.
Further, 
$z_n^{\pm}$ belongs to $\Sigma_{\Lax}$.
Suppose $z_n^{\pm}\in\gamma_j$ where $\gamma_j$ is a band that leaves the real, 
or imaginary axis.
By analyticity of $\D(z)$ we need consider only two cases.
First, 
suppose there exists another point $z_k^{\pm} \neq z_n^{\pm}$ along the spectral band such that 
$|\Re \D(z_k^{\pm})| = 1$. 
Then $z_k^{\pm} \in \R \cup \i\R$. 
Since the real and imaginary $z$-axes are contours such that $\Im \D(z) = 0$ we get a closed curve in the set $\Gamma$ 
which is a contradiction (see Lemmas~\ref{l:closedband} and~\ref{l:symmetries}).
Second, 
suppose $|\Re\D(z)|<1$ for all $z\in\gamma_j$. 
It then follows that $\gamma_j$ extends to infinity.
By Theorem~\ref{t:resolvent} the only element of $\Sigma_{\mathrm{Lax}}$ which has this property is the real $z$-axis. 
Again, 
we have a contradiction. 
Hence, 
$\Sigma_{\mathrm{Lax}} \subset \Sigma_{\infty} \subset \R \cup \i\R$. 
\qed

\textit{Proof of Lemma~\ref{l:dirbounds2}.} 
Let $\zeta \in \Sigma_{\Dir}(0)$. 
Without loss of generality assume $\Im\zeta > 0$, 
and $\Re\zeta > 0$. 
Recall that the Dirichlet boundary conditions~\eqref{e:Dirbcs} imply:
\be
|v_{1}(L;\zeta,\ep)|^2 - |v_{1}(0;\zeta,\ep)|^2 = |v_{2}(L;\zeta,\ep)|^2 - |v_{2}(0;\zeta,\ep)|^2\,.
\label{e:dirbcs}
\ee
Further, 
the hypotheses imply $\Im q(0) = 0$.
Write~\eqref{e:zs} in component form,
\bse
\label{e:zscompDir1}
\begin{align}
\i\ep v_1' - \i q(x)v_2 &= \zeta v_1\,,
\label{e:zscomp1}
\\
\i \ep v_2' + \i\overline{q(x)}v_1 &= -\zeta v_2\,.
\label{e:zscomp2}
\end{align}
\ese
Multiply~\eqref{e:zscomp1} by $\overline{v}_{1}'$ and integrate by parts
\be
\i\ep \|v_1'\|_2^2 - \i q(0)\big(\overline{v_1(L)}v_2(L) - \overline{v_1(0)}v_2(0)\big) + \i\int_{0}^{L} (q(x)v_{2})'\overline{v}_{1}\, \d x = \zeta \langle v_1, v_1'\rangle\,.
\label{e:form1}
\ee
Then using the boundary conditions one can write 
\be
\i\ep\|v_1'\|_2^2 + \i q(0)\big(|v_{1}(L)|^2 - |v_{1}(0)|^2 \big) + \i\int_{0}^{L} (q(x)v_{2})'\overline{v}_{1}\, \d x = \zeta\langle v_1, v_1'\rangle\,.
\label{e:form2}
\ee
Take the complex conjugate and add to~\eqref{e:form1}
\be
\i\int_{0}^{L} (q(x)v_{2})'\overline{v}_1 -   (\overline{q(x)}\overline{v}_{2})'v_{1} \, \d x = \zeta\langle v_1,v_1'\rangle + \overline{\zeta} \langle v_1',v_1 \rangle\,.
\label{e:form3}
\ee
Note that
\be
\overline{\zeta} \langle v_1',v_1\rangle = \overline{\zeta}\big( |v_{1}(L)|^{2} - |v_{1}(0)|^{2} \big) - \overline{\zeta}\langle v_1,v_1'\rangle\,.
\label{e:side}
\ee
Let $\alpha(\zeta) := \i\ep\overline{\zeta}(|v_1(L)|^2-|v_1(0)|^2)$. Then using~\eqref{e:form3} and~\eqref{e:side} one gets
\be
\alpha(\zeta) + \ep\int_{0}^{L} q'(x)\overline{v}_{1}v_{2} - \overline{q'(x)}v_{1}\overline{v}_{2}\, \d x = (2\ep\Im \zeta) \langle v_1,v_1'\rangle + \ep\int_{0}^{L} \overline{q(x)}v_{1}\overline{v}_{2}' - q(x)\overline{v}_{1}v_{2}'\, \d x \,,
\label{e:form4}
\ee
where we multiplied through by $\ep$. Using~\eqref{e:zscomp2} one then gets
\be
\ep\int_{0}^{L} \overline{q(x)}v_{1}\overline{v}_{2}' - q(x)\overline{v}_1v_2' \, \d x =
(-2\i\Re \zeta)\Re \langle qv_2,v_1 \rangle + (2\i\Im \zeta)\Im \langle qv_2, v_1 \rangle\,.  
\label{e:form5}
\ee
Then using~\eqref{e:dirstar} we get
\be
(-2\i\Re \zeta)\Re\langle qv_2,v_1 \rangle = -\i\Re \zeta \big(\ep (|v_{1}(L)|^{2} - |v_{1}(0)|^{2}) - (2\Im \zeta)\| v_1\|_2^2 \big)\,.
\label{e:form6}
\ee
Further, by~\eqref{e:zscomp1}
\be
(2\ep\Im \zeta)\langle v_1,v_1' \rangle = (2\Im \zeta)\int_{0}^{L} \overline{q(x)}v_{1}\overline{v}_{2}\, \d x + 2\big(\i\Re\zeta\Im \zeta + (\Im \zeta)^2 \big)\| v_1\|_2^2\,.
\label{e:form7}
\ee
Finally, 
using~\eqref{e:form5} and~\eqref{e:form7} and simplifying gives 
\be
\alpha(\zeta) + (2\i\ep) \Im\langle q'v_2,v_1 \rangle = -\i\ep \zeta\big(|v_{1}(L)|^2 - |v_{1}(0)|^2\big) + (4\i\Re \zeta \Im \zeta)\| v_1\|_2^2\,.
\label{e:form8}
\ee
Next we complete the same series of calculations using~\eqref{e:zscomp2}. First, multiply~\eqref{e:zscomp2} by $\overline{v}_{2}'$ and integrate by parts. 
\be
\i\ep \| v_2'\|_2^2 + \i \overline{q(0)}\big(v_1(L)\overline{v_2(L)} - v_1(0)\overline{v_2(0)}\big) - \i\int_{0}^{L} (\overline{q(x)}v_{1})'\overline{v}_{2}\, \d x = -\zeta \langle v_2,v_2' \rangle\,.
\label{e:form9}
\ee
Then using the boundary conditions one can write 
\be
\i\ep \| v_2\|_2^2 - \i \overline{q(0)}\big(|v_2(L)|^2 - |v_2(0)|^2 \big) - \i\int_{0}^{L} (\overline{q(x)}v_{1})'\overline{v}_{2}\, \d x = -\zeta \langle v_2,v_2' \rangle\,. 
\label{e:form10}
\ee
Take the complex conjugate and add
\be
\i\int_{0}^{L} (q(x)\overline{v}_{1})'v_2 -   (\overline{q(x)}v_{1})'\overline{v}_{2} \, \d x = -\zeta \langle v_2,v_2' \rangle - \overline{\zeta} \langle v_2',v_2 \rangle\,.
\label{e:form11}
\ee
Note that
\be
-\zeta \langle v_2,v_2' \rangle  = -\zeta \big( |v_{2}(L)|^{2} - |v_{2}(0)|^{2} \big) + \zeta \langle v_2',v_2 \rangle\,.
\label{e:side2}
\ee
Let $\beta(\zeta) := -\i\ep\zeta(|v_2(L)|^2-|v_2(0)|^2)$. 
Then using~\eqref{e:form11} and~\eqref{e:side2} one gets
\be
\beta(\zeta) + \ep\int_{0}^{L} q'(x)\overline{v}_{1}v_{2} - \overline{q'(x)}v_{1}\overline{v}_{2}\, \d x = (2\ep\Im \zeta) \langle v_2',v_2 \rangle + \ep\int_{0}^{L} \overline{q(x)}v_{1}'\overline{v}_{2} - q(x)\overline{v}_{1}'v_{2}\, \d x\,,
\label{e:form12}
\ee
where we multiplied through by $\ep$. Using~\eqref{e:zscomp1} and~\eqref{e:dirstar2} one then gets
\be
(2\i)\Im \langle v_1',qv_2 \rangle  =
\i\Re \zeta\big(\ep(|v_{1}(L)|^2 - |v_{1}(0)|^2) + (2\Im \zeta)\n v_2\n_2^2 \big)
+ (2\i \Im \zeta) \Im \langle v_1,qv_2 \rangle\,.
\label{e:form13}
\ee
Further, by~\eqref{e:zscomp2}
\be
(2\ep\Im \zeta) \langle v_2',v_2 \rangle = (-2\Im z)\langle v_1,qv_2 \rangle + 2\big(\i\Re\zeta \Im \zeta - (\Im \zeta)^2 \big)\n v_2 \n_2^2\,.
\label{e:form14}
\ee
Finally, using~\eqref{e:form13} and~\eqref{e:form14} and simplifying gives 
\be
\beta(\zeta)  + \ep\int_{0}^{L}q'(x)\overline{v}_{1}v_{2} - \overline{q'(x)}v_{1}\overline{v}_{2} \d x = \i\ep\overline{\zeta}\big(|v_{1}(L)|^2 - |v_{1}(0)|^2\big) + (4\i\Re \zeta \Im \zeta)\|v_2\|_2^2\,.
\label{e:form15}
\ee
Add~\eqref{e:form15} to~\eqref{e:form8} to get
\be
2\ep \int_{0}^{L} q'(x)\overline{v}_1v_{2} - \overline{q'(x)}v_1\overline{v}_2 \d x = 4\i\Re\zeta \Im \zeta \int_{0}^{L} |v_{1}|^2 + |v_{2}|^2 \d x\,.
\label{e:form16}
\ee
 Noting $2|v_1||v_2| \leq |v_1|^2 + |v_2|^2$ one gets
\begin{align*}
4|\Re \zeta||\Im \zeta| \|\@v\|_{2}^{2} &= 2\ep\big| \int_{0}^{L} q'(x)\overline{v}_1v_{2} - \overline{q'(x)}v_{1}\overline{v}_{2}\, \d x\big|\,
\\
&\leq 2\ep \|q'\|_{\infty} \int_{0}^{L} |\overline{v}_1v_{2}| + |v_{1}\overline{v}_{2}| \d x\, \\
&\leq  2  \ep\|q'\|_{\infty}\|\@v\|_{2}^{2}\,,
\end{align*}
which in turn yields the desired result.
\qed

\textit{Proof of Theorem~\ref{t:dirset2}.}
This result follows from Lemmas~\ref{l:bound1} and~\ref{l:dirbounds2}.
\qed

\textit{Proof of Lemma~\ref{l:wkbbound}.}
We take an approach analogous to the one first employed by Klaus and Shaw for single-lobe potentials decaying as $x\to\pm\infty$. 
Let $z \in \Sigma_{\text{Lax}}\setminus \i\R$.
This implies the solution $\@v \not\equiv 0 \in L^{\infty}(\R)$. 
By Floquet's theorem $\@v = \e^{\i\nu x}\@w$, 
where $\@w(x+L;z,\ep) = \@w(x;z,\ep)$, and $\nu\in\R$. 
Plugging this expression for $\@v(x;z,\ep)$ into the ZS system gives~\eqref{e:zss}
with $\overline{q(x)} = q(x)$. 
Multiply~\eqref{e:zs21} by $\overline{w}_2$, 
and~\eqref{e:zs22} by $\overline{w}_1$.
Then subtract and integrate over the period to get the expression:
\be
\i\ep\int_0^L w_1'\overline{w}_2 - \overline{w}_1w_2'\,\d x - \i \int_{0}^{L} q(x)(|w_1|^2 + |w_2|^2)\, \d x = 2z \Re \langle w_1,w_2 \rangle + 2\i\ep\nu \Im \langle w_1,w_2 \rangle\,.
\label{e:lhsimaginary}
\ee
Note that,
due to periodicity,
integration by parts gives
\be
\int_0^L w_1'\overline{w}_2 - \overline{w}_1w_2'\,\d x = \int_0^L \overline{w}_1'w_2 - w_1\overline{w}_2'\,\d x\,.
\label{e:lhsimaginary2}
\ee
Hence,
\be
\overline{\int_0^L w_1'\overline{w}_2 - \overline{w}_1w_2'\, \d x} = \int_0^L w_1'\overline{w}_2 - \overline{w}_1w_2'\, \d x\,,
\label{e:imaginary}
\ee
and it follows that the left-hand side of~\eqref{e:lhsimaginary} is purely imaginary. 
Thus, 
$(\Re z) \Re \langle w_1, w_2 \rangle = 0$. 
By assumption $\Re z > 0$ which implies $\Re \langle w_1,w_2 \rangle = 0$. 
Next, 
multiply~\eqref{e:zs21} by $\overline{w}_1$ and write as:
\be
\overline{w}_1w_2 = \frac{\ep w_1'\overline{w}_1}{q(x)} + \i(z+\ep\nu)\frac{|w_1|^2}{q(x)}.
\label{e:ks1}
\ee
Take the complex conjugate of~\eqref{e:ks1}, 
add,
and integrate over the period:
\begin{align*}
2 \Re\langle w_1,w_2\rangle &= \ep \int_{0}^{L} \frac{(w_1\overline{w}_1)'}{q(x)}\,\d x - 2\Im z \int_{0}^{L} \frac{|w_1|^2}{q(x)}\,\d x \\
&= \ep \int_{0}^{L} \frac{|w_1|^2q'(x)}{q^2(x)}\,\d x - 2\Im z \int_{0}^{L} \frac{|w_1|^2}{q(x)}\,\d x
\end{align*}
Note that since $q(x) > 0$ it follows easily that
\be
\int_{0}^{L} \frac{|w_1|^2}{q(x)}\,\d x \ge \min\limits_{x\in[0,L]}\Big\{\frac{1}{q(x)}\Big\} \|w_1\|_{2}^{2} > 0\,.
\label{e:ks2}
\ee
Finally, 
since $\Re \langle w_1,w_2 \rangle = 0$ one gets
\be
|\Im z|  \int_{0}^{L} \frac{|w_1|^2}{q(x)}\d x  = \bigg|\frac{\ep}{2} \int_{0}^{L} \frac{|w_1|^2q'(x)}{q^2(x)}\d x \bigg| \leq \frac{\ep}{2}\| (\ln q)' \|_{\infty}\int_{0}^{L} \frac{|w_1|^2}{q(x)}\d x\,,
\ee
which gives~\eqref{e:wkbbound}.

Next,
without loss of generality assume $x_o=0$.
Let $\zeta\in\Sigma_{\rm Dir}(0)\setminus \i\R$.
Write the ZS system in component form as in~\eqref{e:zscompDir1}.
Multiply~\eqref{e:zscomp1} by $\overline{v}_2$,
and~\eqref{e:zscomp2} by $\overline{v}_1$. 
Subtract and integrate over the period to get the expression:
\be
\i\ep \int_0^L v_1'\overline{v}_2 - \overline{v}_1v_2'\, \d x - \i\int_{0}^{L} q(x)(|v_1|^2+|v_2|^2)\,\d x = 2\zeta\Re\langle v_1,v_2 \rangle\,.
\label{e:integratedDir}
\ee
Then note integration by parts gives
\be
\int_0^L v_1'\overline{v}_2 - \overline{v}_1v_2'\, \d x = \int_0^L\overline{v}_1'v_2 - v_1\overline{v}_2'\, \d x \,,
\label{e:lhsimaginaryDir}
\ee
where since $\@v(x;z,\ep)$ is an eigenfunction corresponding to Dirichlet BCs~\eqref{e:Dirbcs} it follows
\be
(v_1\overline{v}_2-\overline{v}_1v_2)|_{0}^{L} = -|v_2(L)|^2 + |v_1(L)|^2 + |v_2(0)|^2 - |v_1(0)|^2 = 0\,.
\ee
Hence,
\be
\overline{\int_0^L v_1'\overline{v}_2 - \overline{v}_1v_2'\, \d x} = \int_0^L v_1'\overline{v}_2 - \overline{v}_1v_2'\, \d x\,,
\ee
and it follows that the left-hand side of~\eqref{e:integratedDir} is purely imaginary.
The rest of the argument is identical to that above for the Lax spectrum.
This completes the proof.
\qed

\section{Examples}
\label{s:examples}

In this section we illustrate the results of the previous sections by considering certain interesting classes of periodic potentials and 
by computing their Lax spectrum analytically, 
or numerically using Floquet-Hill's method~\cite{DK2006}.
(We refer the reader to~\ref{s:fhm} for some details on Floquet-Hill's method.)
All results were checked for numerical convergence. 

\subsection{Plane wave potentials}
\label{s:const}
As our first example we consider a potential whose spectrum can be computed exactly (see~\ref{s:planewave} for details). 
Namely, 
we consider the plane wave potential
\begin{equation}
q(x) = A\e^{\i Vx}\,,
\label{e:plane_wave}
\end{equation}
where $A\in\Real$ is constant amplitude, 
and $V\in\Real$ is the wave number. 

When $V=0$ we get a constant background.
In this case:
\be
\Sigma_{\Lax} = \R \cup \i[-A, A]\,.
\label{e:specset1}
\ee
When $V \ne 0$, and $L= 2\pi/|V|$ the Lax spectrum is given by:
\begin{equation}
\Sigma_{\Lax} = \Real \cup \big[-\ep V/2 - \i A \;,\; -\ep V/2 + \i A \big]\,.
\label{e:spec_plnwave}
\end{equation}
Thus, 
the Lax spectrum is composed of two bands in the complex plane, 
as shown in Fig.~\ref{f:4} below. 
Note that in this case the number of bands is not proportional to $1/\ep$ as $\ep\downarrow 0$.
Indeed, 
there are only two spectral bands for any $\ep > 0$. 
Also, 
for any $\ep>0$ and $V\ne 0$ there are no spectral bands on the imaginary axis. 
It is only in the limit $\ep \downarrow 0$ that the complex band becomes purely imaginary. 

\subsection{Piecewise continuous potentials}
As a second example, 
we look at whether the assumptions in Theorem~\ref{c:mainresult} are necessary or merely sufficient. 
To this end, 
suppose that the potential is only piecewise smooth. 
That is $q$ and $q'$ have finite left and right limits for all $x$,
and finitely many jump discontinuities on any bounded interval.
Specifically, 
we consider the $2$-periodic extension of the signum function
\be
q(x) = \text{sgn}(x) := \begin{cases} +1 &\text{if} \;\;\; x > 0 \\ -1 &\text{if} \;\;\; x < 0 \end{cases} 
\label{e:signum}
\ee
where $x \in [-1, 1)$.                                
In this case the potential has a jump discontinuity, 
and therefore the hypotheses of Lemma~\ref{l:bound2} are not satisfied.

On the other hand, 
the analytically calculated Floquet discriminant is given by~\eqref{e:discrimpiece},
and the numerically computed Lax spectrum of the ZS system~\eqref{e:zs} for the potential~\eqref{e:signum} is shown in Fig.~\ref{f:5} below.
In this example, 
the Floquet eigenvalues arise in symmetric quartets due to the symmetry of the potential \eqref{e:signum} (specifically, odd and real).    

\begin{figure}[t!]
\vglue-2ex
\centerline{\includegraphics[width=7.5cm]{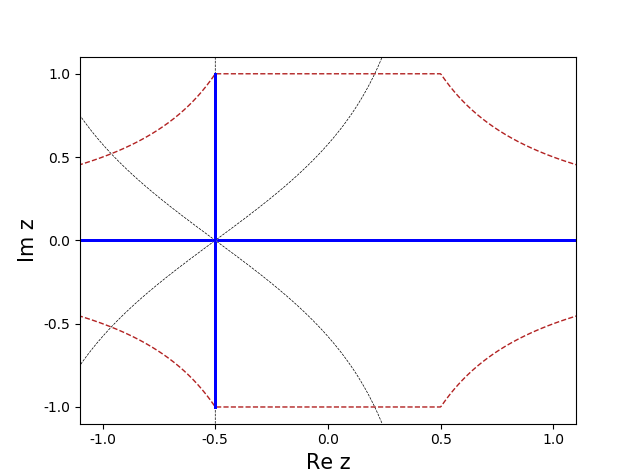}
\includegraphics[width=7.5cm]{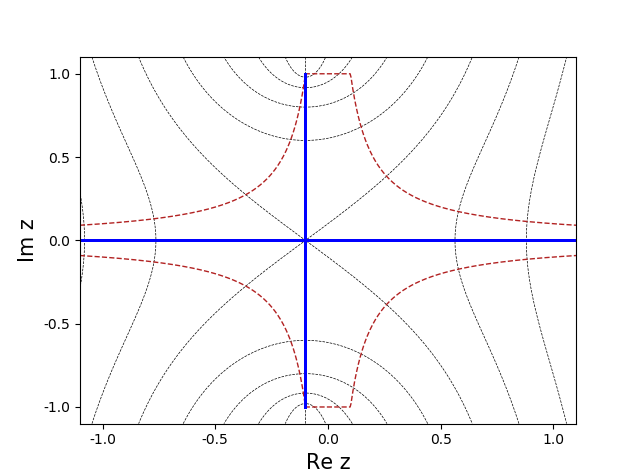}}
\caption{Lax spectrum for the potential $q(x)=\e^{\i x}$ with minimal period $L=2\pi$ (blue). Contours $\Gamma=\{z\in\C:\Im\D(z)=0\}$ (black dashed). The curve which bounds the imaginary component of elements in the spectrum $|\Im z| = \min\{\| q\|_{\infty},\, \ep\|q'\|_{\infty}/2|\Re z|\}$ (red dashed).
Left: $\ep=1$. Right: $\ep=0.2$.}
\label{f:4}
\end{figure}

\begin{figure}[t!]
\vglue-2ex
\centerline{\includegraphics[width=7.5cm]{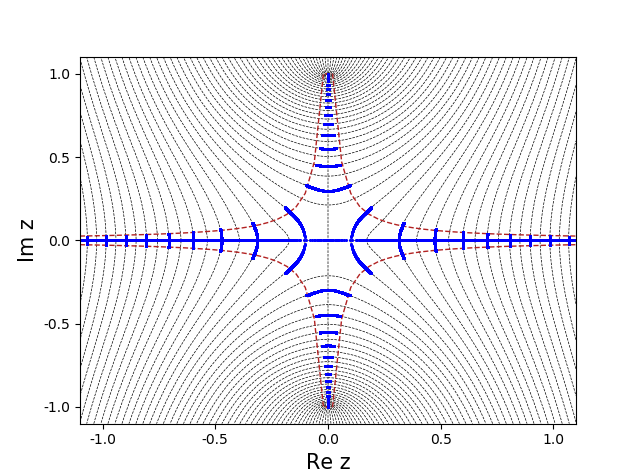}
\includegraphics[width=7.5cm]{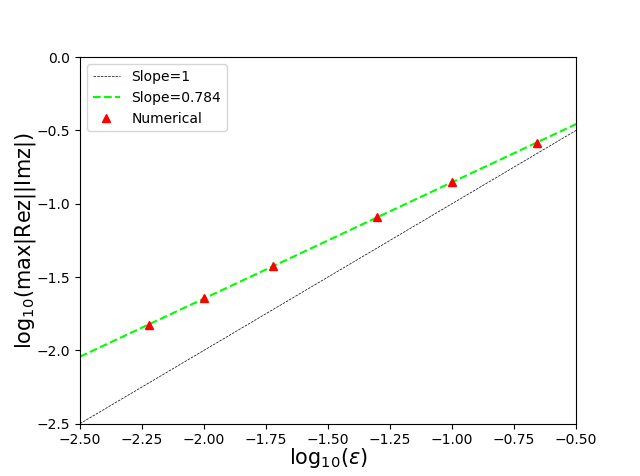}}
\caption{Left: Same as Fig.~\ref{f:4} but for the periodic potential $q(x) = \text{sgn}(x)$ with $L = 2$, and $\ep=0.019$. 
Right: Convergence of eigenvalues to $\Sigma_{\infty}$ as $\ep\downarrow 0$ for ``sgn'' potential. 
Numerical computation of ${\rm max}_{\nu\in[0,2\pi/L)}|\Re z(\nu)||\Im z(\nu)|$ as $\ep\downarrow 0$ (red triangles). Least-squares fit through the data points (light green dashed).} 
\label{f:5}
\end{figure}

In this case, 
one could consider a bound on the imaginary component of the eigenvalue similar to \eqref{e:maininequality}, 
but obtained 
using the numerically calculated eigenvalues at $\ep=1$, 
and then examine (numerically) how the spectrum changes as a function of $\ep$.
The numerical results indicate that, 
while the spectrum is still confined to a region of the complex plane similar to~\eqref{e:eigenvalue_set}, 
the size of the region is no longer simply proportional to $\ep$. 
Indeed, 
a numerical study shown in Fig.~\ref{f:5} (right) suggests that
one has $\max_{\nu\in [0,2\pi/L)}|\Re z(\nu)||\Im z(\nu)| = O(\epsilon^{\alpha})$ as $\epsilon \downarrow 0$, 
with $\alpha = 0.784<1$ (see Fig.~\ref{f:5}).
(Whereas, 
if the potential satisfied the hypotheses of Corollary~\ref{c:mainresult}, 
one would have $\alpha=1$.)
Nonetheless, 
the numerical results clearly indicate that the spectrum tends to the real and imaginary axes in the limit $\ep\downarrow 0$.
An interesting future direction would be to find $\alpha$ rigorously for piecewise continuous potentials.
Similar results are also obtained for potentials with steps of arbitrary magnitude (owing to the invariant properties of the scattering problem under scaling transformations)
as well as for potentials with asymmetric steps (e.g., $q(x) = 0$ for $x<0$ and $q(x) = 1$ for $x>0$).
A detailed discussion is omitted for brevity.

\subsection{Real-valued periodic single-lobe potentials}
\label{s:singlelobe}
Next we consider a periodic analog of the Klaus and Shaw single-lobe potentials~\cite{KS2002}. 
Specifically, 
we take the potential $q$ to be the $L$-periodic extension of a real-valued continuously differentiable function on the interval $(-L/2, L/2)$.
We further assume:
(i) $q(x)\geq0$ for all $x\in\Real$,
(ii) $q(-x) = q(x)$,
and (iii) $q(x)$ is increasing on $[-L/2, 0)$ and is decreasing on $(0, L/2]$.
In particular, 
below we consider two examples:
\bse
\label{e:IC}
\begin{align}
q(x) &= \e^{-\sin^2x}\,, \quad L = \pi\,,
\label{e:IC1}
\\
q(x) &= \dn(x;m)\,, \quad L = 2K(m)\,,
\label{e:IC2}
\end{align}
\ese
where $\dn(x;m)$ is a Jacobi elliptic function, 
$m\in(0,1)$ the corresponding elliptic parameter, 
and $K(m)$ the complete elliptic integral of the first kind (see~\cite{NIST} for details).
Moreover, 
when $m=1$ the problem reduces to that studied in~\cite{SatsumaYajima}.

\begin{figure}[t!]
\vglue-2ex
\centerline{\includegraphics[width=7.5cm]{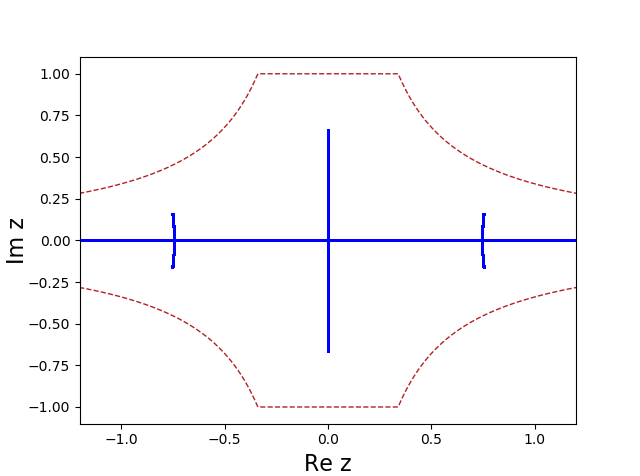}
\includegraphics[width=7.5cm]{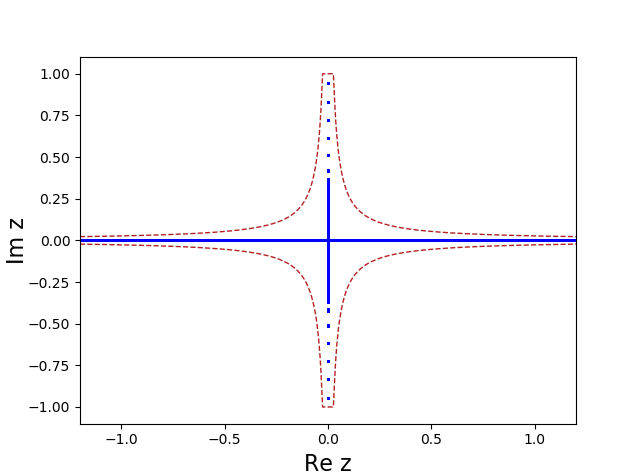}}
\caption{Numerically computed Lax spectrum for the periodic potential $q(x)=\e^{-\sin^{2} x}$ (blue). The curve which bounds the imaginary component of the eigenvalue $|\Im z| = \min \{\|q\|_{\infty},\, \ep\|q'\|_{\infty} / 2|\Re z| \}$ (dark red dashed). Left: Semiclassical parameter $\ep=1$. Right: Semiclassical parameter $\ep=0.079$.}
\label{f:6}
\medskip
\centerline{\includegraphics[width=7.5cm]{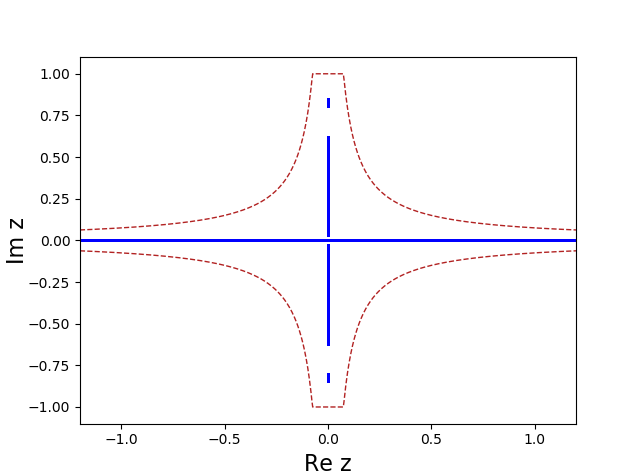}
\includegraphics[width=7.5cm]{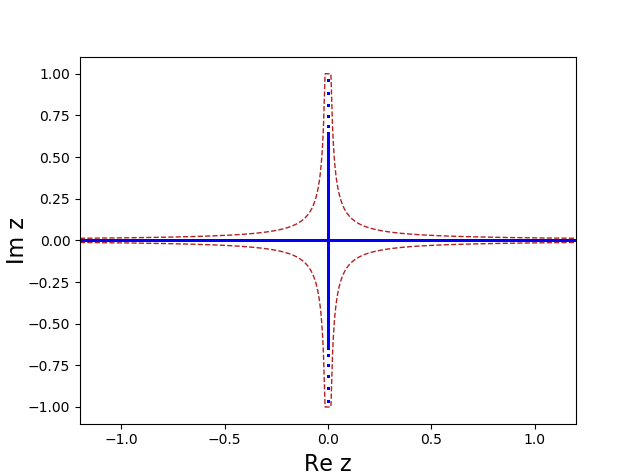}}
\caption{Same as Fig.~\ref{f:6}, but for the periodic potential $q(x)=\dn(x;m)$, where $m=0.6$. 
Left: $\ep=0.5$. Right: $\ep=0.1$.}
\label{f:7}
\end{figure}

Figures~\ref{f:6} and \ref{f:7} show the results of numerical computations of the Lax spectrum of the ZS system~\eqref{e:zs} for the two potentials in~\eqref{e:IC}. 
The results clearly demonstrate that, 
in the case of periodic single-lobe potentials, 
there exists complex Floquet eigenvalues of~\eqref{e:zs} off of the imaginary axis. 
This is in contrast to single-lobe potentials decaying as $x\to \pm\infty$, 
whose eigenvalues were proven to be only real and purely imaginary, 
both for potentials with zero BCs~\cite{KS2002} and for potentials with symmetric non-zero BCs~\cite{BL2018}.
In contrast, 
Figs.~\ref{f:6} and~\ref{f:7} show that, 
in the periodic case, 
the confinement to the real and imaginary axes does not hold in general.
Indeed, 
it would be surprising if such a feature were present in the periodic case as this would form a large class of finite-band potentials by Theorem~\ref{t:finiteband2}.
Another interesting feature of the spectrum for single-lobe periodic potentials is the formation of spectral bands and gaps along the intervals $\pm (\i q_{\min}, \i q_{\max})$ of the spectral plane. 
Using WKB approximation it was shown formally in~\cite{BiondiniOregero} that the number of spectral bands in this interval is $O(1/\ep)$ as $\ep \downarrow 0$.
Another interesting question is whether spines can occur along the imaginary axis for single-lobe periodic potentials. 
The numerical results suggest that the Lax spectrum is purely imaginary outside of a neighborhood of the real $z$-axis. 
This property was shown to hold in the limit $\ep \downarrow 0$ (see Theorem~\ref{l:wkbbound}) where the neighborhood was proportional to $\ep$.

\begin{figure}[t!]
\vglue-2ex
\centerline{\includegraphics[width=7.5cm]{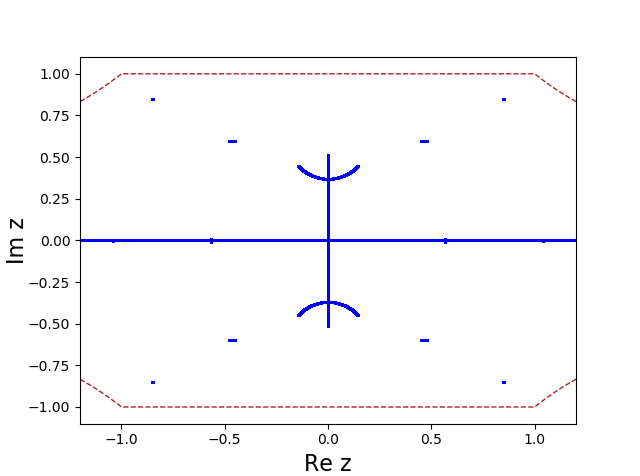}
\includegraphics[width=7.5cm]{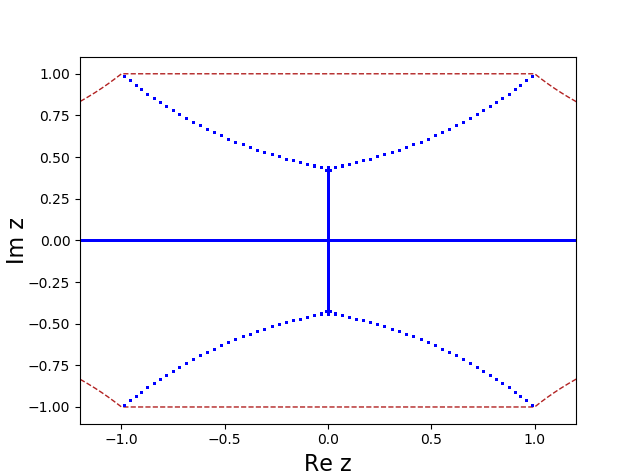}}
\caption{Same as Fig.~\ref{f:6} but for the periodic potential $q(x;\ep) = \e^{\i\cos(2x)/\ep}$. 
Left: $\ep=0.22$. Right: $\ep=0.019$.}
\label{f:8}
\medskip
\centerline{\includegraphics[width=7.5cm]{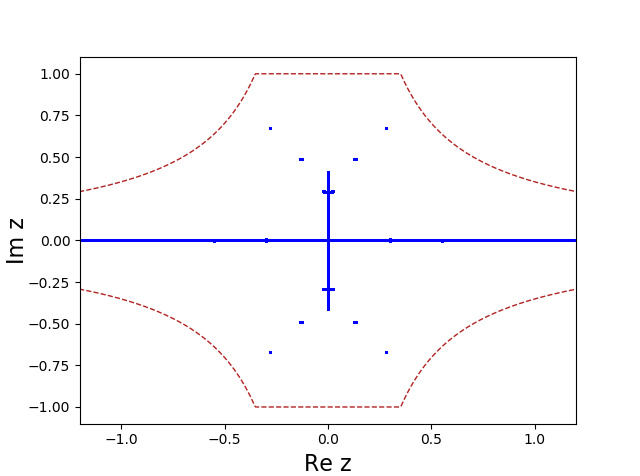}
\includegraphics[width=7.5cm]{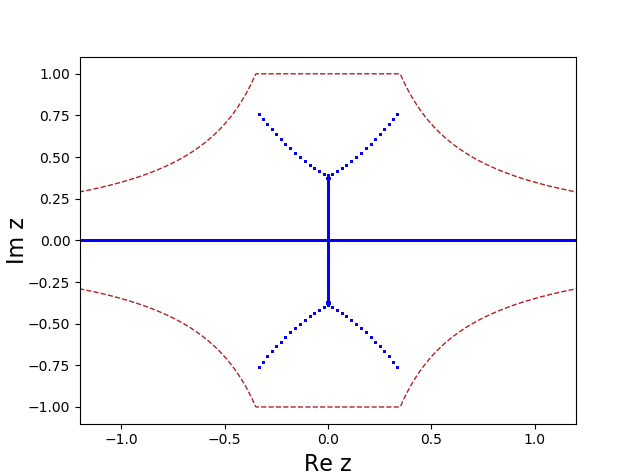}}
\caption{Same as Fig.~\ref{f:6} but for the periodic potential $q(x;\ep) = \text{dn}(x;m)\e^{\i2\text{dn}(x;m)/\ep}$, where $m=0.88$. 
Left: $\epsilon=0.2$. Right: $\ep=0.03$.}
\label{f:9}
\end{figure}
\subsection{Periodic potentials with rapid phase variations}

Here we let the potential depend on the semiclassical parameter,
that is, $q = q(x;\epsilon)$. 
In particular, 
we consider a periodic potential with rapidly varying phase,
i.e., 
\be
q(x;\ep) = A(x)\e^{\i S(x)/\ep}\,,
\label{e:rapid_phase}
\ee
where $A(x)$ and $S(x)$ are continuously differentiable real functions independent of $\ep$. 
Specifically, 
we consider the following two potentials:
\bse
\label{e:fast_phase}
\begin{align}
q(x;\ep) &= \e^{\i\cos(2x)/\ep}\,, \quad L = \pi\,,
\\
q(x;\ep) &= \dn(x;m)\e^{2\i \dn(x;m)/\ep}\,, \quad L = 2K(m)\,.
\end{align}
\ese
Importantly,
the bound~\eqref{e:maininequality} still holds.
The key difference from the previous cases, 
however, 
is that now
$\|q'\|_{\infty} = O(1/\ep)$ as $\ep\downarrow 0$.
Hence, 
one does not expect Theorem~\ref{c:mainresult} to hold for potentials of the
form~\eqref{e:fast_phase}. 
Moreover,
using turning point curves and ideas of Deift,
Venakides,
and Zhou (see~\cite{difranco} for details) one expects that sharper bounds beyond what is available in~\eqref{e:maininequality} can be obtained.
This is an interesting direction for future work.
Indeed, 
the results of numerical computations of the spectrum produced by the potentials~\eqref{e:fast_phase}, 
and depicted in Figs~\ref{f:8} and~\ref{f:9},
show that 
there are complex spectral bands off of the real and imaginary axes that persist in the limit $\ep\downarrow 0$. 
Importantly the numerical results are still consistent with the inequality~\eqref{e:maininequality}.
More interestingly, 
the spectrum appears to accumulate on a set of well-defined curves in the complex plane as $\ep\downarrow 0$. 
This is qualitatively the same Y-shaped curve observed in~\cite{bronski} for the zero-background potential $q(x;\ep) = \sech(2x)\e^{\i \sech(2x)/\ep}$ 
on the infinite line. 
This is another example in which the differences in the spectrum due to different choices of boundary conditions seem to become negligible in 
the semiclassical limit of the focusing problem.

\section*{Acknowledgments}

We thank the anonymous referee for the careful review and the many helpful comments.\break
This work was partially supported by the National Science Foundation under grant number DMS-2009487.

\section*{Appendix}
\setcounter{section}1
\setcounter{subsection}0
\setcounter{equation}0
\def\thesection{\Alph{section}}
\def\theequation{\Alph{section}.\arabic{equation}}
\def\thetheorem{\Alph{section}.\arabic{theorem}}
\def\thefigure{\Alph{section}.\arabic{figure}}

\subsection{Definitions}
\label{s:notations}
At various points in this work we make use of the Pauli matrices which are defined as:
\be
\sigma_{1} := \begin{pmatrix} 0 & 1 \\ 1 & 0 \end{pmatrix}\,, \quad \sigma_{2} := \begin{pmatrix*}[r] 0 & -\i \\ \i & 0 \end{pmatrix*}\,, \quad \sigma_{3} := \begin{pmatrix*}[r] 1 & 0 \\ 0 & -1 \end{pmatrix*}\,.
\label{e:pauli}
\ee
Next we discuss the normed linear spaces used throughout the work.  
Let $g$ be a Lebesgue measurable function. 
Then the essential supremum is defined as:
\be
\n g\n_{\infty}:=\inf \{ C : m( \{x\in\Real : |g(x)| > C \} ) = 0 \}\,,
\label{e:inf}
\ee
where $[m]$ is the Lebesgue measure.
Note we always identify functions equal almost everywhere ``a.e." with respect to Lebesgue measure.
Specifically, 
$f\equiv g$ if 
$m(\{x\in\Real:f(x) \neq g(x)\}) = 0$.
Thus a Lebesgue measurable function is in $L^{\infty}(\Real)$ if it's bounded a.e. [$m$]. 
If $g \in L^{\infty}(\R)$ is periodic one easily gets $g \in L^1_{\text{loc}}(\R)$,
that is,
$g$ is Lebesgue integrable on compact subsets of the real numbers.

Also,
we define the inner product:
\be
\langle f,g \rangle := \int_{0}^{L} f(x)\overline{g(x)}\,\d x\,,
\label{e:scalarprod}
\ee
where $f$, $g$ are scalar Lebesgue measurable functions.
This gives
\be
\|f\|_{2} := \sqrt{\langle f,f \rangle} = \Big(\int_{0}^{L} |f(x)|^2\, \d x\Big)^{1/2}\,,
\label{e:scalarl2}
\ee
and $L^2([0,L])$ is the space of scalar Lebesgue measurable functions which are square integrable,
that is,
$\|f\|_2<\infty$. 
Similarly, 
we define the inner product:

\be
\langle \@f, \@g \rangle := \int_{0}^{L} f_{1}(x)\overline{g_{1}(x)} + f_{2}(x)\overline{g_2(x)}\, \d x\,,
\label{e:inproduct}
\ee
where $\@f$, $\@g$ are two-component Lebesgue measurable vector functions. 
This gives
\be
\n \@f \n_{2} := \sqrt{\langle \@f, \@f \rangle} = \Big( \int_{0}^{L} |f_{1}(x)|^2 + |f_{2}(x)|^2\, \d x \Big)^{1/2}\,,
\label{e:2norm}
\ee
and $L^2([0,L],\C^2)$ is the space of two-component vector functions which are square integrable,
that is,
$\|\@f\|_2<\infty$.

Finally, 
a function $f : [a,b] \to \C$ is absolutely continuous on $[a,b]$ if $\forall \eta > 0$, $\exists \delta > 0$ such that whenever a finite sequence of pairwise disjoint sub-intervals $[a_k, b_k]$ of $[a,b]$ satisfies $\sum_{k}|b_k-a_k| < \delta$ then $\sum_{k}|f(b_k)-f(a_k)| < \eta$. 
The collection of such functions is denoted $AC([a,b])$.
Importantly, 
$f \in AC([a,b])$ implies $f'$ exists a.e $[m]$,
and $f'$ is Lebsegue integrable. 
A function $f\in AC_{\text{loc}}(\R)$ if it is absolutely continuous on compact subsets of the real axis.

\subsection{Symmetries of the ZS scattering problem}
\label{s:symmetries}

Recall next some symmetries of the solutions of~\eqref{e:zs}, 
all of which are easily verified by direct computation. 
Let $Y(x;z,\ep)$ be a fundamental matrix solution of the ZS scattering problem~\eqref{e:zs}.
Then another solution of \eqref{e:zs} is
\be
\widetilde{Y}(x;z,\ep) = \i\sigma_{2} \,\overline{Y(x;\overline{z},\ep)}\,.
\label{e:sym1}
\ee
Moreover,
if the potential is real, 
a further solution of \eqref{e:zs} is given by
\be
\widetilde{Y}_{(r)}(x;z,\ep) = \overline{Y(x;-\overline{z},\ep)}\,.
\label{e:sym2}
\ee
Alternatively suppose the potential satisfies a generalized reflection symmetry,
that is,
$q(-x) = \e^{2\i\theta}q(x)$ for all $x\ge0$, 
for some $\theta \in \Real$.
(Obviously for $\theta=0$, 
and $\theta = \pi/2$ one has the cases of even and odd potentials, 
respectively.)
Then another solution of \eqref{e:zs} is given by
\be
\widetilde{Y}_{(g)}(x;z,\ep) = (\cos\theta\sigma_1 + \sin\theta\sigma_2)\,\overline{Y(-x,-\overline{z},\ep)}\,.
\label{e:sym3}
\ee
Finally, 
if the potential is PT-symmetric, 
i.e., 
if $q(-x) = \overline{q(x)}$, 
then another solution of~\eqref{e:zs} is given by 
\be
\widetilde{Y}_{(pt)}(x;z,\ep) = \sigma_3\overline{Y(-x;\overline{z},\ep)}\,.
\label{e:sym4}
\ee

\subsection{Comparison between various definitions of Dirichlet spectrum}
\label{e:defndirichlet}

In this section we show how to connect the Dirichlet spectrum
to the zeros of an analytic function.
Recall the Dirichlet BCs:
\begin{gather}
v_1(0;\zeta,\ep) + v_2(0;\zeta,\ep) = 0\,,\qquad
v_1(L;\zeta,\ep) + v_2(L;\zeta,\ep) = 0\,,\qquad
\label{e:BCs}
\end{gather}
where $\{\zeta_n\}$ denote Dirichlet eigenvalues and $x_o=0$ without loss of generality.
Let $\@v(x;z,\ep)=\Phi(x;z,\ep)\@c(z;\ep)$,
where, 
as before, 
$\Phi$ is a fundamental matrix solution of~\eqref{e:zs} normalized so that $\Phi(0;z,\ep)\equiv\I$,
and $\@c(z;\ep) = (c_1, c_2)^T$.
Clearly, 
$\@v(0;z,\ep) = \@c(z;\ep)$.
If $\zeta\in\Sigma_{\rm Dir}(0)$ and $\@v(x;\zeta,\ep)$ is the eigenfunction associated to $\zeta$,
by the first of~\eqref{e:BCs} one can write
\be
\@v(L;\zeta,\ep) = M(\zeta;\ep)\@c(\zeta;\ep) = v_{1}(0;\zeta,\ep)\begin{pmatrix*} M_{11} - M_{12} \\ M_{21} - M_{22} \end{pmatrix*}\,,
\label{e:Dir_M21}
\ee
where $M(z;\ep) = \Phi(L;z,\ep)$ is the monodromy matrix associated with $\Phi$, 
as per~\eqref{e:monodromy} and~\eqref{e:monodromy2}. 
Using the second of~\eqref{e:BCs} we then have 
\be
\big( M_{11} - M_{12} + M_{21} - M_{22} \big)|_{z=\zeta} = 0\,.
\label{e:dir_zero1}
\ee
Now, 
using the modified fundamental matrix solution $\widetilde{\Phi}$ defined in~\eqref{e:unitary},  
recall that the associated monodromy matrix is $\widetilde{M}^{\ep}(z) = C^{-1} M(z;\ep)C$.
One therefore gets
\be
\widetilde{M}^\ep_{21}(\zeta) = \halfi(M_{11} - M_{12} + M_{21} - M_{22})|_{z=\zeta}=0\,.
\label{e:dir_zero2}
\ee
Hence $\zeta\in\Sigma_{\rm Dir}(0)$ if and only if $\widetilde{M}^\ep_{21}(\zeta) = 0$,
in agreement with~\eqref{e:dirzeros}.

Note that some authors define the Dirichlet spectrum via slightly different boundary conditions, 
namely,
\begin{gather}
v_1(0;\zeta,\ep)-v_2(0;\zeta,\ep)=0\,, \qquad
\quad v_{1}(L;\zeta,\ep)-v_{2}(L;\zeta,\ep)=0\,, \quad
\end{gather}
(see~\cite{djakovmityagin, kapmity, kappeler}). 
Using similar arguments as above, 
with these BCs one gets
\be
\big( M_{11} - M_{21} + M_{12} - M_{22} \big)|_{z=\zeta} = 0\,.
\label{e:dir_zero3}
\ee
We then introduce the modified fundamental matrix solution $\breve\Phi(x;z,\ep) = \Phi(x;z,\ep)\breve C$,
where
\be
\breve C = \frac{1}{\sqrt{2}} \begin{pmatrix*}[r] \i & - 1 \\ \i & 1 \end{pmatrix*}\,.
\label{e:unitary2}
\ee
For the associated monodromy matrix $\breve{M}^{\ep}(z) = \breve C^{-1} M(z;\ep) \breve C$,
one gets 
\be
\breve{M}_{21}^{\ep}(\zeta) = -\halfi(M_{11} - M_{21} + M_{12} - M_{22})|_{z=\zeta}=0\,.
\label{e:dir_zero4}
\ee
Hence,
with either choice of BCs, 
one can relate the Dirichlet spectrum to the zeros of an analytic function which is 
the row 2 column 1 element of a modified monodromy matrix. 
The trace formulae~\eqref{e:trace1} and~\eqref{e:trace2} use definition~\eqref{e:BCs}.
The Dirichlet eigenvalues obey certain ODEs (called Dubrovin equations) with respect to $x_o$ and $t$~\cite{ForestLee,MaAblowitz}.

Importantly, 
it is straightforward to see that the bounds obtained for the Dirichlet spectrum are independent of which choice of boundary conditions is used. 

\subsection{Floquet discriminant for constant and plane wave potentials}
\label{s:planewave}

Here we calculate the Lax spectrum for plane wave potentials, 
namely, 
$q(x)$ given by~\eqref{e:plane_wave}.
In this case the ZS system~\eqref{e:zs} is given by
\be
\ep\@v' = \begin{pmatrix} -\i z & A\e^{\i Vx} \\ -A\e^{-\i Vx} & \i z \end{pmatrix} \@v\,.
\label{e:zs_ap}
\ee
As a special case, 
when $V=0$ we obtain a constant background potential.
In this case~\eqref{e:zs_ap} is a constant-coefficient differential equation, 
and it is easy to obtain a fundamental matrix solution as:
\be
Y(x;z,\ep) = \Big(\mathbb{I} - \frac{\i A}{z - \xi} \sigma_{1}\Big) \e^{\i\xi x\sigma_{3}/\ep}\,,
\ee
where $\mathbb{I}$ is the $2\times 2$ identity matrix, 
$\xi^2 = A^2 + z^2$, 
and the Pauli matrices $\sigma_{1}$ and $\sigma_{3}$ are as in~\eqref{e:pauli}. 
A simple calculation then yields a monodromy matrix~\eqref{e:monodromy} as 
\be
M(z;\ep) = \e^{\i\xi L \sigma_{3}/\ep}\,,
\label{e:monodromy_const}
\ee
and hence,
\be
\D(z) = \cos(\xi L/\ep)\,.
\label{e:delta_0}
\ee
It follows that $z \in \Sigma_{\Lax}$ if and only if $\xi \in \Real$. 
Thus, 
for a constant background potential one gets 
\be
\Sigma_{\Lax} = \Real \cup [-\i A, \i A]\,.
\label{e:constback}
\ee 
Additionally, 
using the similarity transformation~\eqref{e:unitary} one gets 
\be
\widetilde{M}_{21}^{\ep}(z) = \Big(\frac{z-\i A}{2\xi}\Big)\sin(\xi L/\ep)\,.
\label{e:dirconst}
\ee
It then follows that the Dirichlet spectrum is the set of zeros of the analytic function~\eqref{e:dirconst}.

When $V\ne0$,~\eqref{e:zs_ap} is not constant-coefficient. 
Nonetheless,
the substitution $\@v=\e^{\i Vx\sigma_3/2}\@u$ transforms~\eqref{e:zs_ap} to the constant-coefficient system
\be
\ep\@u' = \begin{pmatrix*} -\i(z+\ep V/2) & A \\ -A & \i(z+\ep V/2) \end{pmatrix*}\@u\,. 
\ee
As before,
we obtain a fundamental matrix solution,
namely,
\be
Y(x;z,\ep) = \e^{\i Vx \sigma_{3}/2} \bigg(\mathbb{I}-\frac{\i A}{z-\xi + \ep V/2} \sigma_{1} \bigg) \e^{\i\xi x \sigma_{3}/\ep}\,,
\label{e:jost}
\ee
where $\xi^2 = A^{2} + (z + \ep V/2)^{2}$.
Simple calculations then yield a monodromy matrix
\be
M(z;\ep) = \Big(\cos(VL/2) + \sin(VL/2) \big( -\i(z+\ep V/2)\sigma_{3} - \i A \sigma_2 \big)\Big) \,\e^{\i \xi L\sigma_3/\ep},  
\label{e:planewave_monodromy}
\ee
and hence,
\be
\D(z) = \cos(VL/2)\cos(\xi L/\ep) +  \frac{z+\ep V/2}{\xi}\sin(VL/2)\sin(\xi L/\ep)\,.
\label{e:trm_ap}
\ee
Take $L = 2\pi/|V|$. 
Then~\eqref{e:trm_ap} simplifies to $\D(z) = -\cos(2\pi \xi/\ep V)$. 
It follows that $z\in\Sigma_{\Lax}$ if and only if $\xi\in\Real$. 
Then, 
the Lax spectrum for a plane wave potential is given by~\eqref{e:spec_plnwave}
(see Fig.~\ref{f:4}).

\subsection{Floquet discriminant for piecewise constant potentials}
\label{s:signum}

As one last example we calculate the trace of the monodromy matrix for 
$q(x)=A{\rm sgn}(x)$, 
where $A>0$ is a constant, 
and ``${\rm sgn}$" is given by~\eqref{e:signum}.
In this case the ZS system~\eqref{e:zs} is given by
\be
\ep \@v' = \begin{pmatrix} -\i z & A\text{sgn}(x) \\ -A\text{sgn}(x) & \i z \end{pmatrix} \@v\,.
\label{e:piececonst}
\ee
Solving the above system on the intervals $[-L/2, 0)$, 
and $[0, L/2)$, 
and imposing the normalization $\Phi(0;z,\ep)\equiv\mathbb{I}$ one finds the matrix solution
\be
\Phi(x;z,\ep) = \frac{\xi-z}{2\xi}\Big(\mathbb{I}-\frac{\i A}{z-\xi} \sigma_1 \Big) \e^{\i\xi x \sigma_3/\ep}\Big(\mathbb{I} + \frac{\i A}{z-\xi} \sigma_1 \Big)\,,
\label{e:pms}
\ee
where $\xi^2 = A^2 + z^2$. 
It follows that 
\be
M(z;\ep) = \Phi^{-1}(-L/2;z,\ep)\Phi(L/2;z,\ep)\,,
\label{e:piecemonodromy}
\ee
and hence,
\be
\D(z) = \frac{1}{\xi^2}(A^2 + z^2\cos(\xi L/\ep))\,.
\label{e:discrimpiece}
\ee

\subsection{Floquet-Hill's method}
\label{s:fhm}

Here we provide some details regarding the numerical calculations of the Lax spectrum in the main text. 
Recall that the Zakharov-Shabat scattering problem is given by the first-order system of ODEs \eqref{e:zs}.
Further recall that we can rewrite \eqref{e:zs} in the form of an eigenvalue problem \eqref{e:op}.
Since \eqref{e:op} is non-self-adjoint, 
when computing its spectrum numerically one needs an approach capable of accurately calculating eigenvalues in a large region of the complex plane. 
One such approach is Floquet-Hill's method (see \cite{DK2006} for more details). 
The method is particularly well-suited for calculating the eigenvalues of linear operators with periodic coefficients, 
and provides a global approximation of the spectrum.

Recall that Floquet's theorem~\ref{t:floquet} implies that any bounded solution of~\eqref{e:zs} is necessarily of the form
\be
\@v(x;z,\ep) = \e^{\i\nu x}\@w(x;z,\ep)\,,
\label{e:bounded_a}
\ee
where $\@w(x+L;z,\ep) = \@w(x;z,\ep)$ and $\nu\in [0,2\pi/L)$.
Plugging~\eqref{e:bounded_a} into the ZS system~\eqref{e:zs} yields the modified eigenvalue problem 
\be
\i\sigma_{3} \big( \ep(\partial_{x} + \i\nu) - Q \big) \@w = z\@w\,.
\label{e:mevp}
\ee
The key point is that the modified eigenvalue problem acts on periodic functions. 
Therefore one can expand~\eqref{e:mevp} using a discrete Fourier transform (DFT):
\be
\widehat{\mathcal{L}^{\ep}_{\nu,N}}\hat{\@w}_N = z_N\hat{\@w}_N\,,
\label{e:fs_a}
\ee
where  
\be
\widehat{\mathcal{L}^{\ep}_{\nu,N}} := \begin{pmatrix} -\ep(k + \nu) & -\i \mathfrak{T} \\ -\i \overline{\mathfrak{T}} & \ep(k + \nu) \end{pmatrix}_{2N\times 2N}\,.
\label{e:fs_mat_a}
\ee
Moreover,
$k=\diag(k_{n})$ with $n\in\{-\halfN,\dots, \halfN-1\}$ is a diagonal matrix of Fourier modes such that $k_{n} = 2n\pi/L$, 
and $\mathfrak{T}$ is a $N\times N$ Toeplitz matrix representing the convolution operator generated by the DFT of $\widehat{q\@w}$. 

Floquet-Hill's method then approximates the Floquet spectrum by numerically computing the eigenvalues of the matrix~\eqref{e:fs_mat_a}.
Indeed,
fix $\nu\in[0,2\pi/L)$.
Then choosing a truncation $N=2^j$ of the number of Fourier modes of the eigenfunction $\@w(x;z,\ep)$ results in a matrix system of dimension $2N$.
Moreover,
$z_N$ are approximations in the sense that all $z_N\in\Sigma(\widehat{\mathcal{L}^{\ep}_{\nu,N}})\to z\in\Sigma_{\nu}$ as $N\to\infty$. 
By taking an evenly spaced sequence $\nu_{i}\in[0,2\pi/L)$ one approximates the entire Lax spectrum,
i.e.,
\be
\lim_{N\to\infty} \bigcup_{\nu\in[0,2\pi/L)} \Sigma(\widehat{\mathcal{L}^{\ep}_{\nu,N}}) = \Sigma_{\Lax}.
\label{e:almost}
\ee
(For details on convergence see~\cite{johnson}.)
Numerical accuracy of the approximation is determined by the number of Fourier modes used in the truncation, 
and on the method used to compute the eigenvalues of the matrix~\eqref{e:fs_mat_a}.
The resolution of the spectral bands is determined by how fine we partition the interval $\nu \in [0, 2\pi/L)$. 
All computations in the main text were checked for numerical convergence.

\makeatletter
\makeatother
\def\title#1,{\relax}


\begin{thebibliography}O
\itemsep0pt
\parsep0pt
\small

\bibitem{APT2004}
M.J. Ablowitz, B. Prinari and A.D. Trubatch
\reftitle{Discrete and continuous nonlinear Schr\"odinger systems},
Cambridge University Press, 2004

\bibitem{az2014}
D.S. Agafontsev and V.E. Zakharov, 
\reftitle{Integrable turbulence and formation of rogue waves}, 
\href{http://iopscience.iop.org/article/10.1088/0951-7715/28/8/2791/pdf}
{Nonlinearity, \textbf{28}, 2791--2821 (2015)}

\bibitem{BBEIM}
E.D. Belekolos, A.I. Bobenko,  V.Z. Enol'skii, A.R. Its and V.B. Matveev,
\reftitle{Algebro-geometric approach to nonlinear integrable equations},
Springer-Verlag, 1994 
	
\bibitem{bertolatovbis}
M. Bertola and A. Tovbis,
\reftitle{Universality for the focusing nonlinear Schr\"odinger equation at the gradient catastrophe point: Rational breathers and poles of the tritronqu\'ee solution to Painlev\'e~I},
\href{http://onlinelibrary.wiley.com/doi/10.1002/cpa.21445/epdf}
{Commun. Pure Appl. Math. \textbf{66}, 678--752 (2013)}

\bibitem{BL2018}
G. Biondini and X.-D.\ Luo,
\reftitle{Imaginary eigenvalues of Zakharov-Shabat problems with non-zero background},
\href{https://doi.org/10.1016/j.physleta.2018.06.045}
{Phys. Lett. A \textbf{382}, 2632--2637 (2018)}

\bibitem{BiondiniOregero}
G. Biondini and J. Oregero
\reftitle{Semiclassical dynamics in self-focusing nonlinear media with periodic initial conditions},
\href{https://doi.org/10.1111/sapm.12321}
{Stud.\ Appl.\ Math. \textbf{145}, 325--356 (2020)}

\bibitem{bronski}
J.C. Bronski, 
\reftitle{Semiclassical eigenvalue distribution of the Zakharov-Shabat eigenvalue problem},
\href{\doibase https://doi.org/10.1016/0167-2789(95)00311-8}
{Phys. D \textbf{97}, 376--397 (1996)}

\bibitem{Floquet_int}
B.M. Brown, M.S.P. Eastham and K.M. Schmidt,
\reftitle{Periodic Differential Operators},
Springer Basel, 2013

\bibitem{danilov}
L.I. Danilov,
\reftitle{Absolute continuity of the spectrum of a periodic Dirac operator},
\href{https://doi.org/10.1007/BF02754212}
{Diff. Equat. \textbf{36}, 262--271 (2000)}

\bibitem{DK2006}
B. Deconinck and J.N. Kutz,
\reftitle{Computing spectra of linear operators using the Floquet-Fourier-Hill method},
\href{\doibase 10.1016/j.jcp.2006.03.020}
{J. Comput. Phys. \textbf{219}, 296--321 (2006)}

\bibitem{DVZ}
P. Deift, S. Venakides and X. Zhou,
\reftitle{New results in small dispersion KdV by an extension of the steepest descent method for Riemann-Hilbert problems}
\href{https://doi.org/10.1155/S1073792897000214}
{Int. Math. Res. Notices \textbf{1997}, 285--299 (1997)}
	
\bibitem{deiftmclaughlin}
P. Deift and K.D.T.-R. McLaughlin,
\reftitle{A continuum limit of the Toda lattice}
{Mem. Amer. Math Soc., \textbf{131}, 1998}

\bibitem{pre2017deng}
G. Deng, S. Li, G. Biondini and S. Trillo,
\reftitle{Recurrence due to periodic multi-soliton fission in the defocusing nonlinear Schr\"odinger equation},
\href{http://www.math.buffalo.edu/~biondini/papers/pre2017v96p052213.pdf}
{Phys. Rev. E \textbf{96}, 052213 (2017)}

\bibitem{difranco}
J.C. DiFranco and P.D. Miller,
\reftitle{Semiclassical modified nonlinear Schr\"{o}dinger equation I: Modulation theory and spectral analysis},
\href{https://doi.org/10.1016/J.PHYSD.2007.11.022}
{Phys. D \textbf{237}, 947--997 (2008)}

\bibitem{djakovmityagin}
P. Djakov and B. Mityagin,
\reftitle{Instability zones of periodic 1-dimensional Schr\"odinger and Dirac operators},
\href{http://dx.doi.org/10.1070/RM2006v061n04ABEH004343}
{Russian Mathematical Surveys \textbf{61}:4 663--766 (2006)}

\bibitem{DudleyTaylor}
J.M. Dudley and J.R. Taylor,
\reftitle{Supercontinuum generation in optical fibers}
Cambridge University Press, 2010

\bibitem{Eastham}
M.S.P. Eastham, 
\reftitle{The spectral theory of periodic differential equations},
Scottish Academic Press, 1973
	
\bibitem{elhoefer}
G.A. El and M.A. Hoefer
\reftitle{Dispersive shock waves and modulation theory},
\href{https://doi.org/10.1016/j.physd.2016.04.006}
{Phys. D \textbf{333}, 11--65 (2016)}

\bibitem{elkhamistovbis}
G.A. El, E.G. Khamis and A. Tovbis,
\reftitle{Dam break problem for the focusing nonlinear Schr\"odinger equation and the generation of rogue waves},
\href{https://doi.org/10.1088/0951-7715%2F29%2F9%2F2798}
{Nonlinearity \textbf{29},  2798--2836 (2016)}

\bibitem{eltovbisSG}
G.A. El and A. Tovbis
\reftitle{Spectral theory of soliton and breather gases for the focusing nonlinear Schr\"odinger equation},
\href{https:doi.org/10.1103/PhysRevE.101.052207}
{Phys. Rev. E \textbf{101}, 052207 (2020)}
	
\bibitem{RamyEl}
R. El-Ganainy, K.G. Makris, M. Khajavikhan, Z.H. Musslimani, S. Rotter and D.N. Christodoulides,
\reftitle{Non-Hermitian physics and PT symmetry},
\href{https://doi.org/10.1038/nphys4323}
{Nature Physics \textbf{14}, 11-19 (2018)}
	
\bibitem{Floquet}
G. Floquet, 
\reftitle{Sur les \'equations différentielles linéaires \`a coefficients p\'eriodiques}, 
Ann. École Normale Sup. \textbf{12}, 47--88 (1883)

\bibitem{ForestLee}
M.G. Forest and J.E. Lee,
\reftitle{Geometry and modulation theory for the periodic NLS equation},
\href{\doibase 10.1007/978-1-4613-8689-6_3}
{IMA Volumes in Mathematics and Its Applications, vol 2 (1986)}
	
\bibitem{FT1987}
L.D. Faddeev and L.A. Takhtajan,
\textit{Hamiltonian methods in the theory of solitons},
Springer, Berlin, 1987

\bibitem{fujiiewittsten2018}
S. Fujii\'e and J. Wittsten,
\reftitle{Quantization conditions of eigenvalues for semiclassical Zakharov-Shabat systems on the circle},
\href{https://doi.org/10.3934/dcds.2018167}
{Discr. Cont. Dyn. Syst. \textbf{38}, 3851--3873 (2018)}

\bibitem{Gesztesy_Hill1996}
F. Gesztesy and R. Weikard
\reftitle{Picard potentials and Hill's equation on a torus},
\href{https://doi.org/10.1007/BF02547336}
{Acta Math. \textbf{176}, 73--107 (1996)}

\bibitem{gesztesyweikard_acta1998}
F. Gesztesy and R. Weikard,
\reftitle{A characterization of all elliptic algebro-geometric solutions of the AKNS hierarchy},
\href{https://doi.org/10.1007/BF02392748}
{Acta Math. \textbf{181}, 63--108 (1998)}

\bibitem{gesztesyweikard_bams1998}
F. Gesztesy and R. Weikard,
\reftitle{Elliptic algebro-geometric solutions of the KdV and AKNS hierarchies -- An analytic approach},
\href{https://doi.org/10.1090/S0273-0979-98-00765-4}
{Bull. AMS \textbf{35}, 271--317 (1998)}

\bibitem{kapmity}
B. Gr\'ebert, T. Kappeler and B. Mityagin, 
\reftitle{Gap estimates of the spectrum of the Zakharov-Shabat system}
\href{https://doi.org/10.1016/S0893-9659(98)00063-9}
{Appl. Math. Lett. \textbf{11}:4 95--97 (1998)}

\bibitem{kappeler}
B. Gr\'ebert and T. Kappeler,
\reftitle{Estimates on periodic and Dirichlet eigenvalues for the Zakharov-Shabat system},
\href{https://doi.org/10.5167/uzh-22027}
{Asymptotic Analysis \textbf{25}(3-4) 201--237 (2001)}

\bibitem{hislop}
P.D. Hislop and I.M. Sigal,
\reftitle{Introduction to Spectral Theory with Applications to Schr\"odinger Equations},
Springer, 1996

\bibitem{jenkins2013}
R. Jenkins and K.D.T-R. McLaughlin,
\reftitle{Semiclassical limit of focusing NLS for a family of square barrier initial data},
\href{https://doi.org/10.1002/cpa.21494}
{Comm. Pure Appl. Math. \textbf{67}, 246--320 (2013)}

\bibitem{johnson}
M.A. Johnson and K. Zumbrun,
\reftitle{Convergence of Hill's method for nonselfadjoint operators},
\href{https://www.jstor.org/stable/41582715}
{SIAM J. Numer. Anal. \textbf{50}, 64--78 (2012)}

\bibitem{KMM2003}
S. Kamvissis, K.D.T-R. McLaughlin and P.D. Miller,
\reftitle{Semiclassical soliton ensembles for the focusing nonlinear {S}chr\"odinger equation},
{Princeton 2003}

\bibitem{KS2002}
M. Klaus and J.K. Shaw,
\reftitle{Purely imaginary eigenvalues of Zakharov-Shabat systems},
\href{http://dx.doi.org/10.1103/PhysRevE.65.036607}
{Phys. Rev. E \textbf{65}, 036607 (2002)}

\bibitem{KS2003}
M. Klaus and J.K. Shaw,
\reftitle{On the eigenvalues of Zakharov-Shabat systems},
\href{https://doi.org/10.1137/S0036141002403067}
{SIAM J. Math. Anal. \textbf{34}, 759--773 (2003)}

\bibitem{LaxLevermore}
P.D. Lax and C.D. Levermore,
\reftitle{The small dispersion limit of the {Korteweg-de~Vries} equation {I}, {II} and {III}},
\href{https://doi.org/10.1002/cpa.3160360302}
{Commun. Pure Appl. Math. \textbf{36}, {253--290, 571--593 and 809--829} (1983)}

\bibitem{mclaughlinli}
Y. Li and D.W. McLaughlin,
\reftitle{Morse and Melnikov functions for NLS Pde's},
\href{https://doi.org/10.1007/BF02105191}
{Commun. Math. Phys. \textbf{162} 175--214 (1994)}

\bibitem{lyngmiller}
G.D. Lyng and P.D. Miller,
\reftitle{The N-soliton of the focusing nonlinear Schr\"odinger equation for N large},
\href{https://doi.org/10.1002/cpa.20162}
{Commun. Pure Appl. Math \textbf{40}, 951--1026 (2007)}

\bibitem{MaAblowitz}
Y.-C. Ma and M.J. Ablowitz, 
\reftitle{The periodic cubic NLS equation},
\href{https://doi.org/10.1002/sapm1981652113}
{Stud. Appl. Math. \textbf{63}, 113 (1981)}
	
\bibitem{MW1966}
W. Magnus and S. Winkler, 
\reftitle{Hill's equation},
Dover, 1966

\bibitem{mclaughlinoverman}
D.W. McLaughlin and E.A. Overman II,
\reftitle{Whiskered Tori for Integrable PDE's: Chaotic Behavior in Near Integrable Pde's},
in Surveys in Applied Mathematics, Volume I, 
Ed. J. P. Keller, D. W. McLaughlin and G. P. Papanicolaou,
Plenum, 1995

\bibitem{Messiah}
A. Messiah
\reftitle{Quantum Mechanics},
Elsevier, 1961

\bibitem{miller2001}
P.D. Miller,
\reftitle{Some remarks on a WKB method for the nonselfadjoint Zakharov-Shabat eigenvalue problem with analytic potential and fast phase},
\href{\doibase https://doi.org/10.1016/S0167-2789(01)00166-X}
{Phys. D \textbf{152-53}, 145--162 (2001)}

\bibitem{millerkamvissis1998}
P.D. Miller and S. Kamvissis,
\reftitle{On the semiclassical limit of the focusing nonlinear Schr\"odinger equation},
\href{https://doi.org/10.1016/S0375-9601(98)00565-9}
{Phys. Lett. A \textbf{247}, 75--86 (1998)}
	
\bibitem{trillo2018naturephoton}
A. Mussot, C. Naveau, M. Conforti, A. Kudlinski, F. Copie, P. Szriftgiser and S. Trillo,
\reftitle{Fibre multi-wave mixing combs reveal the broken symmetry of Fermi–Pasta–Ulam recurrence},
\href{https://www.nature.com/articles/s41566-018-0136-1}
{Nature photonics \textbf{12}, 303--308 (2018)}

\bibitem{NMPZ1984}
S.P. Novikov, S.V. Manakov, L.P. Pitaevskii and V.E. Zakharov, 
\reftitle{Theory of solitons: The inverse scattering transform}, 
Plenum, 1984

\bibitem{NIST}
F.W. Olver, D.W. Lozier, R.F. Boisvert and C.W. Clark,
\reftitle{NIST Handbook of Mathematical Functions},
Cambridge University Press, 2010
	
\bibitem{onoratoosborne}
M. Onorato, A.R. Osborne and M. Serio,
\reftitle{Modulational instability in crossing sea states: A possible mechanism for the formation of freak waves},
\href{\doibase 10.1103/PhysRevLett.96.014503}
{Phys. Rev. Lett. \textbf{96}, 014503 (2006)}

\bibitem{rofebek}
F.S. Rofe-Beketov,
\reftitle{On the spectrum of non-selfadjoint differential operators with periodic coefficients},
\href{http://mi.mathnet.ru/eng/dan28724}
{Dokl. Akad. Nauk SSSR \textbf{152}, 1312--1315, (1963)}; English transl. in Soviet Math. Dokl. \textbf{4} (1963)

\bibitem{rofebek2}
F.S. Rofe-Beketov and A.M. Kholkin
\reftitle{Spectral analysis of differential operators: Interplay between spectral and oscillatory properties},
World Scientific, 2005

\bibitem{randoux2014}
S. Randoux, P. Walczak, M. Onorato and P. Suret,
\reftitle{Intermittency in integrable turbulence},
\href{\doibase 10.1103/PhysRevLett.113.113902}
{Phys. Rev. Lett. \textbf{113}, 113902 (2014)}

\bibitem{SatsumaYajima}
J. Satsuma and N. Yajima,
\reftitle{B. Initial value problems of one-dimensional self-modulation of nonlinear waves in dispersive media},
\href{http://ptps.oxfordjournals.org/content/55/284.abstract}
{Progress of Theoretical Physics Supplement, \textbf{55}, 284--306 (1974)}

\bibitem{serov}
M.I. Serov,
\reftitle{On properties of non-selfadjoint differential operators of the second order},
\href{http://mi.mathnet.ru/eng/dan39815}
{Dokl. Akad. Nauk SSSR \textbf{131}, 27--29, (1960)}; English transl. in Soviet Math. Dokl. \textbf{1} (1960)

\bibitem{shin_2004}
K.C. Shin,
\reftitle{On the shape of spectra for non-self-adjoint periodic Schr\"odinger operators},
\href{https://doi.org/10.1088/1361-6447/37/34/007}
{J. Phys. A \textbf{37}, 8287--8291 (2004)}

\bibitem{solli}
D.R. Solli, C. Ropers, P. Koonath and B. Jalali,
\reftitle{Optical rogue waves},
\href{\doibase 10.1038/nature06402} 
{Nature \textbf{450}, 1054--1057 (2007)}

\bibitem{titchmarsh}
E.C. Titchmarsh,
\reftitle{The theory of functions},
Oxford University Press, 1939

\bibitem{tkachenkannals}
V. Tkachenko,
\reftitle{Spectra of non-selfadjoint Hill's operators and a class of Riemann surfaces},
\href{https://doi.org/10.2307/2118642}
{Annals of Math. \textbf{143}, 181--231 (1996)}

\bibitem{tkachenko}
V. Tkachenko,
\reftitle{Non-self-adjoint periodic Dirac operators},
\href{https://doi.org/10.1007/978-3-0348-8247-7_22}
{Operator Theory, System Theory and Related Topics, Birkh\"auser, Basel, pp 485--512, 2001}

\bibitem{tovbisvenakides}
A. Tovbis and S. Venakides,
\reftitle{The eigenvalue problem for the focusing nonlinear Schr\"odinger equation: new solvable cases},
\href{\doibase https://doi.org/10.1016/S0167-2789(00)00126-3}
{Phys. D \textbf{146}, 150--164 (2000)}

\bibitem{TVZ2004}
A. Tovbis, S. Venakides and X. Zhou,
\reftitle{On semiclassical (zero dispersion limit) solutions of the focusing nonlinear Schr\"odinger equation},
\href{https://doi.org/10.1002/cpa.20024}
{Commun. Pure Appl. Math. \textbf{57}, 877--985 (2004)} 

\bibitem{tvz2006}
A. Tovbis, S. Venakides and X. Zhou,
\reftitle{On the long-time limit of semiclassical (zero dispersion limit) solutions of the focusing nonlinear Schr\"odinger equation: Pure radiation case},
\href{ https://doi-org.gate.lib.buffalo.edu/10.1002/cpa.20142}
{Commun. Pure Appl. Math. \textbf{59}, 1379--1432 (2006)}

\bibitem{TV2012}
A. Tovbis and S. Venakides,
\reftitle{Semiclassical limit of the scattering transform for the focusing nonlinear Schr\"odinger equation},
\href{https://doi.org/10.1093/imrn/rnr092}
{Int. Math. Res. Notices \textbf{2012}, 2212--2271 (2012)}

\bibitem{Whitham}
G.B. Whitham,
\reftitle{Linear and nonlinear waves},
Wiley, 1974

\bibitem{ZabuskyKruskal} 
N.J. Zabusky and M.D. Kruskal, 
\reftitle{Interaction of solitons in a collisionless plasma and the recurrence of initial states}, 
\href{https://doi.org/10.1103/PhysRevLett.15.240}
{Phys. Rev. Lett. {\bf15}, 240 (1965)}

\bibitem{zakharov2009}
V.E. Zakharov,
\reftitle{Turbulence in integrable systems},
\href{\doibase 10.1111/j.1467-9590.2009.00430.x}
{{{Stud. Appl. Math.}} \textbf{122},\ {219--234} (2009)}

\bibitem{ZS1972}
V.E. Zakharov and A.B. Shabat,
\reftitle{Exact theory of two-dimensional self-focusing and one-dimensional self-modulation of waves in nonlinear media},
Sov. Phys. JETP \textbf{34}, 62 (1972) 	
\end{thebibliography}
\end{document}